\documentclass[11pt]{article}
\usepackage{mathrsfs}
\usepackage{amsfonts}
\usepackage{amssymb,amsmath,amsthm}
\usepackage{graphicx}
\usepackage{subfigure}
\usepackage{color}
\usepackage{float}

\def\disp{\displaystyle}

\def\dref#1{(\ref{#1})}
\def\crr{\cr\noalign{\vskip2mm}}

\newcommand{\Ial}{{}_0I_t^\alpha}

\numberwithin{equation}{section}

\newtheorem{theorem}{Theorem}[section]
\newtheorem{definition}{Definition}[section]

\newtheorem{lemma}{Lemma}[section]
\newtheorem{remark}{Remark}[section]

\renewcommand{\hline}{{\mathbb H}}

\newcommand{\dd}  {{\rm d}\hbox{\hskip 0.5pt}}

\newcommand{\mm}     {{\hbox{\hskip 0.5pt}}}

\newcommand{\bluff}  {{\hbox{\raise 15pt \hbox{\mm}}}}
\newcommand{\sbluff} {{\hbox{\raise  9pt \hbox{\mm}}}}
\begin{document}
\renewcommand{\thefootnote}{\fnsymbol{footnote}}
\newcommand{\footremember}[2]{%
   \footnote{#2}
    \newcounter{#1}
    \setcounter{#1}{\value{footnote}}%
}
\newcommand{\footrecall}[1]{%
    \footnotemark[\value{#1}]%
}
\makeatletter
\def\blfootnote{\gdef\@thefnmark{}\@footnotetext}
\makeatother

\begin{center}
{\LARGE \bf A weakly singular Volterra integral equation and  \\ [0.4ex]
the fractional Lyapunov asymptotic stability \\[0.5ex]
theorem  
}\\[4ex]
Zuoshunhua Shi,  \ \ \ Hua-Cheng Zhou\;
\blfootnote{This work is supported by the National Natural Science Foundation of China 
(no. 62173348, 12161141013) and Hunan Basic Science Research Center for Mathematical
Analysis (2024JC2002).}
\blfootnote{Z.S.H. Shi (shizsh@csu.edu.cn) and H.-C. Zhou (hczhou@amss.ac.cn)
 are with School of Mathematics and Statistics, HNP-LAMA,
  Central South University, Changsha, 410083,  PR China}
\end{center}
\vspace{3mm}

\noindent {\bf Abstract:} 
We study the weakly singular Volterra integral equation
$$
V(t)=A-\frac{1}{\Gamma(\alpha)}\int_0^t(t-s)^{\alpha-1}h(|V(s)|)\,\mathrm ds,\qquad t\geq0,
$$
where $\alpha\in(0,1)$, $A>0$ and $h$ is a class-$\mathcal{K}$ function.
This equation arises as the scalar reduction of fractional-order Lyapunov inequalities
through the fractional variation-of-constants formula.
We prove that it admits a unique global solution $V$ with $0<V(t)\le A$ for all $t\ge0$,
that $V$ is non-increasing, and that $V(t)\to0$ as $t\to\infty$,
under the pure class-$\mathcal{K}$ hypothesis on $h$.
If, in addition, $h$ is locally H\"older continuous on $(0,\infty)$
with exponent $\beta>1-\alpha$, then the rescaled function
$W(t):=V(t^{1/\alpha})$ is of class $C^1$ on $[0,\infty)$ and satisfies $W'(t)\le0$ for all $t\ge0$.
We determine the optimal regularity of
the solution---it is $\alpha$-H\"older continuous, with the sharp expansion
$$
V(t)=A-\frac{h(A)}{\Gamma(\alpha+1)}t^{\alpha}+o(t^{\alpha})
$$
at the origin---and its sharp decay rate,
$h(V(t))\sim \frac{A}{\Gamma(1-\alpha)}t^{-\alpha}$ as $t\to\infty$.
In the power case $h(x)\sim\kappa x^{p}$ one has
$V(t)\sim(\frac{A}{\kappa\,\Gamma(1-\alpha)})^{1/p}t^{-\alpha/p}$.
As applications,  Theorem \ref{thm-Gen-fl} gives a complete proof
of the fractional-order Lyapunov asymptotic stability theorem
of Li, Chen and Podlubny (Theorem 11, {\it Automatica}, 45(2009), 1965--1969)
under the general class-$\mathcal K$ hypothesis, a case not covered by the linear, power-law and convex arguments available in the literature.
This theorem is to fractional-order stability theory what Lyapunov's classical theorem
is to the stability theory of differential equations, and it ranks among the basic tools of fractional-order control.
Its proof closes a gap that has been open since the theorem first appeared.
We further show that the fractional inequality
${}^C_0D_t^{\alpha}V(t)\le -cV^{p}(t)$ with $0<p<1$ implies Mittag-Leffler stability of $V$,
yet neither entails nor precludes finite-time decay to zero,
in contrast to the integer-order criterion $\dot V(t)\le -cV^{p}(t)$,
which entails finite-time decay to zero. 
 \vspace{3mm}

\vspace{3mm}


\section{Introduction}

Fractional-order differential equations provide a natural framework for
dynamical systems with memory and hereditary effects, and Caputo-type initial
value problems have become a standard model in viscoelasticity, anomalous
diffusion and other applied sciences \cite{Podlubny1999,KilbasJbook2006,Diethelm2010}.
In parallel, a rigorous well-posedness theory for nonlinear time-fractional
PDEs, including Hamilton--Jacobi equations with Caputo's time-fractional
derivative and the time-fractional diffusion equation, has been developed
through the viscosity approach \cite{GigaNamba2017CPDE,GigaMitakeSato2022JDE}.
As in the classical theory, stability is a central question. For integer-order
systems, Lyapunov's direct method reduces the stability analysis to the
construction of a scalar function whose derivative along the flow is
non-positive \cite{Lyapunov1892,YoshizawaBook}. The extension to the fractional
setting is far from immediate. The Caputo derivative satisfies no classical
Leibniz product rule, and the sign of the Caputo derivative of a continuous
function does not, in general, determine its monotonicity: a function may
satisfy ${}^C_0D_t^{\alpha}V(t)\leq0$ everywhere and still fail to be
non-increasing, as the example in Remark \ref{RemarkCaputoDdoesnotimplymono}
shows.

A general Lyapunov criterion for Caputo fractional-order systems was proposed
by Li, Chen and Podlubny in \cite{LiYanCP2009auto}, together with the notion of
Mittag-Leffler stability. For the Caputo fractional-order system
\dref{sys-plantFrac}, their Theorem 11 asserts that if there exist
class-$\mathcal{K}$ functions $\alpha_1,\alpha_2,\alpha_3$ and a Lyapunov
function $V(t,x)$ such that
$$
\mbox{(A)}\qquad
\alpha_1(\|x\|)\leq V(t,x)\leq\alpha_2(\|x\|),\qquad
{}^C_0D_t^{\alpha}V(t,x)\leq-\alpha_3(\|x\|)
$$
hold along every solution of \dref{sys-plantFrac}, then the equilibrium point
$x=0$ is asymptotically stable. The criterion has become one of the basic tools
of fractional-order stability theory. It is routinely invoked in stability
proofs for fractional observers and controllers, in synchronization and consensus
schemes for fractional networks, and in the body of literature built on the
pattern ``Lyapunov function plus fractional derivative inequality''
\cite{LiYanCP2009auto,AguilaCamacho2014cns}.

It is useful to distinguish what Theorem 11 asserts from what its proof
delivers. When the class-$\mathcal{K}$ functions are linear, $\alpha_i(s)=k_is$,
the conclusion follows from the linear fractional comparison argument and the
proof given in \cite{LiYanCP2009auto} is complete. In the genuinely nonlinear
case, the arguments available in the literature, including the one given in
\cite{LiYanCP2009auto}, establish at most
$\liminf_{t\to\infty}\|x(t)\|=0$. The full convergence
$\lim_{t\to\infty}\|x(t)\|=0$ required by asymptotic stability is not obtained
there.

The obstruction is structural. Condition (A) constrains the Caputo derivative
of the scalar function $t\mapsto V(t,x(t))$, not its values, and the standard
device that converts one into the other is the fundamental identity
$$
V(t,x(t))=V(t_0,x(t_0))
+{}_0I_t^{\alpha}\big[{}^C_0D_\tau^{\alpha}V(\tau,x(\tau))\big],
\qquad t\geq t_0,
$$
which turns the differential inequality into an inequality of weakly singular
Volterra type. The object to be analysed is therefore not an ordinary
differential equation but an integral equation whose kernel $(t-s)^{\alpha-1}$
is singular at the upper limit. Two features of this reformulation matter. The
singularity of the kernel forces a cusp at the origin, so that the solution is
merely $\alpha$-H\"older continuous there and its classical derivative is
unbounded. Furthermore, the current value $V(t)$ is tied to the entire history of
$h(V(\cdot))$ through a weight that blows up at the upper limit, so that decay
has to be extracted from the way the singularity accumulates the past, rather
than from a pointwise condition on the derivative.

These observations lead us to isolate the scalar problem. Setting
$A:=V(0,x_0)$ and $h(\cdot):=\alpha_3(\alpha_2^{-1}(\cdot))$, condition (A)
yields the scalar differential inequality
\begin{equation}\label{eq:gen-DV-hV}
{}^C_0D_t^{\alpha}V(t)\leq -h(V(t)),\qquad
h(\cdot):=\alpha_3(\alpha_2^{-1}(\cdot)),
\end{equation}
with a class-$\mathcal{K}$ function $h$, and hence, through the fractional
variation-of-constants formula, to the global analysis of the model weakly
singular Volterra integral equation
\begin{equation}\label{eq:main-VI}
V(t)=A-\frac{1}{\Gamma(\alpha)}\int_0^t(t-s)^{\alpha-1}h(|V(s)|)\,\mathrm ds,
\qquad A>0.
\end{equation}
By the lower bound in (A) one has $A\geq\alpha_1(\|x_0\|)>0$ whenever
$x_0\neq0$, and the case $x_0=0$ is trivial. The modulus inside $h$ is retained
because the positivity of $V$ is not known in advance, since it is one of the
conclusions. The inequality \eqref{eq:gen-DV-hV} is then reduced to
\eqref{eq:main-VI} by comparison, see Theorem \ref{thm-CapDVleq-hV-hinf} and
Theorem \ref{thm-DV-hV-gencase}. The core of this paper is a systematic study of
\eqref{eq:main-VI}: global existence, uniqueness, monotonicity, regularity and
decay of its solution, under two sets of hypotheses on $h$, a local Lipschitz
condition on $h'$ (Section \ref{Sec-h-diff-property}) and the pure
class-$\mathcal{K}$ assumption with no differentiability hypothesis whatsoever
(Section \ref{Sec-h-classK-property}).

Several prior works address parts of this program under stronger or more
specific hypotheses, and it is useful to delimit their scope. Feng, Li, Liu and
Xu \cite{FengLiLiuXu2018} consider the one-dimensional autonomous equation
${}^C D_c^{\gamma}u=f(u)$ with $f\in C^1$ and $f'$ locally Lipschitz, prove that
the solution is monotone on its maximal interval of existence (Theorem 3.4),
obtain comparison principles, and analyse in detail the power case $f(u)=Au^p$
including blow-up. Their Lemma 3.1 also yields $u\in C^1(0,T_b)$ and
$u'(t)=\mathcal O(t^{\gamma-1})$, which is the scalar autonomous counterpart of
the regularity established here. Their analysis is restricted to autonomous
scalar equations, so a general class-$\mathcal{K}$ nonlinearity, the global
solution of \eqref{eq:main-VI} and the differentiability of the rescaled
function $W(t)=V(t^{1/\alpha})$ at $t=0$ are not covered there. Vergara and
Zacher \cite{VergaraZacher2015} treat, among other nonlocal subdiffusion
problems, the scalar equation ${}^C_0D_t^{\alpha}(u-u_0)+\nu u^{\gamma}=0$ with
$\nu>0$, $\gamma>0$ and $u_0>0$, and construct sub- and supersolutions which
give the two-sided bounds
$c_1/(1+t^{\alpha/\gamma})\leq u(t)\leq c_2/(1+t^{\alpha/\gamma})$
(Theorem 7.1). For a power nonlinearity and a non-negative solution the rate
$t^{-\alpha/p}$ with $p=\gamma$ is thus already available, including the range
$0<p<1$. Their comparison argument applies to the equation rather than to a
differential inequality, and it does not cover a general class-$\mathcal{K}$
law. Tuan and Trinh \cite{TuanTrinh2018} prove, for
$u\in x_0+I_{0+}^{\alpha}C([0,T];\mathbb R^d)$ and a convex differentiable $V$,
the fractional inequality
${}^C D_{0+}^{\alpha}V(u(t))\leq\langle\nabla V(u(t)),{}^C D_{0+}^{\alpha}u(t)\rangle$,
without requiring the pseudo-state $u$ to be differentiable, an assumption
present in earlier proofs. Their stability theorem is then obtained under
polynomial bounds $C_1\|x\|^a\leq V(x)\leq C_2\|x\|^b$ and
$\langle\nabla V(x),f(x)\rangle\leq-C_3\|x\|^c$ with $c\geq b$. These bounds are
of power type, so the criterion does not extend to a general class-$\mathcal{K}$
triple $(\alpha_1,\alpha_2,\alpha_3)$, although the differentiability issue in
the convex case is settled there. Rajchakit \cite{Rajchakit2026} derives two-term
asymptotic expansions for nonlinear Caputo relaxation. Under a terminal
uniform-integrability condition and for an orbit already known to be positive,
non-increasing and convergent to zero, the leading term in \S5.1, (26)--(30), of
\cite{Rajchakit2026} is
$\phi(u(t))\sim\frac{u_0}{\Gamma(1-\alpha)}t^{-\alpha}$, with the same leading
constant as in the decay estimate obtained below. The hypotheses on the orbit
and on the restoring law are stronger than the ones used here, so the present
results do not follow from \cite{Rajchakit2026}, but the overlap in the leading
asymptotic term should be recorded.

The contributions of this paper can now be stated. We give a complete proof of
Theorem 11 for general class-$\mathcal{K}$ functions, without requiring the
solution $x(\cdot)$ of \dref{sys-plantFrac} to be differentiable, which covers
the nonlinear case left open in \cite{LiYanCP2009auto} and in the earlier
literature, the linear case being contained in \cite{LiYanCP2009auto} and the
power-law and convex settings having been treated in \cite{VergaraZacher2015}
and \cite{TuanTrinh2018} respectively. Compared with
\cite{FengLiLiuXu2018,VergaraZacher2015}, neither the positivity of the
solution nor its decay is assumed in advance, and $h$ is an arbitrary
class-$\mathcal{K}$ function rather than a power. We then establish, for
\eqref{eq:main-VI} and under two hypotheses on $h$, existence, uniqueness,
positivity, monotonicity, sharp regularity and sharp decay with explicit
constants. Finally, we show that the differential inequality
\eqref{eq:gen-DV-hV}, with a general class-$\mathcal{K}$ $h$, already forces
$V(t)\to0$ for every non-negative solution, and we determine the range of
validity of the frequently used inequality
${}^C_0D_t^{\alpha}V(t)\leq-cV^p(t)$ with $0<p<1$.

In Section \ref{Sec-h-diff-property} we assume that $h$ is a
class-$\mathcal{K}$ function whose derivative $h'$ is locally Lipschitz
continuous on $(0,+\infty)$. This covers the model cases $h(x)=x^p$ with
$p\in(0,1)$ and $h(x)=x^2\sin\frac1x+4x+x^2$, for which $h'(0)$ need not exist.
Theorem \ref{Thm-Dhinfatzero-mono-asympt} shows that every solution of
\eqref{eq:main-VI} in $C[0,+\infty)$ is strictly positive, that the rescaled
function $W(t):=V(t^{1/\alpha})$ belongs to $C^1[0,+\infty)$ with $W'(t)\leq0$,
and that $V(t)\to0$ as $t\to\infty$. Remark \ref{rem-exist-uniq-Thm41} records,
as a by-product of the proof, that \eqref{eq:main-VI} then has exactly one
solution in $C[0,+\infty)$. Theorem \ref{thm-DV-hV-gencase} shows in turn that
the differential inequality \eqref{eq:gen-DV-hV}, with a pure
class-$\mathcal{K}$ $h$, already forces $V(t)\to0$ for every non-negative
solution. Combining these results,
Theorem \ref{thm-Gen-fl} proves Theorem 11 under the general
class-$\mathcal{K}$ hypothesis:
under the sector bounds and the fractional Lyapunov inequality in (A), the
equilibrium of \dref{sys-plantFrac} is asymptotically stable. A second,
independent proof of the attractivity step, which does not use the properties
of the solutions of \eqref{eq:main-VI}, is given in
Remark \ref{rem-alt-proof-Thm34}. A version for
convex Lyapunov functions is given in Theorem \ref{tm-Vconvex-Lya-zero}, on the
basis of the fractional inequality for convex functions of the solution. The
local well-posedness and the regularity of the solutions of
\dref{sys-plantFrac} itself are established in Theorem \ref{thm-gen-localsolution},
together with an example showing that uniqueness may fail when the
Lipschitz-H\"older condition \eqref{eq:fun-f-sat-Lip} is dropped.

In Section \ref{Sec-h-classK-property} we work under the pure
class-$\mathcal{K}$ assumption only. By smoothing $h$ through mollification and
passing to the limit by uniform estimates and compactness, we prove that
\eqref{eq:main-VI} has a global solution
$V\in C^{0,\alpha}_{\rm loc}[0,+\infty)$ satisfying $0<V(t)\leq A$,
non-increasing, with $V(t)\to0$ as $t\to\infty$
(Theorem \ref{thm:global-classK}). A comparison principle for positive
solutions (Lemma \ref{lem:comparison}) yields uniqueness among all continuous
solutions (Theorem \ref{thm:all-solutions}). Theorem \ref{thm:regularity} gives
the optimal regularity, namely $\alpha$-H\"older continuity together with the
sharp expansion
$$
V(t)=A-\frac{h(A)}{\Gamma(\alpha+1)}t^{\alpha}+o(t^{\alpha})
$$
at the origin, so that $V'(0^+)=-\infty$ and no exponent larger than $\alpha$
is possible. The fractional cusp at the origin is intrinsic to Caputo equations
with non-zero initial values and is removed by the rescaling
$W(s)=V(s^{1/\alpha})$, for which $W'(0)=-h(A)/\Gamma(\alpha+1)$. The same
theorem shows that ${}^C_0D_t^{\alpha}V$ exists and is continuous on
$[0,+\infty)$ with ${}^C_0D_t^{\alpha}V(t)=-h(V(t))$. If in addition $h$ is
locally H\"older continuous on $(0,\infty)$ with exponent $\beta>1-\alpha$,
then Theorem \ref{thm:regularity-diffable} gives $V\in C^1(0,\infty)$ together
with an explicit formula for $V'$, as well as $W\in C^1[0,\infty)$ with
$W'(t)\leq0$. Theorem \ref{thm:decay-classK} identifies the
decay rate, $h(V(t))\sim\frac{A}{\Gamma(1-\alpha)}t^{-\alpha}$ as
$t\to\infty$, then
$V(t)\sim(\frac{A}{\kappa\,\Gamma(1-\alpha)})^{1/p}t^{-\alpha/p}$ in the
power case $h(x)\sim\kappa x^p$, and two-sided estimates when $h$ is sandwiched
between two power functions. As a by-product, we clarify the frequently used
inequality ${}^C_0D_t^{\alpha}V(t)\leq-cV^p(t)$ with $0<p<1$: such a $V$ is
Mittag-Leffler stable, yet the inequality neither entails nor precludes
finite-time decay to zero, in contrast with the integer-order analogue
$\dot V(t)\leq-cV^p(t)$, which entails it
(Theorem \ref{thm:CapD-Vine-p-stability}).

The proofs are self-contained and combine several tools: the Schauder (or
Leray--Schauder) fixed point theorem for local existence, fractional comparison
principles, a regularization procedure with estimates uniform in the smoothing
parameter, and the Laplace transform together with the Hardy--Littlewood
Tauberian theorem for the sharp decay rates.

The paper is organized as follows. Section \ref{sec:pre} collects the necessary
background on fractional calculus and the auxiliary results used throughout.
Section \ref{Sec-h-diff-property} contains the main results under the local
Lipschitz assumption on $h'$, in particular the proof of Theorem 11 of
\cite{LiYanCP2009auto}. Section \ref{Sec-h-classK-property} establishes
existence, uniqueness, regularity and decay of the solution under the pure
class-$\mathcal{K}$ assumption.

\section{Preliminaries}
  \label{sec:pre}
Let $\mathbb{R}$ be an Euclidean space.
For $T\in(0,\infty)$,  we define
$C[0,T]=\{\varphi:[0,T]\to\mathbb{R}|\ \varphi(\cdot)$ is continuous$\}$ and
$C[0,\infty)=\{\varphi:[0,\infty)\to\mathbb{R}|\ \varphi(\cdot)$ is continuous$\}$.
Let $\alpha\in(0,1)$ and $x\in C[0,+\infty)$. ${}_0I_t^\alpha$ denotes the Riemann--Liouville fractional integral
\begin{equation}\label{eq:Ia}
(\Ial x)(t):=\frac{1}{\Gamma(\alpha)}\int_0^t(t-s)^{\alpha-1}x(s)\,\mathrm ds.
\end{equation}
The Caputo fractional derivative is defined by
$$
{}^C_0D_t^{\alpha}x(t):=\frac{\mathrm d}{\mathrm dt}({}_0I_t^{1-\alpha}(x-x(0)))(t)=\frac{1}{\Gamma(1-\alpha)}
\bigg[\frac{\mathrm d}{\mathrm dt}\int_0^t(t-s)^{-\alpha}x(s)
{\rm d}s-t^{-\alpha}x(0)\bigg].
$$
Throughout this paper, when a Caputo derivative is evaluated at $t=0$
(as in differential inequalities stated on closed intervals), it is
understood as the continuous extension
$\lim_{t\to0^+}{}^C_0D_t^{\alpha}x(\cdot)$, whenever this limit exists.
%
We consider the Caputo fractional-order system
\begin{equation}\label{sys-plantFrac}
\left\{\begin{array}{l}
{}^C_0D_t^{\alpha}x(t)=f(t,x(t)),\ t>0,\crr
x(0)=x_0
\end{array}\right.
\end{equation}
where $x(t)\in\mathbb{R}^n$ is the system state, $x_0$ is the initial state,
$f:[0,+\infty)\times\Omega\to\mathbb{R}^n$ is continuous on
 $[0,+\infty)\times\Omega$, and $\Omega\subset\mathbb{R}^n$ is a domain containing the origin $x=0$.
 By \cite[Page 323, Theorem 5.15]{KilbasJbook2006},
the Caputo initial value problem \dref{sys-plantFrac} is equivalent,
for continuous $x$, to a second kind Volterra integral equation
\begin{equation}\label{Thm11-IE}
x(t)=x_0+\frac1{\Gamma(\alpha)}\int_0^t(t-s)^{\alpha-1}f(s,x(s))\mathrm{d}s,\qquad t\in[0,T],
\end{equation}
in the sense that $x\in C[0,T]$ solves \dref{Thm11-IE} if and only if $x\in C[0,T]$, ${}^C_0D_t^{\alpha}x=f(\cdot,x(\cdot))$ on $(0,T]$, and $x(0)=x_0$.

\begin{definition}
A continuous function $h:[0, +\infty) \to[0,+\infty)$ is said to
belong to class-$\mathcal{K}$ if it is strictly increasing and $h(0)=0$.
\end{definition}

\begin{definition}\cite[page 56]{MittagLefflerbook}
\sl   The Mittag-Leffler function for an arbitrary
$z\in\mathbb{C}$
is defined as
\begin{equation}\label{ML-d}
E_{\alpha,\beta}(z)=\sum^{\infty}_{k=0}\frac{z^{k}}{\Gamma(\alpha k+\beta)},\ \ \alpha, \ \beta>0.
\end{equation}
In particular, if $\beta=1$, then $E_{\alpha,1}(z)=E_{\alpha}(z)$. If $\alpha=\beta=1$, then $E_{1,1}(z)=e^{z}$.
\end{definition}

\begin{lemma}\label{lem-1d-ML-est}\cite[Lemma 2.6]{2025HuangZhouAMO}
Let $\alpha\in(0,1)$.
The Mittag-Leffler function $E_{\alpha,\alpha}(-s^{\alpha})$ has the following estimate
\begin{equation}\label{Eaa-1}
\frac{C_1}{1+s^{2\alpha}}\leq
E_{\alpha,\alpha}(-s^{\alpha})\leq
\frac{C_2}{1+s^{2\alpha}},\ \ \forall\ s\geq0
\end{equation}
for some positive constants $0<C_1\leq C_2$.
Moreover, $\lim_{s\to\infty}\frac{E_{\alpha,\alpha}(-s)}{(1+s^{2})^{-1}}=-\frac{1}{\Gamma(-\alpha)}$.
\end{lemma}
The Mittag-Leffler function $E_{\alpha}(-\omega)$
satisfies the following estimate
 \cite[Theorem 4]{MLF2014OPEST}:
\begin{equation}\label{4t1-6-ne}
\frac{1}{1+\Gamma(1-\alpha)\omega}\leq E_{\alpha}(-\omega)\leq\frac{1}{1+\Gamma(1+\alpha)^{-1}\omega}, \ \ \forall\ \omega\geq0.
\end{equation}

We shall use the following Leray-Schauder fixed point theorem.
\begin{lemma}\label{Leray-Schauder-fixed-point-lem}
\cite[Page 286, Theorem 11.6]{David-Trudinger-book}
Let $S$  be a compact mapping of a Banach space $X$ into itself, and suppose there exists a constant $M$
such that $\|x\|\leq M$ for all $x\in X$ and $\sigma\in[0,1]$ satisfying $x=\sigma Sx$.
Then $S$  has a fixed point.
\end{lemma}



\begin{lemma}\label{lem-sequ-estim-I0}
Let $T>0$ and let $\{f_n\}_{n=0}^{\infty}\in C[0,T]$ be a sequence of
nonnegative 
functions satisfying
$$
f_{n+1}(t)\leq A\int_0^tf_n(s)ds, \;\;
t\in[0,T], \;n\geq0,
$$
where $A$ is a positive constant. Then
$f_n$ has the following estimation
\begin{equation}\label{lem-con-ft-estm}
|f_{n}(t)|\leq\|f_0\|_{C[0,T]}\frac{(At)^n}{n!},\;
t\in[0,T],\; n=1,2...
\end{equation}
Moreover, $\lim_{n\to\infty}\|f_n\|_{C[0,T]}=0$, and for any $\rho>0$, $\{f_n\}_{n=0}^{\infty}$ satisfies
\begin{equation}\label{lem-con-ft-estm-rho}
\sum\limits_{n=0}^{\infty}
\|f_n\|_{C[0,T]}^{\rho}
\leq \|f_0\|^{\rho}_{C[0,T]}\sum_{k=0}^{\infty}\frac{(AT)^{\rho k}}{(k!)^{\rho}}
<+\infty.
\end{equation}
\end{lemma}
\noindent{\bf Proof.}
First, by induction,  \dref{lem-con-ft-estm} can be directly derived.
From \dref{lem-con-ft-estm}, we have
$\|f_n\|_{C[0,T]}\leq\|f_0\|_{C[0,T]}\frac{(AT)^n}{n!}\to0$ as $n\to\infty$.
For any $\rho>0$,
$
\sum_{n=0}^{\infty}\|f_n\|^{\rho}_{C[0,T]}
\leq \|f_0\|^{\rho}_{C[0,T]}\sum_{n=0}^{\infty}\frac{(AT)^{\rho n}}{(n!)^{\rho}}.
$
The series on the right converges for every $\rho>0$ by the ratio test,
since $\frac{(AT)^{\rho(n+1)}}{((n+1)!)^{\rho}}\big/\frac{(AT)^{\rho n}}{(n!)^{\rho}}
=\left(\frac{AT}{n+1}\right)^{\rho}\to0$ as $n\to\infty$.
Hence \dref{lem-con-ft-estm-rho} follows.
\hfill$\Box$

\begin{lemma}\label{lem-sequ-estim}
Let $T>0$ and $B>0$. Let $\{A_n\}_{n=0}^{\infty}$ be a nonnegative sequence.
Suppose that $\{f_n\}_{n=0}^{\infty}\in C[0,T]$ is a sequence of nonnegative
continuous functions satisfying
$$
f_{n+1}(t)\leq A_n+B\int_0^tf_n(s)\mathrm{d}s, \;\;
t\in[0,T], \;n\geq0.
$$
Then the following statements are true:\\
(i) If $A_n=A$ for all $n\geq0$, then\vspace{-2mm}
\begin{equation}\label{Lm-Inest-lem-f}
\|f_{n+1}\|_{C[0,T]}\leq (A+\|f_0\|_{C[0,T]})e^{BT}.
\end{equation}
(ii) Let $\rho>0$. If  $\{A_n\}_{n=0}^{\infty}$ satisfies
$\sum_{n=0}^{\infty}A_n<+\infty$ and
$\sum_{n=0}^{\infty}A_n^{\rho}<+\infty$, then
\begin{equation}\label{Lm-Inest-lem-fI}
\sum_{n=0}^{\infty}\|f_{n}\|^{\rho}_{C[0,T]}
\leq
C_\rho\sum_{k=0}^{\infty}\frac{(BT)^{k\rho}}{(k!)^{\rho}}<+\infty,
\end{equation}
where $C_\rho=\max\Big\{
2^{\rho-1}\Big[\big(\sum_{j=0}^{\infty} A_j \big)^{\rho}
+2\|f_0\|_{C[0,T]}^{\rho}\Big],
\Big(2\|f_0\|_{C[0,T]}^{\rho}+\sum_{j=0}^{\infty}A_{j}^{\rho}\Big)\Big\}$.
\end{lemma}
{\bf Proof.}
First, by induction, it is straightforward to verify that  for each $n$, 
\begin{equation}\label{Lm-Inest-pr1}
f_{n+1}(t)\leq \sum_{k=0}^n A_{n-k}\frac{(Bt)^k}{k!}
+\|f_0\|_{C[0,T]}\frac{(Bt)^{n+1}}{(n+1)!}, \ \ \ \forall t\in [0,T].
\end{equation}
Hence, taking suprema over $[0,T]$ in \dref{Lm-Inest-pr1}
 and writing $w_k:=\frac{(BT)^k}{k!}$ gives
\begin{equation}\label{short-L2-pr}
\|f_{n+1}\|_{C[0,T]}\leq \sum_{k=0}^{n}A_{n-k}w_k
+\|f_0\|_{C[0,T]}w_{n+1},\qquad n\geq0.
\end{equation}
Since $e^{BT}=\sum_{k=0}^{\infty}\frac{(BT)^k}{k!}$,
it follows from $A_n=A$ for all $n\geq0$ and \dref{short-L2-pr} that \dref{Lm-Inest-lem-f} holds.

For (ii), let $S_A:=\sum_{j=0}^{\infty}A_j$.
For $\rho\geq1$, then the inequality $(a+b)^{\rho}\leq2^{\rho-1}(a^{\rho}+b^{\rho})$ and
H\"older's inequality give
$$
\begin{aligned}
\|f_{n+1}\|^{\rho}_{C[0,T]}
&\leq 2^{\rho-1}\Big[\Big(\sum_{k=0}^{n}A_{n-k}w_k\Big)^{\rho}
+\|f_0\|^{\rho}_{C[0,T]}w_{n+1}^{\rho}\Big]\\
&\leq 2^{\rho-1}\Big[S_A^{\rho-1}\sum_{k=0}^{n}A_{n-k}w_k^{\rho}
+\|f_0\|^{\rho}_{C[0,T]}w_{n+1}^{\rho}\Big],
\end{aligned}
$$
This, together with the identity 
obtained by interchanging the order of summation
$\sum_{n=0}^{\infty}\sum_{k=0}^{n}$ $A_{n-k}w_k^{\rho}
=\sum_{k=0}^{\infty}\sum_{n=k}^{\infty}A_{n-k}w_k^{\rho}
=\sum_{k=0}^{\infty}\big(w_k^{\rho}\sum_{n=k}^{\infty}A_{n-k}\big)
=S_A\sum_{k=0}^{\infty}w_k^{\rho}$, yields
\begin{equation}\label{sum-est-rho-geq1}
\begin{array}{l}\disp
\sum_{n=0}^{\infty}\|f_{n+1}\|^{\rho}_{C[0,T]}
\leq2^{\rho-1}\Big[S_A^{\rho-1}\sum_{n=0}^{\infty}
\sum_{k=0}^{n}A_{n-k}w_k^{\rho}
+\|f_0\|_{C[0,T]}^{\rho}
\sum_{n=0}^{\infty}w_{n+1}^{\rho}\Big]\crr\disp\hspace{0cm}
\leq
2^{\rho-1}\Big[S_A^{\rho}
+\|f_0\|_{C[0,T]}^{\rho}\Big]\sum_{k=0}^{\infty}w_k^{\rho}.
\end{array}
\end{equation}
For $0<\rho\leq1$, the inequality
$(a_1+\cdot +a_k)^{\rho}\leq a_1^\rho+\cdot+a_k^\rho$
gives
$
\|f_{n+1}\|^{\rho}_{C[0,T]}
\leq \sum_{k=0}^{n}A_{n-k}^{\rho}w_k^{\rho}
+\|f_0\|^{\rho}_{C[0,T]}w_{n+1}^{\rho}.
$
Similarly, we can obtain
$\sum_{n=0}^{\infty}\|f_{n+1}\|^{\rho}_{C[0,T]}
\leq(\sum_{j=0}^{\infty}A_j^{\rho}+\|f_0\|^{\rho}_{C[0,T]})\sum_{k=0}^{\infty}w_k^{\rho}$.
Since $w_k^{\rho}=\frac{(BT)^{k\rho}}{(k!)^{\rho}}$ and
$\sum_{k=0}^{\infty}w_k^{\rho}<+\infty$ by the ratio test
($w_{k+1}^{\rho}/w_k^{\rho}=\big(\frac{BT}{k+1}\big)^{\rho}\to0$), and since
$\sum_{k=0}^{\infty}w_k^{\rho}\geq1$,
we know that \dref{Lm-Inest-lem-fI} holds.
\hfill$\Box$

\begin{lemma}\label{Ne-operator-com-lem}
Let $\alpha\in(0,1)$ and $f\in C(\mathbb{R}^2)$.
Define an operator $\mathcal{F}:C[0,T]\to C[0,T]$ as
$$
\mathcal{F}(h)(t)=t^{\frac{1}{\alpha}-1}\int_0^t(t^{\frac{1}{\alpha}}-s^{\frac{1}{\alpha}})^{\alpha-1}f(s,h(s))
\mathrm{d}s
$$
Then the operator $\mathcal{F}$ is a compact operator on $C[0,T]$.
\end{lemma}
\noindent{\bf Proof.}
Let $M\subset C[0,T]$ be a bounded set and $C_1$ be a constant such that
$|h(t)|\leq C_1$ for all $t\in[0,T], h\in M$.
Since $f\in C(\mathbb{R}^2)$, a constant $C_2$ exists such that
$C_2=\max\{|f(t,s)|:(t,s)\in[0,T]\times[-C_1,C_1]\}$.
By the substitution $s=tu$, $\mathcal{F}(h)(t)$ can be written as
$\mathcal{F}(h)(t)=t\int_0^1(1-u^{1/\alpha})^{\alpha-1}f(tu,h(tu))\mathrm{d}u$.
Since
$\int_0^1(1-u^{1/\alpha})^{\alpha-1}\mathrm{d}u=\alpha B(\alpha,\alpha)
=\alpha\Gamma^2(\alpha)/\Gamma(2\alpha)$, we have
$|\mathcal{F}(h)(t)|
\leq C_2t\frac{\alpha\Gamma^2(\alpha)}{\Gamma(2\alpha)}\leq C_2T\frac{\alpha\Gamma^2(\alpha)}{\Gamma(2\alpha)}$, for all
$h\in M$, which yields the uniform boundedness of $\mathcal{F}(M)$.

Next, we show $\mathcal{F}(M)$ is equicontinuous.
Fix $0\le t_1\leq t_2\leq T$.
If $t_2\geq 2t_1>0$, then by
$\mathcal{F}(h)(t)$ $=t\int_0^1(1-u^{1/\alpha})^{\alpha-1}f(tu,h(tu))\mathrm{d}u$
and $(t_1+t_2)/(t_2-t_1)=(t_2/t_1+1)(t_2/t_1-1)^{-1}\le 3$, we obtain
\begin{equation}
|\mathcal{F}(h)(t_1)-\mathcal{F}(h)(t_2)|
\le C_2\frac{\alpha\Gamma^2(\alpha)}{\Gamma(2\alpha)}(t_1+t_2)
\le 3C_2\frac{\alpha\Gamma^2(\alpha)}{\Gamma(2\alpha)}(t_2-t_1).
\end{equation}
It remains to consider $0<t_1<t_2<2t_1$. With $k(t,s):=(t^{1/\alpha}-s^{1/\alpha})^{\alpha-1}$,
we write
$$
\begin{aligned}
\mathcal{F}(h)(t_1)-\mathcal{F}(h)(t_2)
&=(t_1^{\frac{1}{\alpha}-1}-t_2^{\frac{1}{\alpha}-1})\int_0^{t_1}k(t_1,s)f(s,h(s))\mathrm{d}s\\
&\quad+t_2^{\frac{1}{\alpha}-1}\int_0^{t_1}\big[k(t_1,s)-k(t_2,s)\big]f(s,h(s))\mathrm{d}s\\
&\quad-t_2^{\frac{1}{\alpha}-1}\int_{t_1}^{t_2}k(t_2,s)f(s,h(s))\mathrm{d}s=:I_1+I_2+I_3.
\end{aligned}
$$
{\em Estimate the term $I_1$.}
Since $\int_0^{t_1}k(t_1,s)\mathrm{d}s=\alpha\Gamma^2(\alpha)/\Gamma(2\alpha)\cdot t_1^{2-\frac{1}{\alpha}}$
and, by the mean value theorem with $\xi\in(t_1,t_2)$,
$|t_1^{\frac{1}{\alpha}-1}-t_2^{\frac{1}{\alpha}-1}|\le(\frac{1}{\alpha}-1)\max\{2^{\frac{1}{\alpha}-2},1\}t_1^{\frac{1}{\alpha}-2}(t_2-t_1)$
(use $\xi^{\frac{1}{\alpha}-2}\le 2^{\frac{1}{\alpha}-2}t_1^{\frac{1}{\alpha}-2}$ for $\alpha\le\frac12$ and
$\xi^{\frac{1}{\alpha}-2}\le t_1^{\frac{1}{\alpha}-2}$ for $\alpha>\frac12$),
we obtain
\begin{equation}
|I_1|\le C_2\frac{\alpha\Gamma^2(\alpha)}{\Gamma(2\alpha)}
\Big(\frac{1}{\alpha}-1\Big)\max\{2^{\frac{1}{\alpha}-2},1\}(t_2-t_1).
\end{equation}
{\em Estimate the term $I_2$.} Set $r:=(t_2-t_1)/4$. Clearly, $r<t_1/4$. We split $I_2=I_{2,1}+I_{2,2}$,
integrated over $(t_1-r,t_1)$ and $(0,t_1-r)$ respectively.
Since $k(t_1,s)-k(t_2,s)\ge0$ (as $t_1^{1/\alpha}<t_2^{1/\alpha}$ and $\alpha-1<0$),
\begin{equation}
\begin{aligned}
|I_{2,1}|&\le C_2t_2^{\frac{1}{\alpha}-1}
  \int_{t_1-r}^{t_1}(t_1^{1/\alpha}-s^{1/\alpha})^{\alpha-1}\mathrm{d}s
  =C_2t_2^{{1}/{\alpha}-1}
\int_{(t_1-r)^{{1}/{\alpha}}}^{t_1^{{1}/{\alpha}}}
\big(t_1^{{1}/{\alpha}}
-u\big)^{\alpha-1}\alpha u^{\alpha-1}\mathrm{d}u,
\end{aligned}
\end{equation}
which, jointly with $u^{\alpha-1}\leq (t_1-r)^{\frac{\alpha-1}{\alpha}}
\leq 2^{\frac{1-\alpha}{\alpha}}t_1^{\frac{\alpha-1}{\alpha}}$
for $u\in({(t_1-r)^{\frac{1}{\alpha}}},{t_1^{\frac{1}{\alpha}}})$ due to $t_1-r\geq \frac{t_1}{2}$, yields
\begin{eqnarray}\label{compact-J2-2-estim}
|I_{2,1}|&\leq&\nonumber
\alpha C_22^{\frac{1-\alpha}{\alpha}} t_2^{\frac{1}{\alpha}-1}
t_1^{\frac1{\alpha}(\alpha-1)}
\int_{(t_1-r)^{{1}/{\alpha}}}^{t_1^{{1}/{\alpha}}}
\big(t_1^{{1}/{\alpha}}
-u\big)^{\alpha-1}\mathrm{d}u\\
&=&
 C_22^{\frac{1-\alpha}{\alpha}}t_2^{{1}/{\alpha}-1}
t_1^{\frac1{\alpha}(\alpha-1)}
(t_1^{{1}/{\alpha}}-(t_1-r)^{{1}/{\alpha}})^{\alpha}\\
&\leq&\nonumber
C_22^{\frac{2-2\alpha}{\alpha}}
 \frac{T^{1-\alpha}}{\alpha^{\alpha}4^{\alpha}}(t_2-t_1)^{\alpha},
\end{eqnarray}
where the inequalities $t_1^{\frac{1}{\alpha}}-(t_1-r)^{\frac{1}{\alpha}}
\leq\frac{1}{4\alpha}T^{\frac{1}{\alpha}-1}(t_2-t_1)$  and
$t_2^{\frac{1}{\alpha}-1}t_1^{(\alpha-1)/\alpha}\le2^{\frac{1}{\alpha}-1}$
are used.
On $(0,t_1-r)$, by the mean value theorem applied to $u\mapsto u^{\alpha-1}$,
$$
k(t_1,s)-k(t_2,s)=(1-\alpha)\eta_s^{\alpha-2}\big(t_2^{1/\alpha}-t_1^{1/\alpha}\big),
\quad \eta_s^{\alpha-2}\le(t_1^{1/\alpha}-s^{1/\alpha})^{\alpha-2}
\le t_1^{(\frac{1}{\alpha}-1)(\alpha-2)}(t_1-s)^{\alpha-2},
$$
where  $t_1^{1/\alpha}-s^{1/\alpha}\ge t_1^{\frac{1}{\alpha}-1}(t_1-s)$ is used in the last step.
Combining this with $t_2^{1/\alpha}-t_1^{1/\alpha}\le\frac{1}{\alpha}t_2^{\frac{1}{\alpha}-1}(t_2-t_1)$
and $\int_0^{t_1-r}(t_1-s)^{\alpha-2}\mathrm{d}s=\frac{r^{\alpha-1}-t_1^{\alpha-1}}{1-\alpha}\le\frac{r^{\alpha-1}}{1-\alpha}$,
we obtain
\begin{eqnarray}
|I_{2,2}|
&\leq&\nonumber
\frac{(1-\alpha)}{\alpha}C_2t_2^{\frac{2}{\alpha}-2}(t_2-t_1)
t_1^{(\frac1{\alpha}-1)(\alpha-2)}
\int_0^{t_1-r}(t_1-s)^{\alpha-2}\mathrm{d}s\\
&\leq&\label{eq:coompact-est-J21-tem-2}
\frac{C_2}{\alpha}
\big(\frac{t_2}{t_1}\big)^{\frac{2}{\alpha}-2}t_1^{1-\alpha}(t_2-t_1)
r^{\alpha-1}\\
&\leq&\nonumber
\frac{C_2}{\alpha}
2^{\frac{2}{\alpha}-2}\frac{T^{1-\alpha}}{4^{\alpha-1}}(t_2-t_1)^{\alpha}.
\end{eqnarray}
{\em Estimate the term $I_3$.} With $u=s^{1/\alpha}$ and $u^{\alpha-1}\le t_1^{(\alpha-1)/\alpha}$,
$t_2^{\frac{1}{\alpha}-1}t_1^{(\alpha-1)/\alpha}\le2^{\frac{1}{\alpha}-1}$, and
$t_2^{1/\alpha}-t_1^{1/\alpha}=\frac{1}{\alpha}\xi^{1/\alpha-1}(t_2-t_1)\le\frac{1}{\alpha}T^{1/\alpha-1}(t_2-t_1)$, we get
$$
|I_3|\le C_2t_2^{\frac{1}{\alpha}-1}\int_{t_1}^{t_2}k(t_2,s)\mathrm{d}s
\le C_2\,2^{\frac{1}{\alpha}-1}(t_2^{1/\alpha}-t_1^{1/\alpha})^{\alpha}
\le C_2\,2^{\frac{1}{\alpha}-1}\frac{T^{1-\alpha}}{\alpha^{\alpha}}(t_2-t_1)^{\alpha}.
$$
Combining the three estimates $I_1, I_2$ and $I_3$,
we obtain a constant $C_{\alpha,T}$ depending only on $\alpha,T$ such that
$$
|\mathcal{F}(h)(t_1)-\mathcal{F}(h)(t_2)|\le C_{\alpha,T}\big[(t_2-t_1)+(t_2-t_1)^{\alpha}\big],
\quad\forall h\in M,
$$
and the case $t_1=0$ follows from $|\mathcal{F}(h)(t)|\leq C_2t\frac{\alpha\Gamma^2(\alpha)}{\Gamma(2\alpha)}$.
Hence $\mathcal{F}(M)$ is  equicontinuous.
By the Arzel\`a--Ascoli theorem, the operator $\mathcal{F}$ is a compact operator from $C[0,T]$ into itself.
\hfill$\Box$

\begin{lemma}\label{lem-weakly-singular-cont}
{
Let $0<\alpha<1$ and $T>0$, and let $\psi\in C[0,T]\cap C^{0,\beta}_{\rm loc}(0,T]$ for some $\beta\in(1-\alpha,1)$.} Then the map
\begin{equation}\label{Jpsi-def}
t\mapsto(J\psi)(t):=\int_0^t(t-s)^{\alpha-2}(\psi(s)-\psi(t))\,\mathrm{d}s,\qquad 0<t\le T,
\end{equation}
is continuous on $(0,T]$.
\end{lemma}
\noindent{\bf Proof.}
Fix $t_0\in(0,T]$ and put $a:=t_0/2$, $b:=\min\{2t_0,T\}$.
Since $\psi\in C^{0,\beta}_{\rm loc}(0,T]$, there is $C_1>0$ with
$|\psi(x)-\psi(y)|\leq C_1|x-y|^\beta$ for $x,y\in[a/2,b]$.
For $t\in[a,b]$, let $I_a(t)$ and  $I_b(t)$ be given by   
\begin{equation}\label{Vder-cont-split}
I_a(t):=\int_0^a(t-s)^{\alpha-2}(\psi(s)-\psi(t))\mathrm{d}s,\qquad I_b(t):=\int_a^t(t-s)^{\alpha-2}(\psi(s)-\psi(t))\mathrm{d}s,
\end{equation}
where $I_b(a)=0$ and $I_a(a)$ is finite (near $s=a$ the integrand is $\mathcal O((a-s)^{\alpha-2+\beta})\in L^1(0,a)$, as $\alpha-2+\beta>-1$).
Then $(J\psi)(t)$ is split $(J\psi)(t)=I_a(t)+I_b(t)$.

We show the continuity of $I_a$.
For $a<t_1<t_2\le b$, $\tau:=t_2-t_1$,
since $|\frac{\mathrm d}{\mathrm dt}(t-s)^{\alpha-2}|\leq (2-\alpha)(t_1-a)^{\alpha-3}$ for $t\in[t_1,t_2]$ and $s\in[0,a]$,
by the mean value theorem, we have
$$|I_a(t_2)-I_a(t_1)|\le(2-\alpha)(t_1-a)^{\alpha-3}\tau\int_0^a|\psi(s)-\psi(t_2)|\mathrm{d}s
   +|\psi(t_1)-\psi(t_2)|\Big|\int_0^a(t_1-s)^{\alpha-2}\mathrm{d}s\Big|,
$$
which, together with $\psi\in C[0,T]$, gives the continuity of $I_a$.

Next we show the continuity of $I_b$.
With $u:=t-s$, $I_b(t)=\int_0^{t-a}u^{\alpha-2}\big(\psi(t-u)-\psi(t)\big)\mathrm{d}u$,
where $t-u\in[a,t]\subset[a,b]$ gives $|\psi(t-u)-\psi(t)|\le C_1 u^\beta$.
For $a\leq t_1<t_2\leq b$, letting $\phi(t):=\psi(t-u)-\psi(t)$, we have
\begin{equation}\label{Vder-Ib-diff}
I_b(t_2)-I_b(t_1)=\int_0^{t_1-a}u^{\alpha-2}[(\phi(t_2)-(\phi(t_1))]\mathrm{d}u
                   +\int_{t_1-a}^{t_2-a}u^{\alpha-2}\phi(t_2)\mathrm{d}u.
\end{equation}
In the first integral $t_1-u,t_2-u\in[a,b]$, so the integrand is bounded
by $2C_1u^{\alpha-2+\beta}\in L^1(0,t_1-a)$ (as $\alpha-2+\beta>-1$) and converges pointwise to $0$ as $t_2\to t_1$.
Lebesgue dominated convergence theorem yields that the first integral converges to $0$.
The second integral is bounded by
$C_1\int_{t_1-a}^{t_2-a}u^{\alpha-2+\beta}\mathrm{d}u
=C_1(\alpha-1+\beta)^{-1}\big((t_2-a)^{\alpha-1+\beta}-(t_1-a)^{\alpha-1+\beta}\big)
\leq C_1(\alpha-1+\beta)^{-1}\tau^{\alpha-1+\beta}=o(1)$,
where $p:=\alpha-1+\beta\in(0,1)$ and the last inequality is the subadditivity bound
$(x+\tau)^p-x^p\le \tau^p$ with $x=t_1-a\ge0$, valid uniformly as $t_1\to a^+$.
Interchanging $t_1$ and $t_2$ gives left-continuity.
Since $t_0\in(0,T]$ was arbitrary, $(J\psi)(t)$ is continuous on $(0,T]$.
\hfill$\Box$

\section{Properties of solution and a rigorous proof of Theorem 11 of \cite{LiYanCP2009auto}}
\label{Sec-h-diff-property}

In this section, we consider the  properties of the solution of ${}^C_0D_t^{\alpha}V(t)=-h(V(t))$
where $h'(x)$ satisfies local Lipschitz continuity on $(0,+\infty)$.
In this case, $\lim_{x\to0^+}h'(x)$ may not exist.
Examples of such functions include
 $h(x)=x^{p}$ with $p\in(0,1)$, and $h(x)=x^2\sin\frac{1}{x}+4x+x^2$ (with $h(0)=0$). 
Moreover, we provide a rigorous proof of Theorem 11 in \cite{LiYanCP2009auto}.

\begin{theorem}\label{Thm-Dhinfatzero-mono-asympt}
Let $h$ be a class-$\mathcal{K}$ function 
whose derivative $h'$ is local Lipschitz continuous on $(0,+\infty)$, i.e., for any interval $[a,b]$ with $a>0$,
there exists a constant $L_{a,b}$ depending only on $a,b$  such that
$
|h'(x)-h'(y)|\leq L_{a,b}|x-y|$ for all $x,y\in [a,b]$.
Let $A>0$. 
Let $V(t)$ satisfy
\begin{equation}\label{Dh-Inequ-Vh-hinf}
V(t)=A-\frac{1}{\Gamma(\alpha)}\int_0^t(t-s)^{\alpha-1}
h(|V(s)|)\mathrm{d}s.
\end{equation}
Then, (i) $V(t)>0$ for all $t\geq0$. 
(ii) 
The rescaled function $W$ belongs to $C^1[0,+\infty)$ and satisfies
 $W'(t)\leq0$ for $t\geq0$.
(iii) $\lim_{t\to\infty}V(t)=0$.
\end{theorem}
\noindent{\bf Proof.}
For any fixed $\epsilon>0$, define $H_\epsilon(x)$ and $\theta_\epsilon(x)$ as 
\begin{equation}\label{Heps-def-fixhc}
H_\epsilon(x):=h\big(\epsilon^2+\theta_\epsilon(x)\big),
\;\;\;
\theta_\epsilon(x):=\frac{x+\sqrt{x^2+\epsilon^2}}{2}.
\end{equation}
Here $\theta_\epsilon$ is the standard smooth positive part:
$\theta_\epsilon\in C^\infty(\mathbb R)$, $\theta_\epsilon(x)\to x_+=\max\{x,0\}$ as
$\epsilon\to0$, $0\le\theta_\epsilon'\le1$, and
$|\theta_\epsilon''(x)|\le 1/(2\epsilon)$. Since
$\epsilon^2+\theta_\epsilon(x)\ge\epsilon^2>0$, the argument of $h$
stays strictly positive everywhere. We  consider the 
second kind Volterra integral equation
\begin{equation}\label{-U-intequ-mod}
U(t)=A-\frac1{\Gamma(\alpha)}\int_0^t(t-s)^{\alpha-1}H_\epsilon(U(s))\mathrm{d}s,\;t\geq0.
\end{equation}
{\bf Step 1.} Using \dref{-U-intequ-mod}, we construct an approximate solution sequence
$\{V_n\}$ for \dref{-U-intequ-mod}. 
Let $V_0(t)$ be an arbitrary function in $C^{1}[0,\infty)$ with $V_0(0)=A$.
For $n\geq1$, we define $V_n$ by
\begin{equation}\label{hinf-Vn-def}
V_{n+1}(t)=A-\frac1{\Gamma(\alpha)}\int_0^t(t-s)^{\alpha-1}
H_\epsilon(V_n(s))\mathrm{d}s,\;\;n=1,2,\cdots.
\end{equation}
From \dref{hinf-Vn-def}, it is easy to see
that $V_n(t)\leq A$ for all $t\geq 0$ and $n=1,2,3,\cdots$.
On the other hand, it is trivial to see that 
for any fixed $T>0$ and for all $n\geq1$,
\begin{eqnarray*}
V_{n}(t)
&\geq &
 -\frac1{\Gamma(\alpha)}\int_0^t(t-s)^{\alpha-1}
h(\epsilon^2+\theta_\epsilon(B))\mathrm{d}s\\
&\geq &
-\frac{T^{\alpha}}{\Gamma(\alpha+1)}h(\epsilon^2+\theta_\epsilon(B)),\; t\in [0,T],
\end{eqnarray*}
where $B=\max\{A,\|V_0\|_{C[0,T]}\}$.
 Hence $\{V_n\}$ is a uniformly  bounded sequence in $C[0,T]$.
Observe that the error sequence $\{V_{n+1}-V_n\}$ satisfies
\begin{equation}\label{hinf-Vn-def-err-estm}
|V_{n+2}(t)-V_{n+1}(t)|
\leq
\frac1{\Gamma(\alpha)}\int_0^t(t-s)^{\alpha-1}
|H_{\epsilon}(V_{n+1}(s))-H_{\epsilon}(V_n(s))|\mathrm{d}s
\end{equation}
with $H_{\epsilon}(x)=h(\epsilon^2+\theta_\epsilon(x))$. 
Since $h\in C^1(0,\infty)$ with $h'$ locally Lipschitz on $(0,\infty)$ while $h'(0)$ may not exist,
the map $H_\epsilon$ is of class $C^{1,1}_{\rm loc}(\mathbb R)$ because $\theta_\epsilon$
is smooth, and the shift $\epsilon^2+\theta_\epsilon(x)\ge\epsilon^2$ keeps
the argument of $h$ strictly positive. Moreover $H_\epsilon'$ is Lipschitz
continuous on every bounded interval, with Lipschitz constant bounded in
terms of $\epsilon$, the local Lipschitz constant of $h'$ on
$[\epsilon^2/2,\epsilon^2+\theta_\epsilon(A)]$, and $\|\theta_\epsilon''\|_\infty\le
1/(2\epsilon)$. The derivative is
\begin{equation}\label{Heps-weak-der-fix}
H_\epsilon'(x)
=h'\big(\epsilon^2+\theta_\epsilon(x)\big)\,
 \frac12\Big(1+\frac{x}{\sqrt{x^2+\epsilon^2}}\Big),
\qquad x\in\mathbb R.
\end{equation}
Let $C_{\epsilon}=\max\{h'(x):
x\in[\epsilon^2,\epsilon^2+\theta_\epsilon(A)]\}$.
Choose $r$ so that $1<r<\frac1{1-\alpha}$.
Noticing that
$$
\int_0^t(t-s)^{(\alpha-1)r} \mathrm{d}s
=t^{1+(\alpha-1)r}\int_0^1(1-s)^{(\alpha-1)r}\mathrm{d}s
\leq \frac{T^{1+(\alpha-1)r}}{{1+(\alpha-1)r}},
$$
by \dref{hinf-Vn-def-err-estm}  and H\"{o}lder's inequality, we obtain
\begin{eqnarray*}
|V_{n+2}(t)-V_{n+1}(t)|
&\leq&
\frac{C_{\epsilon}}{\Gamma(\alpha)}\int_0^t(t-s)^{\alpha-1}
|V_{n+1}(s)-V_n(s)|\mathrm{d}s\\
&\leq&
\frac{C_{\epsilon}}{\Gamma(\alpha)}
\left(\int_0^t(t-s)^{(\alpha-1)r} \mathrm{d}s  \right)^{1/r}
\left(\int_0^t |V_{n+1}(s)-V_n(s)|^{r'} \mathrm{d}s \right)^{1/r'}\\
&\leq&
\frac{C_{\epsilon}T^{\frac{1}{r}+(\alpha-1)}}
{(1+(\alpha-1)r)^{\frac{1}{r}}\Gamma(\alpha)}
\left(\int_0^t |V_{n+1}(s)-V_n(s)|^{r'} \mathrm{d}s \right)^{1/r'},
\ \ \ t\in [0,T],
\end{eqnarray*}
where $r'=\frac{r}{r-1}>1$.
Raising both sides to the $r'$-th power and
by Lemma \ref{lem-sequ-estim-I0}, 
for any $\rho>0$, we have
\begin{equation}\label{rapid-decay-error-sequence-small-p}
\sum_{n=1}^{\infty}\|V_{n+1}-V_n\|_{C[0,T]}^{\rho}
\leq C_1\|V_{2}-V_1\|_{C[0,T]}^{\rho}
\end{equation}
for some $C_1>0$. Particularly, taking $\rho=1$, then \dref{rapid-decay-error-sequence-small-p}
 implies that  $V_n$ is uniformly convergent, and
let $V^{\epsilon}$ be the limit of $V_n$ in $C[0,T]$ for any $T>0$.
From \dref{hinf-Vn-def}, it is seen that
$V^{\epsilon}$ satisfies
\begin{equation}\label{equation-approximate-solution-small-p}
V^{\epsilon}(t)=A-\frac1{\Gamma(\alpha)}\int_0^t(t-s)^{\alpha-1}
H_{\epsilon}(V^{\epsilon}(s))\mathrm{d}s,~~~t\in [0,T].
\end{equation}
{\bf Step 2.} For any fixed $\epsilon>0$, we claim that $V^{\epsilon}(t)$ is decreasing with respect to $t$.
For this purpose, we turn to the estimate of the derivative $(V_n(t^{1/\alpha}))'$. Let $W_n(t)=V_n(t^{1/\alpha})$. 
It follows from \dref{hinf-Vn-def} that
\begin{equation*}
W_{n+1}^{\epsilon}(t)=A-\frac1{\Gamma(\alpha)}\int_0^{t^{\frac{1}{\alpha}}}(t^{\frac{1}{\alpha}}-s)^{\alpha-1}
H_{\epsilon}(W_n(s^{\alpha}))\mathrm{d}s.
\end{equation*}
Since $t\mapsto t^{1/\alpha}\in C^1[0,T^\alpha]$ with $(t^{1/\alpha})'=\frac1\alpha t^{(1-\alpha)/\alpha}\to0$ as $t\to0^+$
 ($\alpha\in(0,1)$), for every $V_0\in C^1[0,\infty)$ the composition
 $W_0(t):=V_0(t^{1/\alpha})$ belongs to $C^1[0,T^\alpha]$ and satisfies $W_0'(0)=0$.
By induction, we are able to show that each $W_n$ is continuously differentiable
  and its derivative $W_n'$  satisfies the following recurrence relationship
\begin{equation}\label{dotWn-recur-ship-1}
\begin{array}{l}\disp
W_{n+1}'(t)\crr\disp
=
-\frac{H_{\epsilon}(A)}{\Gamma(\alpha)}
t^{\frac{\alpha-1}{\alpha}}(t^{1/\alpha})'
-
\frac{(t^{1/\alpha})'}{\Gamma(\alpha)}
\int_0^{t^{1/\alpha}}\hspace{-2mm}
s^{\alpha-1}
H_{\epsilon}'\big( W_n((t^{1/\alpha}-s)^\alpha) \big)\crr\disp\hspace{6cm}\times W_n'((t^{1/\alpha}-s)^\alpha)
\alpha(t^{1/\alpha}-s)^{\alpha-1}\mathrm{d}s\crr\disp
=
-\frac{H_{\epsilon}(A)}{\Gamma(\alpha+1)}
-\frac{ t^{\frac{1}{\alpha} -1 }}{\Gamma(\alpha+1)}
\int_0^{t}(t^{1/\alpha}-s^{1/\alpha})^{\alpha-1}
H_{\epsilon}'(W_{n}(s))W_{n}^{'}(s)\mathrm{d}s,
\end{array}
\end{equation}
which shows $W_n(t)\in C^1[0,+\infty)$ for all $n\geq0$.
Note that the identity \dref{dotWn-recur-ship-1} is obtained by integrating by parts
in the formula for $W_{n+1}$ above, which is legitimate
since $\phi(s):=H_\epsilon(W_n(s^\alpha))$ is absolutely continuous on $[0,T]$
 ($W_n\in C^1[0,T^\alpha]$, $s\mapsto s^\alpha$ is absolutely continuous, and $H_\epsilon$ is Lipschitz on every bounded interval).
This integration by parts replaces the weakly singular kernel $(u-s)^{\alpha-1}$
by the smooth kernel $(u-s)^\alpha$, and the resulting expression is then differentiated pointwise with respect to $t$,
applying the chain rule $\frac{\mathrm{d}}{\mathrm{d}s}H_\epsilon(W(s))=H_\epsilon'(W(s))W'(s)$
to $\phi$.

Let $X_n(t)=-W_n'(t)$. Then \dref{dotWn-recur-ship-1} becomes
\begin{eqnarray}\label{equation-derivative-sequence-small-p}
X_{n+1}(t)
=
\frac{H_{\epsilon}(A)}{\Gamma(\alpha+1)}
-\frac{ t^{\frac{1}{\alpha} -1 }}{\Gamma(\alpha+1)}
\int_0^{t}(t^{1/\alpha}-s^{1/\alpha})^{\alpha-1}
H_{\epsilon}'(W_{n}(s))X_{n}(s)\mathrm{d}s.
\end{eqnarray}
Using a simple change of variables, from \dref{equation-derivative-sequence-small-p}, we get  
\begin{eqnarray*}
X_{n+1}(t)
=
\frac{H_{\epsilon}(A)}{\Gamma(\alpha+1)}
-\frac{t}{\Gamma(\alpha+1)}
\int_0^{1}(1-s^{1/\alpha})^{\alpha-1}
H_{\epsilon}'(W_{n}(ts))X_{n}(ts)\mathrm{d}s.
\end{eqnarray*}
Choosing an exponent $r$ satisfying $r\in(1,\frac1{1-\alpha})$, we apply H\"{o}lder's inequality to obtain
\begin{eqnarray*}
|X_{n+1}(t)|
&\leq&
\frac{H_{\epsilon}(A)}{\Gamma(\alpha+1)}
+\frac{t\|H_{\epsilon}'\|_{L^{\infty}}}{\Gamma(\alpha+1)}
\left(\int_0^1(1-s^{\frac1{\alpha}})^{(\alpha-1)r}\mathrm{d}s\right)^{\frac1{r}}
\left(\int_0^1 |X_n(ts)|^{r'}\mathrm{d}s\right)^{\frac1{r'}}\\
&\leq&
\frac{H_{\epsilon}(A)}{\Gamma(\alpha+1)}
+
C_2\left(\int_0^t |X_n(s)|^{r'}\mathrm{d}s\right)^{\frac1{r'}},~~~0\leq t \leq T,
\end{eqnarray*}
where $C_2$ is a constant depending only on $\epsilon,\alpha,r,$ and $T$. Hence
\begin{eqnarray*}
|X_{n+1}(t)|^{r'}
&\leq&
2^{r'-1}\frac{H_{\epsilon}(A)^{r'}}{\Gamma(\alpha+1)^{r'}}
+
2^{r'-1}C_2^{r'}\int_0^t |X_n(s)|^{r'}\mathrm{d}s,~~~0\leq t \leq T^{\alpha}.
\end{eqnarray*}
It follows from Lemma \ref{lem-sequ-estim} that
$$
|X_{n+1}(t)|^{r'}
\leq \Big(2^{r'-1}\frac{H_{\epsilon}(A)^{r'}}{\Gamma(\alpha+1)^{r'}}
+\|X_{0}\|^{r'}_{C[0,T^{\alpha}]}
\Big)e^{2^{r'-1}C_2^{r'}T^{\alpha}},\ \ \ \forall t\in[0,T^{\alpha}],
$$
which shows that  $\{X_n\}$ are uniformly bounded in $C[0,T^{\alpha}]$.

Observe that $H_{\epsilon}'$ is non-negative and bounded. We can estimate the error sequence $X_{n+1}-X_n$ by
\begin{equation}
\begin{array}{l}\disp
|X_{n+2}(t)-X_{n+1}(t)|\crr\disp
\leq
\frac{ t^{\frac{1}{\alpha} -1 }}{\Gamma(\alpha+1)}
\int_0^{t}(t^{1/\alpha}-s^{1/\alpha})^{\alpha-1}
|H_{\epsilon}'(W_{n+1}(s))X_{n+1}(s)-
H_{\epsilon}'(W_{n}(s))X_{n}(s)|\mathrm{d}s\crr\disp
\leq
\frac{ t^{\frac{1}{\alpha} -1 }}{\Gamma(\alpha+1)}
\int_0^{t}(t^{1/\alpha}-s^{1/\alpha})^{\alpha-1}
H_{\epsilon}'(W_{n+1}(s))|X_{n+1}(s)-X_n(s)|\mathrm{d}s\crr\disp\hspace{0.3cm}
+
\frac{ t^{\frac{1}{\alpha} -1 }}{\Gamma(\alpha+1)}
\int_0^{t}(t^{1/\alpha}-s^{1/\alpha})^{\alpha-1}
|H_{\epsilon}'(W_{n+1}(s))-H_{\epsilon}'(W_{n}(s))||X_n(s)|\mathrm{d}s.
\end{array}
\end{equation}
Now $H_\epsilon'$ is genuinely Lipschitz (by the smooth construction \eqref{Heps-def-fixhc}),
so by the Lipschitz continuity of $H_{\epsilon}'$ and the uniform boundedness of $\{X_n\}$ , it follows that
\begin{eqnarray*}
|X_{n+2}(t)-X_{n+1}(t)|
&\leq&
C_3\left(\int_0^t |X_{n+1}(s)-X_n(s)|^{r'}\mathrm{d}s\right)^{\frac1{r'}}
+
C_4\|W_{n+1}-W_n\|_{C[0,T^{\alpha}]},
\end{eqnarray*}
for $0\leq t \leq T^\alpha$, where $C_3$ and $C_4$ are constants depending only on $\epsilon,\alpha,r,$ and $T$. This implies
\begin{eqnarray*}
|X_{n+2}(t)-X_{n+1}(t)|^{r'}
&\leq&
2^{r'-1}C_1^{r'}\int_0^t |X_{n+1}(s)-X_n(s)|^{r'}\mathrm{d}s
+
2^{r'-1}C_2^{r'}\|W_{n+1}-W_n\|_{C[0,T^{\alpha}]}.
\end{eqnarray*}
The convergence result (\ref{rapid-decay-error-sequence-small-p}) implies that for any $\rho>0$ there holds
$
\sum_{n=1}^{\infty}
\|W_{n+1}-W_n\|_{C[0,T^{\alpha}]}^{\rho}<+\infty.
$ 
As a consequence of Lemma \ref{lem-sequ-estim}, we get
\begin{equation*}
\sum_{n=1}^{\infty}
\|X_{n+1}-X_n\|_{C[0,T^{\alpha}]}^{\rho}<+\infty.
\end{equation*}
for any exponent $\rho>0$. In particular, $\sum_{n=1}^{\infty}\|X_{n+1}-X_n\|_{C[0,T^{\alpha}]}<\infty$,
which leads to the existence of the uniform limit of $X_n$ in $C[0,T^{\alpha}]$.
Let $X^{\epsilon}$ be the uniform limit of $X_n$ in $C[0,T^{\alpha}]$. Taking $n\rightarrow\infty$ in (\ref{equation-derivative-sequence-small-p}) gives
\begin{equation}\label{Xepsi-int-equ}
X^{\epsilon}(t)
=
\frac{H_{\epsilon}(A)}{\Gamma(\alpha+1)}
-\frac{ t^{\frac{1}{\alpha} -1 }}{\Gamma(\alpha+1)}
\int_0^{t}(t^{1/\alpha}-s^{1/\alpha})^{\alpha-1}
H_{\epsilon}'(W^{\epsilon}(s))X^{\epsilon}(s)\mathrm{d}s,
\end{equation}
where $W^{\epsilon}\in C[0,T^{\alpha}]$ is the limit of $W_n$ in $C[0,T^{\alpha}]$, so that $W^{\epsilon}(t)=V^{\epsilon}(t^{1/\alpha})$.
On the other hand, we show that $X^{\epsilon}$ is a non-negative continuous function.
To this end, for the fixed $T>0$, we define the operator
 $\mathcal{S}_\epsilon:C[0,T^\alpha]\to C[0,T^\alpha]$ by
\begin{equation}\label{S-eps-def}
\mathcal{S}_\epsilon(U)(t)=\frac{H_\epsilon(A)}{\Gamma(\alpha+1)}
-\frac{t^{1/\alpha-1}}{\Gamma(\alpha+1)}\int_0^t(t^{1/\alpha}-s^{1/\alpha})^{\alpha-1}
H_\epsilon'(W^\epsilon(s))\,U_+(s)\mathrm{d}s.
\end{equation}
Since $u\mapsto u_+:=\max\{u,0\}$ is Lipschitz continuous and $H_\epsilon'(W^\epsilon(\cdot))\in L^\infty[0,T^\alpha]$
(bounded measurable: $h'$ is bounded on $[\epsilon^2,\epsilon^2+\theta_\epsilon(A)]$),
the operator $\mathcal{S}_\epsilon$ is continuous.
By Lemma \ref{Ne-operator-com-lem}, we know that
the operator $\mathcal{S}_\epsilon$ is compact on $C[0,T^\alpha]$.

We claim that every $U\in C[0,T^\alpha]$ and $\sigma\in[0,1]$
with $U=\sigma\mathcal{S}_\epsilon(U)$
satisfies $U\ge0$ on $[0,T^\alpha]$.
Indeed, if $U$ took a negative value, then since $U(0)=\sigma H_\epsilon(A)/\Gamma(\alpha+1)\ge0$,
there would exist $t_0>0$ and $\delta>0$ with $U(t_0)=0$ and $U(t)<0$ for $t\in(t_0,t_0+\delta]$.
Hence $U_+(s)=0$ for $s\in[t_0,t_0+\delta]$ and \dref{S-eps-def} gives, for $t\in[t_0,t_0+\delta]$,
$$
\begin{array}{l}\disp
U(t)=\frac{\sigma H_\epsilon(A)}{\Gamma(\alpha+1)}
-\frac{\sigma t^{1/\alpha-1}}{\Gamma(\alpha+1)}\int_0^t(t^{1/\alpha}-s^{1/\alpha})^{\alpha-1}
H_\epsilon'(W^\epsilon(s))\,U_+(s)\mathrm{d}s\crr\disp
\hspace{0.9cm}=\frac{\sigma H_\epsilon(A)}{\Gamma(\alpha+1)}
-\frac{\sigma}{\Gamma(\alpha+1)}\int_0^{t_0}\Big(1-\big(s/t\big)^{1/\alpha}\Big)^{\alpha-1}
H_\epsilon'(W^\epsilon(s))\,U_+(s)\mathrm{d}s,
\end{array}
$$
whose right hand side is nondecreasing in $t$ (because $H_\epsilon'(W^\epsilon(s))U_+(s)\ge0$
and $t\mapsto\big(1-(s/t)^{1/\alpha}\big)^{\alpha-1}$ is decreasing).
Hence $U(t)\ge U(t_0)=0$ for $t\in(t_0,t_0+\delta]$, a contradiction,
so $U\ge0$ on $[0,T^\alpha]$, and consequently
$0\le U(t)\le\sigma H_\epsilon(A)/\Gamma(\alpha+1)$ for all $t\in[0,T^\alpha]$.
By Lemma \ref{Leray-Schauder-fixed-point-lem}, $\mathcal{S}_\epsilon$ admits a fixed point $U_\epsilon$,
which is non-negative and satisfies \dref{Xepsi-int-equ} (since $U_{\epsilon,+}=U_\epsilon$).
Since $X^\epsilon$ satisfies the same equation, the difference $Q:=X^\epsilon-U_\epsilon$ obeys
$$
Q(t)=-\frac{t^{1/\alpha-1}}{\Gamma(\alpha+1)}\int_0^t(t^{1/\alpha}-s^{1/\alpha})^{\alpha-1}
H_\epsilon'(W^\epsilon(s))\,Q(s)\mathrm{d}s.
$$
H\"older's inequality with $r\in(1,\frac{1}{1-\alpha})$ and $r'=\frac{r}{r-1}$
gives $|Q(t)|\le C_5\big(\int_0^t|Q(s)|^{r'}\mathrm{d}s\big)^{1/r'}$ for $t\in[0,T^\alpha]$ and some $C_5>0$. 
Raising to the power $r'$ and applying Gronwall's inequality yields $Q\equiv0$ on $[0,T^\alpha]$.
Since $T>0$ is arbitrary, $X^\epsilon=U_\epsilon\ge0$ on $[0,\infty)$.
Since $W_n\rightrightarrows W^{\epsilon}$ and $W_n'=-X_n\rightrightarrows -X^{\epsilon}$ in $C[0,T^\alpha]$,
the classical differentiation theorem for uniformly convergent $C^1[0,T^\alpha]$ sequences yields $W^{\epsilon}\in C^1[0,T^\alpha]$
and $(W^{\epsilon})'=-X^{\epsilon}$, which is precisely the identity stated below:
\begin{equation}
X^{\epsilon}(t)=-\frac{d}{dt}\Big(V^{\epsilon}(t^{1/\alpha}) \Big),~~~t\geq 0.
\end{equation}
Therefore, for any $\epsilon>0$, $V^{\epsilon}$ is decreasing with respect to $t$.  \\
{\bf Step 3.} We claim that $\{V^{\epsilon}:0<\epsilon<1\}$ is precompact in $C[0,T]$. 
By  the non-negativity of $X^{\epsilon}(t)$ for all $t\geq0$ and   \dref{Xepsi-int-equ}, it follows that
$0\leq X^{\epsilon}(t)\leq \frac{H_{\epsilon}(A)}{\Gamma(\alpha+1)}$.
Observe that, for $0<\epsilon<1$, we have the uniform boundedness of $H_{\epsilon}(V^{\epsilon})$,
i.e. $0\leq H_{\epsilon}(V^{\epsilon})\leq h(2+A)$. Hence for $0<\epsilon<1$, we get
\begin{equation*}
|V^{\epsilon}(t_1^{\frac{1}{\alpha}})-V^{\epsilon}(t_2^{\frac{1}{\alpha}})|
=\Big|\int_{t_1}^{t_2}X^{\epsilon}(s)\mathrm{d}s\Big|
\leq
\frac{H_{\epsilon}(A)}{\Gamma(\alpha+1)}
|t_1-t_2|
\leq
\frac{h(2+A)}{\Gamma(\alpha+1)}
|t_1-t_2|,
\end{equation*}
which, together with the inequality $|a^{\alpha}-b^{\alpha}|\leq|a-b|^{\alpha}$ for all $a,b\geq0$ and $\alpha\in(0,1)$, gives
\begin{equation}\label{eq:Veps-Holder-inequality}
|V^{\epsilon}(t_1)-V^{\epsilon}(t_2)|
\leq
\frac{h(2+A)}{\Gamma(\alpha+1)}|t_1-t_2|^{\alpha},
\end{equation}
for all $t_1,t_2\in [0,T]$. By the Arzel\`a--Ascoli theorem, we see that $\{V^{\epsilon}:0<\epsilon<1\}$ is precompact in $C[0,T]$.

We now construct a common limit on $[0,\infty)$ by a diagonal argument.
For each $k\in\mathbb{N}$, in view of the precompactness in $C[0,k]$,
we choose a sequence $\{\epsilon^{(k)}_n\}_{n=1}^{\infty}\subset(0,1)$ with $\epsilon^{(k)}_n\to0$ as $n\to\infty$
such that $\{V^{\epsilon^{(k)}_n}\}_{n=1}^{\infty}$ is a Cauchy sequence in $C[0,k]$,
and we may assume, by passing to a further subsequence if necessary,
that $\{\epsilon^{(k)}_n\}_{n=1}^{\infty}$ is a subsequence of $\{\epsilon^{(k-1)}_n\}_{n=1}^{\infty}$ for each $k\geq2$.
Then the diagonal sequence $\epsilon_n:=\epsilon^{(n)}_n$, $n\in\mathbb{N}$, satisfies $\epsilon_n\to0$,
and for every fixed $k\in\mathbb{N}$ the tail $\{V^{\epsilon_n}\}_{n\geq k}$ is a subsequence of
$\{V^{\epsilon^{(k)}_n}\}_{n=1}^{\infty}$.
Hence $\{V^{\epsilon_n}\}$ is a Cauchy sequence in $C[0,k]$ for every $k\in\mathbb{N}$. Therefore the pointwise limit
\begin{equation}\label{tildeV-def}
\widetilde{V}_A(t):=\lim_{n\to\infty}V^{\epsilon_n}(t),\qquad t\geq0,
\end{equation}
is well defined, and $\widetilde{V}_A\in C[0,\infty)$ with $V^{\epsilon_n}\to\widetilde{V}_A$
uniformly on $[0,k]$ for every $k\in\mathbb{N}$.
Taking the limit $n\to\infty$ in \dref{equation-approximate-solution-small-p} with $\epsilon=\epsilon_n$
on the fixed interval $[0,T]$, and using the uniform convergence $V^{\epsilon_n}\to\widetilde{V}_A$ in $C[0,T]$,
the Lipschitz continuity of $u\mapsto u_+:=\max\{u,0\}$, the continuity of $h$,
and the integrability of the weakly singular kernel $(t-s)^{\alpha-1}$, we obtain 
\begin{equation}\label{tildeV-Vol-equ}
\widetilde{V}_A(t)=A-\frac1{\Gamma(\alpha)}\int_0^t(t-s)^{\alpha-1}
h(\widetilde{V}_{A,+}(s))\mathrm{d}s,~~~t\in [0,T].
\end{equation}
for every $T>0$, that is, \dref{tildeV-Vol-equ} holds for all $t\geq0$.

Next, we use the monotonicity argument 
to verify that $\widetilde{V}_A$ is a positive function in $C[0,\infty)$,
i.e.,
$\widetilde{V}_A(t)>0$ for all $t\in[0,T]$
with any $T>0$.
Indeed, since each $V^{\epsilon}(t)$ is decreasing with respect to $t$
and $\widetilde{V}_A$ is the uniform limit of $V^{\epsilon_n}$,
we conclude that $\widetilde{V}_A(t)$ is decreasing with respect to $t$ on $[0,T]$.
Now, suppose there exists $t_1$ such that $\widetilde{V}_A(t_1)=0$.
Then by the fact that $\widetilde{V}_A(t)$ is decreasing with respect to $t$,
we have $\widetilde{V}_A(t)\leq0$ for all $t\geq t_1$,
which leads $h(\widetilde{V}_{A,+}(s))=0$ for all $s\geq t_1$.
Hence, \dref{tildeV-Vol-equ} is reduced to
\begin{equation}\label{tildeV-Vol-equ-I1}
\begin{array}{l}\disp
A=\widetilde{V}_A(t)+\frac1{\Gamma(\alpha)}\int_0^{t_1}(t-s)^{\alpha-1}
h(\widetilde{V}_{A,+}(s))\mathrm{d}s,\crr\disp\hspace{5mm}
\leq
\frac{C_6}{\Gamma(\alpha)}\int_0^{t_1}(t-s)^{\alpha-1}
\mathrm{d}s
=\frac{C_6}{\Gamma(\alpha+1)}(t^{\alpha}-(t-t_1)^{\alpha}),
\ \ \
\forall t\geq t_1,
\end{array}
\end{equation}
where $C_6=\max\{h(\widetilde{V}_{A,+}(s)):s\in[0,t_1]\}$.
Clearly, $\lim_{t\to\infty}(t^{\alpha}-(t-t_1)^{\alpha})=0$. From this, 
by taking $t\to\infty$ in \dref{tildeV-Vol-equ-I1}, we get $0<A\leq0$, 
which is a contradiction. Hence $\widetilde{V}_A(t)>0$ for all $t\geq0$, which gives $\widetilde{V}_{A,+}(t)=\max\{\widetilde{V}_{A}(t),0\}=\widetilde{V}_{A}(t)$ for all $t\geq0$.
\\
{\bf Step 4.} We claim that $V(t)>0$ for all $t\geq0$.
Let $Y(t)=V(t)-\widetilde{V}_{A/2}(t)$.
Since $Y(0)=A-A/2=A/2>0$, by the continuity of $Y(t)$,
there exists $\delta>0$ such that $Y(t)>0$ for all $t\in[0,\delta]$.
Now, we claim that $Y(t)>0$ for all $t\geq0$.
Otherwise, there exists $t_1$ such that $Y(t_1)=0$ and $Y(t)>0$ for all $t\in[0,t_1)$.
By
\dref{Dh-Inequ-Vh-hinf} and \dref{tildeV-Vol-equ}, and noting $V(t)\geq \widetilde{V}_{A/2}(t)>0$ for $t\in[0,t_1]$, we have
\begin{equation}\label{Dh-Inequ-hinf-Z}
Y(t)=\frac{A}{2}-\frac{1}{\Gamma(\alpha)}
\int_0^t(t-s)^{\alpha-1}
[h(V(s))-h(\widetilde{V}_{A/2}(s))]\mathrm{d}s,\ \ \ t\in[0,t_1],
\end{equation}
which is equivalent to
${}^C_0D_t^{\alpha}Y(t)=-[h(V(t))-h(\widetilde{V}_{A/2}(t))]$ with $Y(0)=A/2$.
Since  $h(V(t))-h(\widetilde{V}_{A/2}(t))=h'(\xi_t)Y(t)$ with
     $\xi_t\in[\widetilde{V}_{A/2}(t),V(t)]$ for $t\in[0,t_1]$,
we define $L_0=2\max\{h'(x):x\in[\widetilde{V}_{A/2}(t_1),\max\{V(s):s\in[0,t_1]\}]\}$. Note that
$L_0<+\infty$ due to $\widetilde{V}_{A/2}(t_1)>0$.
Then, it is easy to see that
${}^C_0D_t^{\alpha}Y(t)=-L_0Y(t)+(L_0-h'(\xi_t))Y(t)$ with $Y(0)=A/2$. Thus,
\dref{Dh-Inequ-hinf-Z} is further equivalent to
\begin{equation}\label{Z-equi-L0-intequ}
Y(t)=E_{\alpha}(-L_0t^{\alpha})\frac{A}{2}+\int_0^t(t-s)^{\alpha-1}E_{\alpha,\alpha}(-L_0 (t-s)^{\alpha})(L_0-h'(\xi_s))Y(s)\mathrm{d}s,\ \ \
t\in[0,t_1],
\end{equation}
Taking $t=t_1$ in \dref{Z-equi-L0-intequ}
and noting $E_{\alpha,\alpha}(-L_0 (t_1-s)^{\alpha})>0$ (see Lemma \ref{lem-1d-ML-est})
and $L_0-h'(\xi_s)\geq0$ for $s\in[0,t_1]$,
we derive $0=Y(t_1)\geq E_{\alpha}(-L_0t_1^{\alpha})\frac{A}{2}>0$, which is a contradiction.
Therefore, $Y(t)>0$ for all $t\geq0$, which yields that $V(t)\geq \widetilde{V}_{A/2}(t)>0$ for all $t\geq0$.  This proves (i).\\
{\bf Step 5.}
Next, we show that $V(t)=\widetilde{V}_A(t)$ for all $t\geq0$ and
${V}(t)$ is decreasing with respect to $t$ on $[0,+\infty)$.
Suppose that $V(t)\not\equiv\widetilde V_A(t)$ on $[0,+\infty)$, and set
$Z(t)=V(t)-\widetilde V_A(t)$, so that $Z(0)=0$. Fix an arbitrary $T>0$.
Since $V(t)>0$ by Step~4 and $\widetilde V_A(t)>0$ by Step~3, and both are
bounded above by $A$, the number
$
m_T:=\min_{t\in[0,T]}\min\{V(t),\widetilde V_A(t)\}>0
$
is positive. As $h'$ is bounded on the compact interval $[m_T,A]$, there is
$L_T<\infty$ such that $h$ is $L_T$-Lipschitz on $[m_T,A]$. Hence, by the
mean value theorem,
\[
|h(V(s))-h(\widetilde V_A(s))|\le L_T|Z(s)|,\qquad s\in[0,T].
\]
Subtracting \dref{tildeV-Vol-equ} from \dref{Dh-Inequ-Vh-hinf} (and using
$|V(s)|=V(s)$ and $\widetilde V_{A,+}(s)=\widetilde V_A(s)$) gives
\begin{equation}\label{Dh-Inequ-hinf-Z2}
Z(t)=-\frac1{\Gamma(\alpha)}\int_0^t(t-s)^{\alpha-1}
\big[h(V(s))-h(\widetilde V_A(s))\big]\,\mathrm ds,\qquad t\in[0,T].
\end{equation}
Therefore, $|Z(t)|\le \frac{L_T}{\Gamma(\alpha)}\,({}_0I_t^\alpha|Z|)(t)$, $t\in[0,T]$.
Iterating this inequality and using that the fractional integration operator
${}_0I_t^\alpha$ is positive with
${}_0I_t^{n\alpha}1=t^{n\alpha}/\Gamma(n\alpha+1)$, we obtain, for every
$n\in\mathbb N$,
\[
|Z(t)|\leq{L_T^n}({}_0I_t^{n\alpha}|Z|)(t)
\le \frac{(L_TT^\alpha)^n}{\Gamma(n\alpha+1)}
\|Z\|_{C[0,T]},\ \ \  t\in[0,T].
\]
The right-hand side tends to $0$ as $n\to\infty$, since $\Gamma(n\alpha+1)$
grows faster than exponentially. Thus $Z(t)=0$ for all $t\in[0,T]$. Since $T>0$
was arbitrary, $Z(t)\equiv0$, i.e.\ $V(t)=\widetilde V_A(t)$ for all $t\ge0$.
Now, from the fact that $\widetilde{V}_A(t)$ is decreasing with respect to $t$ on $[0,T]$ with any fixed $T>0$, we derive that
${V}(t)$ is decreasing with respect to $t$ on $[0,+\infty)$.\\
{\bf Step 6.} We prove the differentiability of $W(t)$ and $W'(t)\leq0$. 
 Since $V(t)>0$ on $t\in[0,\infty)$ (established in Step 4)
and $V$ satisfies \dref{Dh-Inequ-Vh-hinf}, we have $h(|V(s)|)=h(V(s))$.
We write $\psi(t):=h(V(t))$,  so that $V=A-{}_0I_t^\alpha\psi$ on every interval $[0,T]$.
First, we show $V\in C^{0,\alpha}[0,T]$. Indeed, since $V(t)=\widetilde{V}_A(t)$,
by \dref{eq:Veps-Holder-inequality} and \dref{tildeV-def},
for $0\le t_1<t_2\le T$, we have
\begin{equation}\label{V-Holder-reg-B1}
|V(t_2)-V(t_1)|
\leq\frac{h(2+A)}{\Gamma(\alpha+1)}(t_2-t_1)^\alpha.
\end{equation}
Since $V$ is continuous and positive on $[0,T]$, $\delta_T:=\min_{t\in[0,T]}V(t)>0$, $M_T:=\max_{t\in[0,T]}V(t)$,
and $h\in C^1[\delta_T,M_T]$ with $h'$ Lipschitz continuous there,  we know that $\psi\in C^{0,\alpha}[0,T]$.
We use the local smoothing property of the Abel operator  (\cite[Theorems 3.1]{SamkoKilbasMarichev1993}):
if $\psi\in C^{0,\gamma}$ on a neighbourhood of $t_0>0$, then ${}_0I_t^\alpha\psi\in C^{0,\gamma+\alpha}$
on a neighbourhood of $t_0$, provided $\gamma+\alpha<1$.
Iterating in \dref{Dh-Inequ-Vh-hinf} finitely many times (possible since $\alpha>0$),
we obtain $V\in C^{0,\beta}_{\rm loc}(0,T]$ for some $\beta\in(1-\alpha,1)$,
hence $\psi\in C^{0,\beta}_{\rm loc}(0,T]$.
To obtain the differentiability, fix $t>0$ and let $\tau>0$ be small. From \dref{Dh-Inequ-Vh-hinf},
\begin{equation}\label{V-dq-B1}
\begin{array}{l}\disp
\frac{V(t+\tau)-V(t)}{\tau}=-\frac{1}{\Gamma(\alpha)\tau}\Big[\int_0^t\big((t+\tau-s)^{\alpha-1}-(t-s)^{\alpha-1}\big)\psi(s)\mathrm{d}s\crr\disp
\hspace{5cm}+\int_t^{t+\tau}(t+\tau-s)^{\alpha-1}\psi(s)\mathrm{d}s\Big].
\end{array}
\end{equation}
Write $\psi(s)=\psi(t)+(\psi(s)-\psi(t))$.
The terms carrying the constant $\psi(t)$ add up exactly to $\frac{\psi(t)}{\Gamma(\alpha)\tau}\cdot\frac{(t+\tau)^\alpha-t^\alpha}{\alpha}$,
since $\int_0^t((t+\tau-s)^{\alpha-1}-(t-s)^{\alpha-1})\mathrm{d}s+\int_t^{t+\tau}(t+\tau-s)^{\alpha-1}\mathrm{d}s
=\frac{(t+\tau)^\alpha-t^\alpha}{\alpha}$, the $\tau^\alpha$-terms cancelling. %
Split the first integral in the remainder at $s=t-\tau$:
on $[0,t-\tau]$ the integrand is $\mathcal{O}(\tau)(t-s)^{\alpha-2}$
times the bounded function $\psi(s)-\psi(t)$, and dominated convergence
 (the dominating function $(t-s)^{\alpha-2}|\psi(s)-\psi(t)|$ is integrable on $[0,t)$
 because $\beta>1-\alpha$) yields the limit $(\alpha-1)\int_0^t(t-s)^{\alpha-2}(\psi(s)-\psi(t))\mathrm{d}s$.
On $[t-\tau,t]$ and $[t,t+\tau]$ the local H\"older bound $|\psi(s)-\psi(t)|\le C|s-t|^\beta$
gives $\mathcal{O}(\tau^{\alpha+\beta})$, hence $\mathcal{O}(\tau^{\alpha+\beta-1})$
after division by $\tau$, which vanishes as $\tau\to0^+$.
Passing to the limit $\tau\to0^+$ (the same computation with $\tau<0$, $t+\tau>0$,
gives the left derivative) yields
\begin{equation}\label{V-der-B1}
V'(t)=-\frac{1}{\Gamma(\alpha)}\Big[h(V(t))\,t^{\alpha-1}+(\alpha-1)\int_0^t(t-s)^{\alpha-2}\big(h(V(s))-h(V(t))\big)\mathrm{d}s\Big],\qquad t>0,
\end{equation}
where the integral is absolutely convergent because $\psi=h(V)\in C^{0,\beta}_{\rm loc}(0,T]$ with $\beta>1-\alpha$.
By Lemma \ref{lem-weakly-singular-cont}, applied on any interval $(0,T]$
with $\psi=h(V)\in C[0,T]\cap C^{0,\beta}_{\rm loc}(0,T]$,
the map $t\mapsto\int_0^t(t-s)^{\alpha-2}(\psi(s)-\psi(t))\mathrm{d}s$ is continuous on $(0,\infty)$.
This, together with \dref{V-der-B1} and the continuity of $h(V(t))t^{\alpha-1}$ on $(0,\infty)$, gives $V\in C^1(0,\infty)$.
Consequently, by the chain rule, $W(t)=V(t^{1/\alpha})$ is differentiable on $(0,\infty)$ with
\begin{equation}\label{W-der-B1}
W'(t)=V'(t^{1/\alpha})\,\frac{1}{\alpha}\,t^{\frac{1}{\alpha}-1},\qquad t>0.
\end{equation}
By using $u^\alpha=t$ for $u=t^{1/\alpha}$ and the continuity of $h(V(t))$ at $t=0$
with $h(V(0))=h(A)$, 
we obtain
$$
\begin{array}{l}\disp
\frac{W(t)-W(0)}t=-\frac{1}{\Gamma(\alpha)\,t}\int_0^{t^{1/\alpha}}(t^{1/\alpha}-s)^{\alpha-1}h(V(s))\mathrm{d}s\crr\disp
\hspace{2.4cm}=-\frac{1}{\Gamma(\alpha)\,u^\alpha}\int_0^u(u-s)^{\alpha-1}h(V(s))\mathrm{d}s
\longrightarrow-\frac{h(A)}{\Gamma(\alpha+1)},\ \ \mbox{ as }t\to0^+.
\end{array}
$$
Hence $W$ is differentiable at $t=0$ with $W'_+(0)=-h(A)/\Gamma(\alpha+1)\le0$.
Finally, since $V$ is decreasing on $[0,\infty)$ (established in Step 5) and $V\in C^1(0,\infty)$,
we have $V'(t)\le0$ for $t>0$, and \dref{W-der-B1} together with $W'_+(0)\le0$ shows $W'(t)\le0$ for all $t\ge0$.
It remains to verify $\lim_{t\to0^+}W'(t)=W'_+(0)$. Set
$I(t):=\int_0^t(t-s)^{\alpha-2}(\psi(s)-\psi(t))\mathrm ds$.
From \dref{V-der-B1} and \dref{W-der-B1}, we have
$$
W'(t^\alpha)=-\frac{\psi(t)}{\Gamma(\alpha+1)}
+\frac{1-\alpha}{\alpha\Gamma(\alpha)}\,t^{1-\alpha}I(t),\qquad t>0.
$$
Since $\psi(t)\to \psi(0)=h(A)$ as $t\to0^+$, it suffices to show
$t^{1-\alpha}I(t)\to0$.
 Since $V$ and $\psi$ are decreasing and $h\geq0$,   we have for $0\leq s<t$,
\begin{equation}\label{eq:V-t-s-alpha}
\begin{array}{l}\disp
0\leq V(s)-V(t)=\frac{1}{\Gamma(\alpha)}
\bigg[
\int_0^s((t-r)^{\alpha-1}-(s-r)^{\alpha-1})\psi(r)\mathrm dr
+\int_s^t(t-r)^{\alpha-1}\psi(r)\mathrm dr
\bigg]\crr\disp\hspace{2.8cm}
\leq\frac{\psi(s)}{\Gamma(\alpha)}
\bigg[
\int_0^s((t-r)^{\alpha-1}-(s-r)^{\alpha-1})\mathrm dr
+\int_s^t(t-r)^{\alpha-1}\mathrm dr
\bigg]\crr\disp\hspace{2.8cm}
\leq  \frac{\psi(s)}{\Gamma(\alpha+1)}(t^\alpha-s^\alpha)
\leq \frac{h(A)}{\Gamma(\alpha+1)}(t^\alpha-s^\alpha).
\end{array}
\end{equation}
By the local Lipschitz continuity of $h$ on $[\delta_T,A]$, with Lipschitz constant $L_{T,A}$,
$
|\psi(s)-\psi(t)|=|h(V(s))-h(V(t))|
\le L_{T,A}|V(s)-V(t)|
\le C_7(t^\alpha-s^\alpha),
$
where $C_7=L_{T,A}h(A)/\Gamma(\alpha$ $+1)$. Substituting this into $I(t)$
and writing $s=tr$, we obtain
$$
|I(t)|
\le C_7\int_0^t(t-s)^{\alpha-2}(t^\alpha-s^\alpha)\,\mathrm ds
= C_7t^{2\alpha-1}
  \int_0^1(1-r)^{\alpha-2}(1-r^\alpha)\,\mathrm dr.
$$
The last integral is finite: near $r=1$, $1-r^\alpha\sim\alpha(1-r)$,
so the integrand behaves like $(1-r)^{\alpha-1}$, which is integrable
because $\alpha>0$. Near $r=0$ it is bounded. Hence
$
t^{1-\alpha}|I(t)|
\le C_7K_\alpha\,t^{\alpha}\longrightarrow0
$ as $t\to0^+$,
where $K_\alpha=\int_0^1(1-r)^{\alpha-2}(1-r^\alpha)\,\mathrm dr<+\infty$.
Therefore $W'(t)\to -\psi(0)/\Gamma(\alpha+1)=W'_+(0)$, so $W'$ is
continuous at $0$ and $W\in C^1[0,+\infty)$.
This proves (ii).\\
{\bf Step 7.} From $V(t)>0$ (established in Step 4) and the fact $V$ is decreasing on $[0,\infty)$ (established in Step 5),
we know that
$l:=\lim_{t\to\infty}V(t)$ exists and $l\geq0$.
Now, we claim that $l=0$. Indeed, if $l>0$, noting that $V(t)$ is monotonically decreasing,
we have $V(t)\geq l$ for $t\geq0$.
Since $V$ satisfies the integral equation \dref{Dh-Inequ-Vh-hinf}, we obtain
\begin{equation*}
V(t)
\leq A-\frac{h(l)}{\Gamma(\alpha)}
\int_0^{t}(t-s)^{\alpha-1}\mathrm{d}s=A-\frac{h(l)}{\Gamma(\alpha+1)}t^{\alpha},
\end{equation*}
which leads to $V(t)<0$ if $t>(A \Gamma(\alpha+1)/h(l))^{1/\alpha}$.
This is a contradiction.
Hence $l=0$ and 
we derive $\lim\limits_{t\rightarrow \infty} V(t)=0$.  This proves (iii).\\

\begin{remark}\label{rem-exist-uniq-Thm41}
As a by-product of the proof of Theorem \ref{Thm-Dhinfatzero-mono-asympt},
for each fixed $A>0$ the integral equation \dref{Dh-Inequ-Vh-hinf} admits a unique solution in $C[0,\infty)$.
The existence is established in Step 3: the diagonal limit $\widetilde{V}_A$ of the sequence $\{V^{\epsilon_n}\}$
satisfies the equation and is strictly positive, so that $h(|\widetilde{V}_A|)=h(\widetilde{V}_A)$.
The uniqueness follows in Step 5, where the Gronwall-type argument applied to $Z:=V-\widetilde{V}_A$
yields $V=\widetilde{V}_A$ for every (continuous) solution $V$. Thus the solution described in (i)--(iii) is unique.
\end{remark}

\begin{remark}\label{rem-W-reg-not-inherited}
The differentiability of $W(t)=V(t^{1/\alpha})$ asserted in (ii) is not inherited from the $\epsilon$-approximate level.
Step 2 shows, for each fixed $\epsilon>0$, that $W^{\epsilon}\in C^1[0,T]$ and $(W^{\epsilon})'=-X^{\epsilon}$,
where $W^{\epsilon}(t):=V^{\epsilon}(t^{1/\alpha})$ is the uniform limit of $W_n$ in $C[0,T]$.
The passage $\epsilon\to0$ performed in Steps 3--4, however, concerns the functions $V^{\epsilon}$ themselves:
it yields the uniform convergence $V^{\epsilon_n}\to V=\widetilde{V}_A$ on compact intervals,
but no statement about the derivatives $X^{\epsilon}$
(neither their convergence, nor the matching of any limit with the derivative of $V$) is established there.
In particular, the $\epsilon$-wise identity $(W^{\epsilon})'=-X^{\epsilon}$ cannot be passed to the limit directly.
We also note that $X^{\epsilon}$ involves the 
 factor $t^{\frac1\alpha-1}$ and depends on $\epsilon$
through $H_\epsilon'$, so any derivative-level limit argument would require additional justification.
This is why the regularity of $V$ itself (which satisfies the $\epsilon$-independent equation
\dref{Dh-Inequ-Vh-hinf} and is strictly positive on $[0,\infty)$,
so that $h'$ is Lipschitz on the range of $V$ on compact intervals) is established separately in Step 6 above.
\end{remark}

\begin{remark}\label{RemarkCaputoDdoesnotimplymono}
It is seen from the proof of Theorem \ref{Thm-Dhinfatzero-mono-asympt} that
in order to show $V(t)$ converges to zero, we spend a lot of space proving the monotonicity of the function $V(t)$ with respect to $t$.
It should be pointed out that
${}^C_0D_t^{\alpha}V(t)\leq0$ generally does not imply that $V(t)$ is non-increasing,
which is different from the simple fact that $V'(t)\leq0$ implies $V(t)$ is non-increasing.
To observe the above phenomena, we construct a continuously differentiable
function $V(t)$ satisfying ${}^C_0D_t^{\alpha}V(t)\leq0$ on $[0,+\infty)$ and
being not monotonic on $[0,+\infty)$.
For any set $E$, let $\chi_{E}(t)$ be the characteristic function of the set $E$.
Define $f(t)=\exp\big(-\frac{1}{(1-t)^2}\big)\chi_{[0,1)}(t)$.
Let $g_{k,N}(t)$ be defined by
$$g_{k,N}(t)=k\exp\left(-\frac{N}{(1-t)^2} - Nt^2   \right)\chi_{(1,+\infty)}(t),$$
where $k>0$ is a small real number and $N\geq 1$ is a sufficiently large real number. 
Then both $f$ and $g_{k,N}$ are smooth functions on $(0,\infty)$.
Construct $V$ as
\begin{equation*}
V(t)=\left\{\begin{array}{ll}\big( 2-( 1+t )^{\alpha}  \big)f(t) + g_{k,N}(t),\  &t\in[0,1)\cup(1,+\infty),\crr
0,\ & t=1,
\end{array}\right.
\end{equation*}
We claim that for small $k$ and large $N$ there holds
\begin{equation}\label{negative-Caputo-derivative}
\frac{1}{\Gamma(1-\alpha)}
\int_0^t(t-s)^{-\alpha} V'(s)\mathrm{d}s<0,\ \ \ t>0.
\end{equation}
For $t\in (0,1)$, $V'(t)=[-\alpha(1+t)^{\alpha-1}
+[ 2-( 1+t )^{\alpha}][2(t-1)^{-3}]]\exp\big(-\frac{1}{(1-t)^2}\big)< 0$, and clearly $V'(1)=0$,
so we can see that $\int_0^t(t-s)^{-\alpha} V'(s)\mathrm{d}s<0$ for all $t\in(0,1]$.
For $t> 1$, we find that
\begin{equation}\label{CapDV-tge1}
\begin{array}{l}\disp
\frac{1}{\Gamma(1-\alpha)}
\int_0^1(t-s)^{-\alpha} V'(s)\mathrm{d}s
\leq
\frac{1}{\Gamma(1-\alpha)}\int_0^1(t-s)^{-\alpha} (-\alpha)(1+s)^{\alpha-1}f(s)\mathrm{d}s\crr\disp
\leq
\frac{1}{\Gamma(1-\alpha)}\int_0^{1/2}(t-s)^{-\alpha} (-\alpha)(1+s)^{\alpha-1}f(s)\mathrm{d}s\crr\disp
\leq
\left(\frac32\right)^{\alpha-1}\frac{-\alpha e^{-4}}{\Gamma(1-\alpha)}
\int_0^{1/2} (t-s)^{-\alpha}\mathrm{d}s
\leq
\left(\frac32\right)^{\alpha-1}\frac{-\alpha e^{-4}}{\Gamma(1-\alpha)} \frac{t^{-\alpha}}{2}:=
-C_{\alpha} t^{-\alpha},
\end{array}
\end{equation}
where 
$f'(t)\leq 0$ and $f(t)\geq e^{-4}$ for $0<t<1/2$
are used. For $t>1$, a simple calculation gives
\begin{equation}\label{firstDV-tge1}V'(t)=k( -2N(1-t)^{-3}- 2Nt)\exp(-N(1-t)^{-2}-Nt^2).\end{equation}
Clearly, $\lim_{t\to1^+}( -2N(1-t)^{-3}- 2Nt)\exp(-N(1-t)^{-2}-Nt^2/2)=
\lim_{t\to+\infty}( -2N(1-t)^{-3}- 2Nt)\exp(-N(1-t)^{-2}-Nt^2/2)=0$.
This implies that there exists an absolute constant $C_0$ such that $|V'(t)|\leq kC_0\exp(-Nt^2/2)$ for all $t>1$.
 Hence it follows that for $t\in(1,2]$,
$\int_1^t(t-s)^{-\alpha} |V'(s)|\mathrm{d}s
\leq \frac{kC_0e^{-N/2}}{1-\alpha}\cdot(t-1)^{1-\alpha}
\leq \frac{kC_0e^{-N/2}}{1-\alpha}2^{\alpha}t^{-\alpha}=a_{k,N}t^{-\alpha} $
and for $t>2$,
\begin{equation}\label{CapDV-tge2}
\begin{array}{l}\disp
\int_1^t(t-s)^{-\alpha} |V'(s)|\mathrm{d}s\crr\disp
\leq
kC_0\int_1^{t/2}(t-s)^{-\alpha}\exp(-Ns^2/2)\mathrm{d}s
+
kC_0\int_{t/2}^{t}(t-s)^{-\alpha} \exp(-Ns^2/2)\mathrm{d}s
\crr\disp
\leq
kC_0 \left(\frac{t}{2}\right)^{-\alpha}\int_1^{t/2}\exp(-Ns^2/2)\mathrm{d}s
+
kC_0 \exp(-Nt^2/8)
\int_{t/2}^{t}(t-s)^{-\alpha} \mathrm{d}s\crr\disp
\leq
b_{k,N}t^{-\alpha}+c_{k,N}t^{-\alpha},
\end{array}
\end{equation}
where
$a_{k,N}={kC_0e^{-N/2}}2^{\alpha}/(1-\alpha)$,
$b_{k,N}=kC_02^{\alpha}\int_1^{\infty}\exp(-Ns^2/2)\mathrm{d}s$
and $c_{k,N}=kC_02^{\alpha-1}C_1$ with $C_1=\sup$ $\{t\exp(-Nt^2/8)/(1-\alpha):t\geq1\}$. It is clear that
$\lim_{k\to0+,N\to+\infty}a_{k,N}=\lim_{k\to0+,N\to+\infty}b_{k,N}
=\lim_{k\to0+,N\to+\infty}c_{k,N}=0$.
Combining this with \dref{CapDV-tge1} and \dref{CapDV-tge2}, we see that 
 $\int_0^t(t-s)^{-\alpha} V'(s)\mathrm{d}s<0$ for all $t>1$
is true for small $k$ and large $N$. Hence,
${}^C_0D_t^{\alpha}V(t)<0$ on $(0,+\infty)$.
It is obvious that $V(t)$ is not monotonic on $[0,+\infty)$ due to $V(t)>0$ for $t\in[0,1)\cup(1,+\infty)$ and $V(1)=0$.
Moreover, from \dref{firstDV-tge1}, there exists two constants $t_1,t_2>0$
such that $V'(t)>0$ for $t\in(1,t_1)$ and $V'(t)<0$ for $t\in(t_2,+\infty)$.
\end{remark}

\begin{lemma}\label{lem-fracsys-compare}
Let $f(t,v)$ be a continuous function, 
and let $f(t,v)$ be a nondecreasing function with respect to $v$.
Suppose $v(t)$ is a continuous function satisfying
\begin{equation}\label{comthm-Vineq1}
{}^C_0D_t^{\alpha}v(t)\leq \lambda v(t)+f(t,v(t)),\; v(0)=v_0.
\end{equation} Suppose
also that $w(t)$ is a continuous function satisfying
\begin{equation}\label{comthm-Vineq2}
{}^C_0D_t^{\alpha}w(t)\geq  \lambda w(t)+ f(t,w(t)),\; w(0)=w_0.
\end{equation}
Then $v(t)< w(t)$ on their common interval of existence if $v_0<w_0$.
\end{lemma}
\noindent{\bf Proof.}
By \dref{comthm-Vineq1}, there exists a non-negative function $m_1(t)\geq0$ such that
${}^C_0D_t^{\alpha}v(t)=\lambda v(t)+ f(t,v(t))-m_1(t)$,
 which, by \cite[Page 323, Theorem 5.15]{KilbasJbook2006}, is equivalent to the following Volterra integral equation:
\begin{equation}
v(t)=E_{\alpha}(\lambda t^{\alpha})v(0)
   +\int_0^t(t-s)^{\alpha-1}E_{\alpha,\alpha}(\lambda (t-s)^{\alpha})[f(s,v(s))-m_1(s)]\mathrm{d}s.
\end{equation}
Analogously, there exists a non-negative function $m_2(t)\geq0$ such that
\begin{equation}
w(t)=E_{\alpha}(\lambda t^{\alpha})w(0)+\int_0^t(t-s)^{\alpha-1}E_{\alpha,\alpha}(\lambda (t-s)^{\alpha})[f(s,w(s))+m_2(s)]\mathrm{d}s.
\end{equation}
Let $z(t)=w(t)-v(t)$ with $z(0)=w_0-v_0>0$. Clearly, $z(t)$ satisfies
\begin{equation}
z(t)=E_{\alpha}(\lambda t^{\alpha})z(0)
+\int_0^t(t-s)^{\alpha-1}E_{\alpha,\alpha}(\lambda (t-s)^{\alpha})[f(s,w(s))-f(s,v(s))+m_3(s)]\mathrm{d}s,
\end{equation}
where $m_3(t)=m_1(t)+m_2(t)\geq0$.
We show $z(t)>0$ for all $t\geq0$. Indeed, if this is not true, there exists $t_1>0$ such that
$z(t)>0$ for $t\in[0,t_1)$ and $z(t_1)=0$.
Since $f(t,v)$ is nondecreasing function with respect to $v$, by $z(t)=w(t)-v(t)>0$ for $t\in[0,t_1)$, we have $g(s):=f(s,w(s))-f(s,v(s))\geq0$ for $s\in[0,t_1)$.
It follows from $E_{\alpha}(\lambda t_1^{\alpha})z(0)>0$ that
$$
\begin{array}{l}\disp
0=z(t_1)=E_{\alpha}(\lambda t_1^{\alpha})z(0) +\int_0^{t_1}(t_1-s)^{\alpha-1}E_{\alpha,\alpha}(\lambda (t_1-s)^{\alpha})[g(s)+m_3(s)]\mathrm{d}s>0,
\end{array}
$$
which completes the proof.
\hfill$\Box$

\begin{theorem}\label{thm-CapDVleq-hV-hinf}
Let $h$ be a class-$\mathcal{K}$ function and let its derivative function $h'$
satisfy local Lipschitz continuity on $(0,+\infty)$.  
Let $V(t)$ be a non-negative function on $[0,+\infty)$
and satisfy ${}^C_0D_t^{\alpha}V(t)\leq-h(V(t))$. Then $V(t)\to0$ as $t\to\infty$.
\end{theorem}
\noindent{\bf Proof.}
First, consider $V(0)>0$.
By Theorem \ref{Thm-Dhinfatzero-mono-asympt} and Remark \ref{rem-exist-uniq-Thm41},
we know that ${}^C_0D_t^{\alpha}W(t)=-h(W(t))$ with $W(0)=2V(0)$ admits a unique solution $W(t)\in C[0,+\infty)$
satisfying $W(t)>0$ for all $t\in[0,+\infty)$ and  $\lim_{t\to\infty}W(t)=0$.

By ${}^C_0D_t^{\alpha}V(t)\leq-h(V(t))$,
 there exists a non-negative function $m(t)$ such that ${}^C_0D_t^{\alpha}V(t)=-h(V(t))-m(t)$. Then
let $Z(t)=W(t)-V(t)$. By Theorem \ref{Thm-Dhinfatzero-mono-asympt},  $W(t)>0$ for all $t\geq0$.
Now we claim  $Z(t)>0$ for all $t\geq0$.
Indeed, if not, there exists a $t_1>0$ such that $Z(t_1)=0$ and $Z(t)>0$ for all $t\in[0,t_1)$.
Clearly, $Z(t)$ satisfies
\begin{equation}\label{hinf-Z-defequ-S}
{}^C_0D_t^{\alpha}Z(t)=-[h(W(t))-h(V(t))]+m(t)
=-\eta(t)Z(t)+m(t),
\end{equation}
where $\eta(t)$ is defined by $\eta(t)=\frac{h(W(t))-h(V(t))}{W(t)-V(t)}$
 for $t\in[0,t_1)$ and $\eta(t)=h'(W(t_1))$ for $t=t_1$.
Note that $W(t)-V(t)>0$ for $t\in[0,t_1)$, the function $\eta(t)$ is well-defined.
Clearly, $\lim_{t\to t_1^-}\eta(t)=h'(W(t_1))$, which implies that $\eta(t)$ is continuous on $t\in[0,t_1]$.
Hence, there exists a constant $K>0$ such that $|\eta(t)|\leq K$ for all $t\in[0,t_1]$.  \dref{hinf-Z-defequ-S} gives
\begin{equation}\label{hinf-Z-defequ-I1-S}
{}^C_0D_t^{\alpha}Z(t)
=-2KZ(t)+(2K-\eta(t))Z(t)+m(t)
\geq -2KZ(t)+(2K-\eta(t))Z(t).
\end{equation}
It is seen that $(2K-\eta(t))x$ is increasing with respect to $x\geq0$ when $t\in[0,t_1]$.
Apply Lemma \ref{lem-fracsys-compare} with $\lambda=-2K$ and $f(t,x):=(2K-\eta(t))x$
(nondecreasing in $x\ge0$ because $\eta(t)\le K$), to $v\equiv0$ (which satisfies \dref{hinf-Z-defequ-I1-S} with equality)
and $w=Z$, since $v_0=0<Z(0)=V(0)$,
we know that $Z(t)>0$ for all $t\in[0,t_1]$ and thus $Z(t_1)>0$,
which is a contradiction. Hence, $W(t)>V(t)\geq0$ for all $t\geq0$.
It follows from Theorem \ref{Thm-Dhinfatzero-mono-asympt} that $\lim_{t\to\infty}W(t)=0$ and thus
$\lim_{t\to\infty}V(t)=0$.
Finally, we consider $V(0)=0$. From ${}^C_0D_t^{\alpha}V(t)\leq0$ and the non-negativity of $V(t)$,
 it is trivial to see that $V(t)\leq V(0)$ and thus $V(t)\equiv0$ for $t\geq0$.
This completes the proof.
\hfill$\Box$

\begin{lemma}\label{double-h1h2-compresh-app}
Let $h$  be a class-$\mathcal{K}$ function.
Then there exist two  class-$\mathcal{K}$ $h_1$, $h_2$ such that $h_1(x)\leq h(x)\leq h_2(x)$ for all $x\geq0$
and $h_1,h_2\in C^2(0,+\infty)$.
\end{lemma}
\noindent{\bf Proof.} For $x=0$, we set $h_1(0)=h_2(0)=l_1(0)=l_2(0)=0$. For $x>0$ we define
$$
l_1(x):=\frac{1}{x}\int_0^x h(s)\,\mathrm{d}s,\quad h_1(x):=\frac{1}{x}\int_0^x l_1(s)\,\mathrm{d}s,
$$
and
$$
l_2(x):=\frac{1}{x}\int_x^{2x} h(s)\,\mathrm{d}s,\quad h_2(x):=\frac{1}{x}\int_x^{2x} l_2(s)\,\mathrm{d}s,
$$
A simple computation shows that
$h_1(x)\leq l_1(x)\leq h(x)\leq l_2(x)\leq h_2(x)$ for all $x\ge0$,
and $h_1,h_2$ are class-$\mathcal{K}$ functions which are $C^2(0,\infty)$ while $l_1,l_2\in C^1(0,\infty)$.
\hfill$\Box$

\begin{theorem}\label{thm-DV-hV-gencase}
Let $h$ be a class-$\mathcal{K}$ function. Let $V$ be a non-negative function and satisfy
 ${}^C_0D_t^{\alpha}V(t)\leq-h(V(t))$. Then $V(t)\to0$ as $t\to\infty$.
\end{theorem}
\noindent{\bf Proof.}
By Lemma \ref{double-h1h2-compresh-app}, there exists a class-$\mathcal{K}$ function $h_1$
such that $h_1\in C^2(0,+\infty)$ and $h_1(x)\leq h(x)$ for $x\geq0$.
Hence, we have ${}^C_0D_t^{\alpha}V(t)\leq-h(V(t))\leq-h_1(V(t))$.
By Theorem \ref{thm-CapDVleq-hV-hinf}
we derive $V(t)\to0$ as $t\to+\infty$.
\hfill$\Box$

With the preparations above in place, we now provide a rigorous proof of Theorem 11 in \cite{LiYanCP2009auto}.
\begin{theorem}\label{thm-Gen-fl} (Theorem 11 of \cite{LiYanCP2009auto}) 
Let $x = 0$ be an equilibrium point for the system \dref{sys-plantFrac} and $\Omega$ be a domain containing the origin.
Assume that there exist a Lyapunov function $V(t, x(t))$ and three class-$\mathcal{K}$ functions $\alpha_i$, $i = 1, 2, 3$
satisfying
\begin{equation}\label{tm-fv-Frac-Ly-cond-R}
\begin{array}{l}
\alpha_1(\|x\|)\leq V(t, x(t))
\leq\alpha_2(\|x\|),\\
{}^C_0D_t^{\alpha}V(t, x(t))\leq-\alpha_3(\|x\|)
\end{array}
\end{equation}
Then the equilibrium point of system \dref{sys-plantFrac} is asymptotically stable.
\end{theorem}
\noindent{\bf Proof.}
\noindent{\bf Step 1:} We show the Lyapunov stability.
Set $\alpha_4(s):=\alpha_3(\alpha_2^{-1}(s))$, which is a
class-$\mathcal K$ function. From $\|x(t)\|\geq\alpha_2^{-1}(V(t,x(t)))$, we have
$$
{}^C_0D_t^\alpha V(t,x(t))
\le-\alpha_3(\|x(t)\|)
\le-\alpha_3(\alpha_2^{-1}(V(t,x(t))))
=-\alpha_4(V(t,x(t))).
$$
In particular, ${}^C_0D_t^\alpha V(t,x(t))\leq0$. By the fundamental
identity for the Caputo derivative,
$$
V(t,x(t))=V(0,x(0))+{}_0I_t^\alpha\big[{}^C_0D_\tau^\alpha V(\tau,x(\tau))\big].
$$
Since the integrand is non-positive and ${}_0I_t^\alpha$ is a positive
operator, the fractional integral is non-positive, and therefore
$V(t,x(t))\le V(0,x(0))$ for all $t\ge0$.
Combining with the upper and lower bounds in \dref{tm-fv-Frac-Ly-cond-R}, we obtain
$$
\alpha_1(\|x(t)\|)\le V(t,x(t))\le V(0,x(0))\le\alpha_2(\|x(0)\|).
$$
Given any $\varepsilon>0$, choose $\delta>0$ such that
$\alpha_2(\delta)<\alpha_1(\varepsilon)$ (possible because
$\alpha_1,\alpha_2\in\mathcal K$ are continuous and strictly increasing
with $\alpha_1(0)=\alpha_2(0)=0$). Then $\|x(0)\|<\delta$ implies
$$
\alpha_1(\|x(t)\|)\le\alpha_2(\|x(0)\|)<\alpha_1(\varepsilon),
$$
and since $\alpha_1$ is strictly increasing, $\|x(t)\|<\varepsilon$
for all $t\ge0$. This proves Lyapunov stability.

\noindent{\bf Step 2:} We show the attractivity.
It follows as a by-product of Theorem \ref{thm-DV-hV-gencase}, that is, of the
analysis of the Volterra equation \eqref{eq:main-VI} carried out in this
section. A second, independent proof, which uses none of the properties of
the solutions of \eqref{eq:main-VI}, is given in Remark \ref{rem-alt-proof-Thm34};
both routes are kept, since the analysis of \eqref{eq:main-VI} is a
contribution of independent interest.
\noindent From ${}^C_0D_t^\alpha V(t,x(t))\le-\alpha_4(V(t,x(t)))$ and
Theorem~\ref{thm-DV-hV-gencase} (which states that
${}^C_0D_t^\alpha V\le-h(V)$ with $h$ being class-$\mathcal K$ function implies
$V(t)\to0$), we obtain
$\lim_{t\to+\infty}V(t,x(t))=0$.
Using the lower bound $\alpha_1(\|x(t)\|)\le V(t,x(t))$, we have
$\lim_{t\to+\infty}\alpha_1(\|x(t)\|)=0$.
Since $\alpha_1^{-1}$ is continuous at $0$ with $\alpha_1^{-1}(0)=0$, we derive
$
\lim_{t\to+\infty}\|x(t)\|=0.
$
This proves attractivity.
Hence, the equilibrium point of system \dref{sys-plantFrac} is asymptotically stable.
\hfill$\Box$

\begin{remark}\label{rem-alt-proof-Thm34}
The proof of Theorem \ref{thm-Gen-fl} obtains the attractivity in Step 2 as a
special case of Theorem \ref{thm-DV-hV-gencase}, hence as a by-product of the
analysis of the Volterra equation \eqref{eq:main-VI} carried out in Theorems
\ref{Thm-Dhinfatzero-mono-asympt} and \ref{thm-DV-hV-gencase}. That analysis is
a contribution of independent interest, pursued further in
Section \ref{Sec-h-classK-property}, and it is kept as the main route. We record
here a second, independent proof of the same step, which uses none of the
existence, positivity, monotonicity, uniqueness and regularity properties of the
solutions of \eqref{eq:main-VI}: it combines only the fundamental identity of
the Caputo derivative, the comparison lemma for linear fractional equations
(Lemma \ref{lem-fracsys-compare}) and the decay of the Mittag-Leffler function.
The second argument is shorter, and it is self-contained for a reader interested
in Theorem \ref{thm-Gen-fl} alone.

\smallskip
\noindent{\bf Second proof of the attractivity.}
Let $v(t):=V(t,x(t))$ and $v_0:=v(0)\ge0$. If $v_0=0$, then
$\alpha_1(\|x(0)\|)\le v_0=0$, so $x(0)=0$ and the equilibrium is already at
rest. Hence assume $v_0>0$. Since $v(t)\le\alpha_2(\|x(t)\|)$, monotonicity of
$\alpha_2^{-1}$ and $\alpha_3$ gives $\alpha_2^{-1}(v(t))\le\|x(t)\|$, and
therefore
\begin{equation}\label{eps-scalar-ineq}
{}^C_0D_t^{\alpha}v(t)\le-\alpha_3(\|x(t)\|)
\le-\alpha_3(\alpha_2^{-1}(v(t)))=-\alpha_4(v(t)).
\end{equation}
It remains to show that \dref{eps-scalar-ineq} forces
$\lim_{t\to+\infty}v(t)=0$. Since $\alpha_4=\alpha_3(\alpha_2^{-1}(s))\ge0$,
the fundamental identity for the Caputo derivative yields
\begin{equation}\label{eps-upper-v0}
v(t)=v_0+(\Ial\,{}^C_0D_t^{\alpha}v)(t)\le v_0.
\end{equation}
For fixed $\varepsilon\in(0,v_0)$, set
$c_\varepsilon:=\alpha_4(\varepsilon)/v_0>0$. We claim that
\begin{equation}\label{eps-linear-bound}
\alpha_4(z)\ge c_\varepsilon(z-\varepsilon),\qquad 0\le z\le v_0.
\end{equation}
Indeed, if $0\le z\le\varepsilon$, then $(z-\varepsilon)/v_0\le0$ and the
right-hand side is non-positive, while $\alpha_4\ge0$. If
$\varepsilon<z\le v_0$, monotonicity of $\alpha_4$ gives
$\alpha_4(z)\ge\alpha_4(\varepsilon)$, while $(z-\varepsilon)/v_0\le1$, so that
$c_\varepsilon(z-\varepsilon)\le\alpha_4(\varepsilon)\le\alpha_4(z)$.

Since the Caputo derivative of the constant $\varepsilon$ vanishes,
\dref{eps-scalar-ineq} and \dref{eps-linear-bound} give
${}^C_0D_t^{\alpha}(v-\varepsilon)(t)\le-c_\varepsilon(v(t)-\varepsilon)$.
Setting $u(t):=v(t)-\varepsilon$, we have
\begin{equation}\label{eps-comparison-ineq}
{}^C_0D_t^{\alpha}u(t)\le-c_\varepsilon u(t).
\end{equation}
Apply Lemma \ref{lem-fracsys-compare} with $\lambda=-c_\varepsilon$ and
$f\equiv0$ to the linear equation
${}^C_0D_t^{\alpha}w(t)=-c_\varepsilon w(t)$, $w(0)=v_0$, whose solution is
$w(t)=v_0E_\alpha(-c_\varepsilon t^\alpha)$. Since $u(0)=v_0-\varepsilon<v_0$,
the comparison principle gives $u(t)<v_0E_\alpha(-c_\varepsilon t^\alpha)$, and
hence
\begin{equation}\label{eps-ML-estimate}
v(t)<\varepsilon+v_0E_\alpha(-c_\varepsilon t^\alpha).
\end{equation}
By \dref{4t1-6-ne}, $E_\alpha(-c_\varepsilon t^\alpha)\to0$ as $t\to\infty$,
so that $\limsup_{t\to\infty}v(t)\le\varepsilon$. Since $v(t)\ge0$ and
$\varepsilon\in(0,v_0)$ is arbitrary, $\lim_{t\to\infty}v(t)=0$. Finally,
$\alpha_1(\|x(t)\|)\le v(t)$ gives
$\|x(t)\|\le\alpha_1^{-1}(v(t))\to0$ as $t\to\infty$, which is the attractivity
of the equilibrium.
\end{remark}

\begin{lemma}\label{lem:Vx-Caputo-ine-min-HC}
Let $V:\mathbb{R}^n\to\mathbb{R}^+$ be a finite convex and differentiable function on $\mathbb{R}^n$.
Let $x\in C[0,T]\cap C^1(0,T]$ satisfy $x'(t)=\mathcal{O}(t^{\alpha-1})$ as $t\to0^+$, and suppose that the limit
$\xi_0:=\lim_{t\to0^+}{}^C_0D_t^{\alpha}x(t)\in\mathbb{R}^n
$
exists.  Then, with the convention
${}^C_0D_0^{\alpha}f:=\lim_{t\to0^+}{}^C_0D_t^{\alpha}f$, we have
\begin{equation}\label{eq:Vx-Caputo-ine-min}
{}^C_0D_t^{\alpha}V(x(t))\leq
\langle\nabla V(x(t)),{}^C_0D_t^{\alpha}x(t)\rangle,
\quad \forall t\in[0,T],
\end{equation}
and ${}^C_0D_0^{\alpha}V(x(0))=\langle\nabla V(x(0)),\xi_0\rangle$.
In particular \dref{eq:Vx-Caputo-ine-min} holds with equality at $t=0$.
\end{lemma}
\noindent{\bf Proof.}
We first show that \dref{eq:Vx-Caputo-ine-min} holds for $t\in(0,T]$.
Fix $t\in(0,T]$.
Since $V:\mathbb{R}^n\to\mathbb{R}^+$ is a finite convex and  differentiable,
by \cite[Page 246, Corollary 25.5.1]{bookConvexanalysis1970},
 $\nabla V(x)$ is continuous on $\mathbb{R}^n$, and by $x\in C[0,T]$, $\nabla V(x(s))$ is bounded on $[0,T]$.
Set $C_0:=\sup_{u\in[0,T]}\|\nabla V(x(u))\|<+\infty$.
Moreover $x\in C[0,T]\cap C^1(0,T]$
with $x'(u)=\mathcal{O}(u^{\alpha-1})$ is absolutely continuous on $[0,T]$.
Consequently $V(x(u))$ is absolutely continuous with respect to $u$ and
$\frac{\dd}{\dd u}V(x(u))=\langle\nabla V(x(u)),x'(u)\rangle$ on $(0,T]$.
Since $x\in C[0,T]\cap C^1(0,T]$ satisfy $x'(s)=\mathcal{O}(s^{\alpha-1})$, for $t>0$, we have
\begin{equation}\label{eq:Vx-Caputo-x-i-N}
{}^C_0D_t^{\alpha}x(t)=\frac{1}{\Gamma(1-\alpha)}
\bigg[\frac{\partial}{\partial t}\int_0^t(t-s)^{-\alpha}x(s)
{\rm d}s-t^{-\alpha}x(0)\bigg]
=\frac{1}{\Gamma(1-\alpha)}\int_0^ts^{-\alpha}{x}'(t-s){\rm d}s
\end{equation}
and
\begin{equation}\label{eq:Vx-Caputo-Vx-i-N}
{}^C_0D_t^{\alpha}V(x(t))
=\frac{1}{\Gamma(1-\alpha)}\int_0^ts^{-\alpha}\langle\nabla V(x(t-s)),{x}'(t-s)\rangle{\rm d}s,
\end{equation}
which implies that \dref{eq:Vx-Caputo-ine-min} is equivalent to
\begin{equation}\label{eq:lem2-equiv}
I:=\int_0^t s^{-\alpha}\langle\nabla V(x(t))-\nabla V(x(t-s)),x'(t-s)\rangle\,{\rm d}s\geq0.
\end{equation}
Set $G(s)=V(x(t-s))-V(x(t))-\langle\nabla V(x(t)),x(t-s)-x(t)\rangle$ for $s\in[0,t]$.
Then $G\in C[0,t]\cap C^1(0,t)$, $G(0)=0$, and
$
G'(s)=\langle \nabla V(x(t))-\nabla V(x(t-s)),x'(t-s)\rangle,
$
hence \dref{eq:lem2-equiv} reads $I=\int_0^t s^{-\alpha}G'(s)\,{\rm d}s$.

Now we show the integrability of $s^{-\alpha}G'(s)$ on $(0,t)$.
$\|x'(u)\|\le C_1u^{\alpha-1}$ for $u\in(0,T]$ and for some $C_1>0$.
For $s\in(0,t/2]$, since
$t-s\ge t/2$, we have $\|x'(t-s)\|\le C_1(t/2)^{\alpha-1}$, so
$|s^{-\alpha}G'(s)|\le 2C_0C_1(t/2)^{\alpha-1}s^{-\alpha}\in L^1(0,t/2)$.  For
$s\in[t/2,t)$, $\|x'(t-s)\|\le C_1(t-s)^{\alpha-1}$ and $s^{-\alpha}\le(t/2)^{-\alpha}$,
so $|s^{-\alpha}G'(s)|\le 2C_0C_1(t/2)^{-\alpha}(t-s)^{\alpha-1}\in L^1(t/2,t)$
(because $\alpha-1>-1$).  Thus $s^{-\alpha}G'(s)$ is integrable on $(0,t)$ and $I$ is finite.

Integrating by parts on $[\varepsilon,t-\varepsilon]$, $0<\varepsilon<t/2$, gives
\begin{equation}\label{eq:Vcon-Caputo-equiv-II-N}
\int_{\varepsilon}^{t-\varepsilon}s^{-\alpha}G'(s){\rm d}s
=(t-\varepsilon)^{-\alpha}G(t-\varepsilon)-\varepsilon^{-\alpha}G(\varepsilon)
 +\alpha\int_{\varepsilon}^{t-\varepsilon}s^{-\alpha-1}G(s){\rm d}s
\end{equation}
Since $V$ is differentiable at $x(t)$, by the definition of differentiability of $V$,
$G(\varepsilon)=o(\|x(t-\varepsilon)-x(t)\|)$ as $\varepsilon\to0^+$.  For
$\varepsilon\le t/2$, 
$
\|x(t-\varepsilon)-x(t)\|
=\big\|\int_{t-\varepsilon}^t x'(u)\,{\rm d}u\big\|
\le C_1\int_{t-\varepsilon}^t u^{\alpha-1}\,{\rm d}u
\le C_1(t-\varepsilon)^{\alpha-1}\varepsilon
=\mathcal{O}(\varepsilon)
$ as $\varepsilon\to0$.
Hence $G(\varepsilon)=o(\varepsilon)$ as $\varepsilon\to0$, and since $\alpha<1$,
$\varepsilon^{-\alpha}G(\varepsilon)=o(\varepsilon^{1-\alpha})$ as $\varepsilon\to0$.
From $G(s)=o(s)$ as $s\to0^+$, it follows that
$s^{-\alpha-1}G(s)=o(s^{-\alpha})\in L^1(0,t)$.
Since $V$ is convex, by \cite[Page 242, Theorem 25.1]{bookConvexanalysis1970},
we have $G(s)\geq0$. Hence $\int_{\varepsilon}^{t-\varepsilon}s^{-\alpha-1}G(s){\rm d}s\geq0$ for all $\varepsilon\in(0,t/2)$
and the integrals converge as $\varepsilon\to0^+$.
  Letting
$\varepsilon\to0^+$ in \dref{eq:Vcon-Caputo-equiv-II-N}, we obtain
\[
I=t^{-\alpha}G(t)+\alpha\int_0^t s^{-\alpha-1}G(s)\,{\rm d}s\geq0,
\]
which is \dref{eq:lem2-equiv}.  Hence \dref{eq:Vx-Caputo-ine-min} holds for
every $t\in(0,T]$.

Next, we show that \dref{eq:Vx-Caputo-ine-min} is valid if $t=0$.
Since $\nabla V(x)$ is continuous on $\mathbb{R}^n$, and $x\in C[0,T]$, the function
$s\mapsto\nabla V(x(s))$ is bounded on $[0,T]$. We write
$\nabla V(x(s))=\nabla V(x(0))+\rho(s)$, where $\rho(s)\to0$ as $s\to0^+$
and $\|\rho(s)\|\leq 2C_0$ with $C_0=\sup_{s\in[0,T]}\|\nabla V(x(s))\|<+\infty$.
Since $\frac{d}{ds}(V(x(s)))=\langle\nabla V(x(s)),x'(s)\rangle$ on $(0,T]$, we have, as
$t\to0^+$,
\begin{equation}\label{eq:min-V}
\begin{array}{ll}
{}^C_0D_t^{\alpha}V(x(t))
&=\frac{1}{\Gamma(1-\alpha)}\int_0^t(t-s)^{-\alpha}
\langle\nabla V(x(0)),x'(s)\rangle\,{\rm d}s+J(t)\\
&
=\langle\nabla V(x(0)),{}^C_0D_t^{\alpha}x(t)\rangle+J(t),
\end{array}
\end{equation}
where $J(t)=\frac{1}{\Gamma(1-\alpha)}\int_0^t(t-s)^{-\alpha}
\langle\rho(s),x'(s)\rangle\,{\rm d}s$.
Letting $s=tu$ in $J(t)$ and using $\|x'(s)\|\leq C_1 s^{\alpha-1}$, 
it follows from $\|\rho(tu)\|\leq 2C_0$, $\|\rho(tu)\|\to0$ as $t\to0^+$ for $u\in(0,1)$,
 $(1-u)^{-\alpha}u^{\alpha-1}\in L^1(0,1)$ and the dominated convergence theorem that
$$
\Big|\int_0^t(t-s)^{-\alpha}\langle\rho(s),x'(s)\rangle\,{\rm d}s\Big|
\leq C_1\int_0^1(1-u)^{-\alpha}u^{\alpha-1}\|\rho(tu)\|\,{\rm d}u
\longrightarrow0
\quad(t\to0^+),
$$
Combining this with \dref{eq:min-V}, $\xi_0:=\lim_{t\to0^+}{}^C_0D_t^{\alpha}x(t)\in\mathbb{R}^n$
and the continuity of $\nabla V$ at $x(0)$ gives
\begin{equation}\label{eq:min-limits}
\lim_{t\to0^+}{}^C_0D_t^{\alpha}V(x(t))
=\lim_{t\to0^+}\langle\nabla V(x(t)),{}^C_0D_t^{\alpha}x(t)\rangle
=\langle\nabla V(x(0)),\xi_0\rangle.
\end{equation}
Since \dref{eq:Vx-Caputo-ine-min} holds for every $t>0$, passing to the limit
$t\to0^+$ (legitimate by \dref{eq:min-limits}) yields equality at $t=0$:
${}^C_0D_0^{\alpha}V(x(0))=\langle\nabla V(x(0)),{}^C_0D_0^{\alpha}x(0)\rangle$.
Hence \dref{eq:Vx-Caputo-ine-min} holds for all $t\in[0,T]$.
\hfill$\Box$

\begin{remark}\label{remark-Vx-fi}
The inequality \dref{eq:Vx-Caputo-ine-min} was first proved in \cite[Theorem 1]{ChenDaiSongZhangZQ2017iet}
where $x(t)$ is assumed to satisfy $x\in C^1[0,\infty)$.
However, in typical applications of the inequality \dref{eq:Vx-Caputo-ine-min},
the function $x$ is a solution to a fractional differential equation,
and this solution is not continuously differentiable at $t=0$.
Let us look at a simple example: Consider fractional differential equation
${}^C_0D_t^{\alpha}x(t)=-x(t)$ with initial value  $x(0)=1$.
Here $\alpha\in(0,1)$. Its solution can be explicitly found to be
$x(t)=E_{\alpha}(-t^{\alpha})=1-\frac{t^{\alpha}}{\Gamma(\alpha+1)}+\sum_{k=2}^{+\infty}\frac{(-1)^kt^{k\alpha}}{\Gamma(k\alpha+1)}$.
Clearly, this function $x(t)$ is not differentiable at $t=0$.
Hence, this inequality with the requirement $x\in C^1[0,\infty)$
cannot be applied directly to the solutions of fractional differential equations.
\end{remark}

\begin{remark}\label{remark-D-Vx-ineq-exam}
To ensure \dref{eq:Vx-Caputo-ine-min} holds for $t=0$, the assumption that the limit
$
\xi_0:=\lim_{t\to0^+}{}^C_0D_t^{\alpha}x(t)\in\mathbb{R}^n
$ cannot be removed.
Since  \dref{eq:Vx-Caputo-ine-min} holds for every $t>0$, it passes to
the limit at $t=0$ provided both limits exist.  However, under the
assumption alone that $x\in C[0,T]\cap C^1(0,T]$ satisfy $x'(t)=\mathcal{O}(t^{\alpha-1})$ as $t\to0^+$,
the limit $\lim_{t\to0^+}{}^C_0D_t^{\alpha}x(t)$ need not exist.
For example,  take $x(0)=0$ and
$$
x(t)=\int_0^ts^{\alpha-1}\sin\Big(\log\frac{1}{s}\Big){\rm d}s,\ \ \forall t>0
$$
which satisfies $x\in AC[0,T]\cap C^1(0,T]$ and
$x'(t)=t^{\alpha-1}\sin(\log\frac{1}{t})=\mathcal{O}(t^{\alpha-1})\in L^1[0,T]$ as $t\to0$. %
Scaling $s=tu$ and writing $\theta=\log(1/t)$, we obtain
$$
\begin{array}{l}\disp
{}^C_0\!D_t^\alpha x(t)
=\frac1{\Gamma(1-\alpha)}\int_0^1(1-u)^{-\alpha}u^{\alpha-1}
\sin\big(\theta+\log(1/u)\big)\,\mathrm du\crr\disp\hspace{1.4cm}
=\frac1{\Gamma(1-\alpha)}\big(A\sin\theta+B\cos\theta\big),
\end{array}
$$
where $A:=\int_0^1(1-u)^{-\alpha}u^{\alpha-1}\cos\log(1/u)\,\mathrm du$
and $B:=\int_0^1(1-u)^{-\alpha}u^{\alpha-1}\sin\log(1/u)\,\mathrm du$.
It is trivial to see $A+iB=B(1-\alpha,\alpha-i)
=\Gamma(1-\alpha)\Gamma(\alpha-i)/\Gamma(1-i)\neq0$,
hence the amplitude $\sqrt{A^2+B^2}/\Gamma(1-\alpha)$ is strictly positive. As
$t\to0^+$, $\theta\to\infty$, and the expression oscillates with non-zero amplitude, so
it has no limit.
Hence \dref{eq:Vx-Caputo-ine-min} at $t=0$ is understood in the
limiting sense and holds whenever the two limits exist. This is automatic
under the stronger hypothesis $x(t)=F(t^{\alpha})$ with $F\in C^1[0,+\infty;\mathbb{R}^n)$.
It is worth noting that the assumption
(the limit $\xi_0:=\lim_{t\to0^+}{}^C_0D_t^{\alpha}x(t)\in\mathbb{R}^n$) does not impose overly restrictive conditions
on the use of the inequality \dref{eq:Vx-Caputo-ine-min}.
The reason is quite simple:
if $f(t,x(t))$ is continuous and the equation \dref{sys-plantFrac} admits a solution,
that solution satisfies the assumption automatically.
\end{remark}

\begin{remark}
Taking $V(x)=x^2$ in Lemma \ref{lem:Vx-Caputo-ine-min-HC}, the inequality
\dref{eq:Vx-Caputo-ine-min} becomes
${}^C_0\!D^\alpha_t x^2(t)\leq 2x(t){}^C_0\!D^\alpha_t x(t)$.
This inequality was first proved in \cite{AAAlikhanov2008-2010DE} under the assumption
that $x(t)$ is absolutely continuous on $[0,T]$. However, absolute continuity alone
cannot guarantee that the inequality holds at $t=0$.
Consider the same construction as in Remark \ref{remark-D-Vx-ineq-exam}, but take
$x(0)=1$
and $x(t)=1+\int_0^t u^{\alpha-1}\sin(\log\frac1u)\,\mathrm du$.
The added constant has zero Caputo derivative, so
${}^C_0\!D_t^\alpha x(t)
=\frac1{\Gamma(1-\alpha)}(A\sin\theta+B\cos\theta)
$
still oscillates with non-zero amplitude and has no limit as $t\to0^+$. Since
$x(t)\to1$, the product $2x(t)\,{}^C_0\!D_t^\alpha x(t)$ also oscillates with non-zero
amplitude and has no endpoint limit, while $x$ remains absolutely continuous and
$x'(t)=\mathcal{O}(t^{\alpha-1})$. Hence the assumption in Lemma
\ref{lem:Vx-Caputo-ine-min-HC} that the limit
$\xi_0:=\lim_{t\to0^+}{}^C_0D_t^\alpha x(t)\in\mathbb{R}^n$ cannot be replaced by
absolute continuity alone. The inequality
${}^C_0\!D^\alpha_t x^2(t)\leq 2x(t){}^C_0\!D^\alpha_t x(t)$ was reproved in
\cite{AguilaCamacho2014cns} where $x(t)$ is assumed to be differentiable. As pointed
out in Remark \ref{remark-Vx-fi}, this inequality with the differentiability requirement
cannot be applied directly to the solutions of fractional differential equations.
\end{remark}

\begin{theorem}\label{thm-gen-localsolution}
Suppose that $f(t,x)$ in \dref{sys-plantFrac} is continuous with respect to $(t,x)$.
Then \dref{sys-plantFrac} admits at least one local solution $x(t)\in C[0,\delta]$ for some $\delta>0$.
Moreover, if there exist constants  $L,L_1>0$ and $\beta>1-\alpha$ such that $f(t,x)$ satisfies
\begin{equation}\label{eq:fun-f-sat-Lip}
\|f(t_1,x_1)-f(t_2,x_2)\|\leq L\|x_1-x_2\|+L_1|t_2-t_1|^\beta,
\end{equation}
then the solution of \dref{sys-plantFrac} is unique.
Further, $x\in C^{0,\alpha}[0,\delta]\cap C^1(0,\delta]$ and $x'(t)=\mathcal{O}(t^{\alpha-1})$ as $t\to0^+$.
Let $y(t)=x(t^{1/\alpha})$. Then  $y\in C^1[0,\delta^\alpha]$.
\end{theorem}
\noindent{\bf Proof.}
{\bf Step 1.} We prove the first assertion.
For $T>0$ put $X_T:=C[0,T]$ with the supremum norm, and define the operator $\mathcal{T}:X_T\to X_T$ by
\begin{equation}\label{Thm11-op}
(\mathcal{T}x)(t):=x_0+\frac1{\Gamma(\alpha)}\int_0^t(t-s)^{\alpha-1}f(s,x(s))\mathrm{d}s.
\end{equation}
For $x\in X_T$,  the function $g(s):=f(s,x(s))$ is continuous on $[0,T]$. Let $\|g\|_\infty:=\max\{|g(t)|:t\in[0,T]\}$.
Since $\int_0^{t_1}(t_1-s)^{\alpha-1}\mathrm{d}s=t_1^\alpha/\alpha$,
$\int_0^{t_1}(t_2-s)^{\alpha-1}\mathrm{d}s=(t_2^\alpha-(t_2-t_1)^\alpha)/\alpha$ and $\int_{t_1}^{t_2}(t_2-s)^{\alpha-1}\mathrm{d}s=(t_2-t_1)^\alpha/\alpha$,
it is easy to see that for  $0\le t_1\le t_2\le T$,
\begin{equation}\label{Thm11-cont}
\Big\|\int_0^{t_2}(t_2-s)^{\alpha-1}g(s)\mathrm{d}s-\int_0^{t_1}(t_1-s)^{\alpha-1}g(s)\mathrm{d}s\Big\|
\le \frac{2\|g\|_\infty}{\alpha}(t_2-t_1)^{\alpha}, 
\end{equation}
which shows that $t\mapsto(\mathcal{T}x)(t)$ is continuous on $[0,T]$, and thus 
$\mathcal{T}$ is well defined.

Since $\Omega$ is a domain containing $x_0$, pick $R>0$ with $\overline{B}(x_0,R)\subset\Omega$.
The set $[0,T]\times\overline{B}(x_0,R)$ is compact and $f$ is continuous on it, so
$
M:=\sup\{\|f(t,\xi)\|:\ t\in[0,T],\ \xi\in\overline{B}(x_0,R)\}<\infty.
$
Choose $\delta>0$ such that $\frac{M\delta^{\alpha}}{\Gamma(\alpha+1)}\leq R$, and let
$
K_R:=\big\{x\in X_\delta:\ \|x-x_0\|_{C[0,\delta]}\le R\big\}
$
be a nonempty, closed, convex and bounded subset of 
$X_\delta$.

(a) $\mathcal{T}(K_R)\subseteq K_R$. For $x\in K_R$ and $t\in[0,\delta]$,
$$
\|(\mathcal{T}x)(t)-x_0\|\le\frac{M}{\Gamma(\alpha)}\int_0^t(t-s)^{\alpha-1}\mathrm{d}s=\frac{Mt^\alpha}{\Gamma(\alpha+1)}\leq R.
$$
(b) $\mathcal{T}$ is continuous. If $x_n\to x$ in $X_\delta$, then for $n$ large the values $x_n(s)$
and $x(s)$ lie in $\overline{B}(x_0,R)$,
and $f$ is uniformly continuous on $[0,\delta]\times\overline{B}(x_0,R)$,
so $f(\cdot,x_n(\cdot))\to f(\cdot,x(\cdot))$ uniformly on $[0,\delta]$.
Applying \dref{Thm11-cont} to $g=f(\cdot,x_n(\cdot))-f(\cdot,x(\cdot))$ yields $\mathcal{T}x_n\to\mathcal{T}x$ in $X_\delta$.\\
(c) $\mathcal{T}(K_R)$ is relatively compact. By (a) the set $\mathcal{T}(K_R)$ is uniformly bounded,
and for $x\in K_R$, $0\le t_1\le t_2\le\delta$, \dref{Thm11-cont} with $g=f(\cdot,x(\cdot))$, $\|g\|_\infty\le M$, gives
$$
\|(\mathcal{T}x)(t_2)-(\mathcal{T}x)(t_1)\|\leq\frac{2M}{\alpha\Gamma(\alpha)}(t_2-t_1)^{\alpha}\to0
\quad\text{as }|t_2-t_1|\to0,
$$
uniformly in $x\in K_R$. Hence $\mathcal{T}(K_R)$ is equicontinuous.
By the Arzel\`a--Ascoli theorem, $\mathcal{T}(K_R)$ is relatively compact in $X_\delta$.

By Schauder's fixed point theorem, $\mathcal{T}$ admits a fixed point $x\in K_R$,
i.e. a solution of \dref{Thm11-IE} on $[0,\delta]$, hence a local solution of \dref{sys-plantFrac}.

{\bf Step 2.}
Let $x_1,x_2\in C[0,\delta]$ be two solutions of \dref{Thm11-IE}
(shrinking $\delta$ if necessary so that their values stay in $\overline{B}(x_0,R)$),
and set $z:=x_1-x_2$. It is seen that 
$$
\|z(t)\|\le\frac{L}{\Gamma(\alpha)}\int_0^t(t-s)^{\alpha-1}\|z(s)\|\,\mathrm{d}s,\quad t\in[0,\delta].
$$
Iterating this inequality and using that the fractional integration operator $\Ial$ is
positive with ${}_0I^{n\alpha}_t1=\frac{t^{n\alpha}}{\Gamma(n\alpha+1)}$, we obtain, for every $n\in\mathbb{N}$,
$$
\|z(t)\|\leq{L^n}\big({}_0I^{n\alpha}_t\|z\|\big)(t)
\leq\frac{(L\delta^{\alpha})^n}{\Gamma(n\alpha+1)}\|z\|_{C[0,\delta]},\quad t\in[0,\delta].
$$
By Stirling's formula, 
$1/\Gamma(n\alpha+1)$ decays faster than any exponential,
so for sufficiently large $n$,
$\frac{(L\delta^{\alpha})^n}{\Gamma(n\alpha+1)}<1$.
Hence $z(t)\equiv0$ for $t\in[0,\delta]$. The local solution is unique.

{\bf Step 3.} 
Let $K$ be a compact neighbourhood of the trajectory $x([0,\delta])$ in $\Omega$ and
$M:=\max\{\|f(t,\xi)\|:t\in[0,T],\ \xi\in K\}$. For $0\le t_1<t_2\le T$, using \dref{Thm11-IE},
\[
x(t_2)-x(t_1)=\frac1{\Gamma(\alpha)}\Big[\int_0^{t_1}\big((t_2-s)^{\alpha-1}-(t_1-s)^{\alpha-1}\big)f(s,x(s))\mathrm{d}s+\int_{t_1}^{t_2}(t_2-s)^{\alpha-1}f(s,x(s))\mathrm{d}s\Big],
\]
and since $(t-s)^{\alpha-1}$ is decreasing in $t$,
\begin{equation}\label{Thm11-Holder}
\|x(t_2)-x(t_1)\|\le\frac{M}{\Gamma(\alpha+1)}\big(t_1^\alpha-t_2^\alpha
  +2(t_2-t_1)^\alpha\big)\le\frac{2M}{\Gamma(\alpha+1)}(t_2-t_1)^\alpha.
\end{equation}
Thus $x\in C^{0,\alpha}[0,\delta]$.

{\bf Step 4.} Let $\phi(t):=f(t,x(t))$, $t\in[0,\delta]$.
By \dref{eq:fun-f-sat-Lip} and \dref{Thm11-Holder}, $\phi\in C^{0,\eta}[0,\delta]$ with $\eta=\min\{\alpha,\beta\}$.
We use the local smoothing property of the Abel operator  (\cite[Theorems 3.1]{SamkoKilbasMarichev1993}):
if $\psi\in C^{0,\eta}$ on a neighbourhood of $t_0>0$, then ${}_0I_t^\alpha\psi\in C^{0,\eta+\alpha}$
on a neighbourhood of $t_0$, provided $\eta+\alpha<1$.
Iterating in $x(t)=x_0+{}_0I_t^\alpha\phi$ finitely many times (note that $\alpha>0$),
we obtain $x\in C^{0,\beta}_{\rm loc}(0,\delta]$ with $\beta\in(1-\alpha,1)$.
Fix $t\in(0,T]$ and let $\tau>0$ be small. From \dref{Thm11-IE},
$$
\frac{x(t+\tau)-x(t)}{\tau}
=\frac{1}{\Gamma(\alpha)\tau}\Big[\int_0^t\big((t+\tau-s)^{\alpha-1}-(t-s)^{\alpha-1}\big)\phi(s)\mathrm{d}s
  +\int_t^{t+\tau}(t+\tau-s)^{\alpha-1}\phi(s)\mathrm{d}s\Big].
$$
From this and 
arguing as in Step 6 of the proof of Theorem \ref{Thm-Dhinfatzero-mono-asympt}, we obtain $x\in C^1(0,\delta]$ and
\begin{equation}\label{Thm11-der}
x'(t)=\frac{1}{\Gamma(\alpha)}\Big[\phi(t)t^{\alpha-1}+(\alpha-1)\int_0^t(t-s)^{\alpha-2}\big(\phi(s)-\phi(t)\big)\mathrm{d}s\Big].
\end{equation}
%
%
{
Next, we show that $x'(t)=\mathcal O(t^{\alpha-1})$ as $t\to0^+$. To this end,
consider the Picard iterates
\[
x_0(t)\equiv x_0,\qquad
x_{n+1}(t)=x_0+\frac1{\Gamma(\alpha)}\int_0^t(t-s)^{\alpha-1}f(s,x_n(s))\,\mathrm ds,
\]
and set $q_n(t):=\|x_n'(t)\|$ for $t>0$ and $\phi_n(t):=f(t,x_n(t))$ for $t\geq0$.
On a small interval $[0,\delta]$ all iterates remain in a compact set,
so $\|f(s,x_n(s))\|\leq M$ for some constant $M>0$.
Applying the Volterra derivative identity to $x_{n+1}=x_0+{}_0^CI^\alpha_t\phi_n$
(which applies on $(0,\delta]$ because $\phi_n$ is locally H\"older of exponent
$\beta>1-\alpha$), and using
$$
\|f(s,x_n(s))-f(t,x_n(t))\|\le L\|x_n(s)-x_n(t)\|+L_1|s-t|^\beta,
$$
we obtain
\begin{equation}\label{eq:qn-seq-estimate-I}
q_{n+1}(t)\le \frac{M}{\Gamma(\alpha)}t^{\alpha-1}
+\frac{1-\alpha}{\Gamma(\alpha)}\int_0^t(t-s)^{\alpha-2}
\big[L\|x_n(s)-x_n(t)\|+L_1(t-s)^\beta\big]\,\mathrm ds.
\end{equation}
For $0<s<t$, $\|x_n(s)-x_n(t)\|\le\int_s^t q_n(r)\mathrm dr$. Substituting
this and swapping the order of integration, we have
\[
(1-\alpha)\int_0^t(t-s)^{\alpha-2}\int_s^t q_n(r)\,\mathrm dr\,\mathrm ds
=\int_0^t\big[(t-r)^{\alpha-1}-t^{\alpha-1}\big]q_n(r)\,\mathrm dr
\le \Gamma(\alpha)({}_0I_t^\alpha q_n)(t).
\]
The $\beta$-H\"older part in \dref{eq:qn-seq-estimate-I} is finite because $\alpha+\beta>1$ and
$
\frac{(1-\alpha)L_1}{\Gamma(\alpha)}
\int_0^t(t-s)^{\alpha+\beta-2}\,\mathrm ds
=\frac{(1-\alpha)L_1}{\Gamma(\alpha)(\alpha+\beta-1)}t^{\alpha+\beta-1}
\le \frac{(1-\alpha)L_1\delta^\beta}{\Gamma(\alpha)(\alpha+\beta-1)}t^{\alpha-1}.
$
Set
$C_1:=\frac{M}{\Gamma(\alpha)}
+\frac{(1-\alpha)L_1\delta^\beta}{\Gamma(\alpha)(\alpha+\beta-1)}$. Then,
\[
q_{n+1}(t)\leq C_1 t^{\alpha-1}+L\,({}_0I_t^\alpha q_n)(t),\qquad t\in(0,\delta].
\]
Iterating this recursion and using
${}_0I_t^{j\alpha}(t^{\alpha-1})=\frac{\Gamma(\alpha)}{\Gamma((j+1)\alpha)}t^{(j+1)\alpha-1}$,
together with $t^{j\alpha}\le\delta^{j\alpha}$, gives
\[
q_n(t)\leq C_1\Gamma(\alpha)
\sum_{j=0}^{\infty}\frac{(L\delta^\alpha)^j}{\Gamma((j+1)\alpha)}t^{\alpha-1}.
\]
It is trivial to see that the series
$\sum_{j=0}^{\infty}\frac{(L\delta^\alpha)^j}{\Gamma((j+1)\alpha)}$
converges by the ratio test. Define
$C_2:=C_1\Gamma(\alpha)
\sum_{j=0}^{\infty}\frac{(L\delta^\alpha)^j}{\Gamma((j+1)\alpha)}<+\infty$,
and we have, uniformly in $n$,
$
q_n(t)\leq C_2t^{\alpha-1}$, $t\in(0,\delta],
$
and thus
$
\|x_n(t)-x_n(s)\|\leq\int_s^t q_n(r)\,\mathrm dr
\leq C_2\int_s^t r^{\alpha-1}\mathrm dr
=C_2(t^\alpha-s^\alpha)/{\alpha}$, for all  $s<t$ and $s,t\in(0,\delta].
$
Passing to the limit in this inequality 
gives 
\begin{equation}\label{eq:sys-solution-x-t-s-est}
\|x(t)-x(s)\|\leq\frac{C_2}{\alpha}(t^\alpha-s^\alpha),\qquad 0< s<t\le\delta.
\end{equation}
For any fixed $t\in(0,\delta]$, dividing both sides of \dref{eq:sys-solution-x-t-s-est} by
$t-s$, letting $s\to t$, and noticing $(t^\alpha-s^\alpha)/(t-s)\to\alpha t^{\alpha-1}$ as $s\to0$,
we obtain
$\|x'(t)\|\leq\frac{C_2}{\alpha}t^{\alpha-1}$ on $(0,\delta]$.

We define $K(t):=\int_0^t(t-s)^{\alpha-2}\big(\phi(s)-\phi(t)\big)\,\mathrm ds$
 and split $K(t)=K_1(t)+K_2(t)$ at $s=t/2$.
Since $\|\phi(s_1)-\phi(s_2)\|\leq C_2|s_1-s_2|^\eta$, for all $s_1,s_2\in[0,\delta]$ and some $C_2>0$,
we have
$\|\phi(s)-\phi(u)\|\le \|\phi(s)-\phi(0)\|+\|\phi(0)-\phi(u)\|
\leq C_2\big(s^\eta+u^\eta\big)\leq 2C_3 u^\eta$ for $s\in[0,u/2]$.
Together with $(t-s)^{\alpha-2}\le(t/2)^{\alpha-2}$, this yields
\begin{equation}\label{eq:thm-K1-estim}
\|K_1(t)\|
\leq 2(t/2)^{\alpha-2}C_3t^\eta\cdot\frac t2
=C_3\,2^{2-\alpha}t^{\alpha-1+\eta}=o(t^{\alpha-1}).
\end{equation}
Since $t^\alpha-s^\alpha\leq2\alpha t^{\alpha-1}(t-s)$ for $s\in[t/2,t]$,
by \dref{eq:fun-f-sat-Lip} and \dref{eq:sys-solution-x-t-s-est}, we obtain
\begin{equation}\label{eq:thm-K2-estim}
\begin{aligned}
\|K_2(t)\|
&\le\int_{t/2}^t(t-s)^{\alpha-2}
\Big[L\frac{C_2}{\alpha}(t^\alpha-s^\alpha)+L_1(t-s)^\beta\Big]\mathrm ds\\
&\le \frac{2^{1-\alpha}LC_2}{\alpha}t^{2\alpha-1}
+\frac{L_1}{2^{\alpha+\beta-1}(\alpha+\beta-1)}\,t^{\alpha+\beta-1}
=o(t^{\alpha-1}),
\end{aligned}
\end{equation}
Combining this with $\|K_1(t)\|
\leq C_2\,2^{2-\alpha}t^{\alpha-1+\eta}$ and
\dref{Thm11-der}, we conclude
\[
x'(t)=\frac{f(0,x_0)}{\Gamma(\alpha)}t^{\alpha-1}
+\mathcal O(t^{\alpha-1+\eta}+t^{2\alpha-1}+t^{\alpha+\beta-1}),
\]
in particular $x'(t)=\mathcal O(t^{\alpha-1})$ as $t\to0^+$.

}

{\bf Step 5.} By the chain rule, for $t>0$, we have
\begin{equation}\label{Thm11-yder}
y'(t)=x'(t^{1/\alpha})\,\frac{1}{\alpha}\,t^{\frac{1}{\alpha}-1}
=\frac{1}{\alpha\Gamma(\alpha)}\Big[\phi(t^{1/\alpha})
+(\alpha-1)\,t^{\frac{1}{\alpha}-1}
\int_0^{t^{1/\alpha}}(t^{1/\alpha}-s)^{\alpha-2}\big(\phi(s)-\phi(t^{1/\alpha})\big)\,\mathrm{d}s\Big],
\end{equation}
where we used \dref{Thm11-der}. Put $u=t^{1/\alpha}$,
so that $t=u^\alpha$ and $t^{\frac{1}{\alpha}-1}=u^{1-\alpha}$.
Split the integral in \dref{Thm11-yder} at $s=u/2$:
$$
\int_0^u(u-s)^{\alpha-2}(\phi(s)-\phi(u))\mathrm ds
=J_1(u)+J_2(u),
$$
where
$
J_1(u):=\int_0^{u/2}(u-s)^{\alpha-2}(\phi(s)-\phi(u))\mathrm ds$,
$J_2(u):=\int_{u/2}^u(u-s)^{\alpha-2}(\phi(s)-\phi(u))\mathrm ds$.
From \dref{eq:thm-K1-estim} and \dref{eq:thm-K2-estim}, it is seen that
$\|J_1(u)\|=o(u^{\alpha-1})$ as $u\to0^+$ and $\|J_2(u)\|=o(u^{\alpha-1})$ as $u\to0^+$.
Hence,
$$
u^{1-\alpha}\int_0^u(u-s)^{\alpha-2}(\phi(s)-\phi(u))\mathrm ds\to0,\ \ \mbox{as }u\to0^+,
$$
which, together with \dref{Thm11-yder} and the continuity of $\phi$, gives
$$
\Big\|y'(t)-\frac{\phi(0)}{\alpha\Gamma(\alpha)}\Big\|
\le\frac{\|\phi(t^{1/\alpha})-\phi(0)\|}{\alpha\Gamma(\alpha)}
+\frac{(1-\alpha)u^{1-\alpha}}{\alpha\Gamma(\alpha)}
   \Big\|\int_0^u(u-s)^{\alpha-2}(\phi(s)-\phi(u))\mathrm ds\Big\|\to0,
$$
as $t\to0^+$.
Hence $\lim_{t\to0^+}y'(t)=\phi(0)/\Gamma(\alpha+1)=f(0,x_0)/\Gamma(\alpha+1)$.
It remains to identify this limit with the right derivative of $y$ at $0$.
Proceeding exactly as in Step 6 of the proof of Theorem \ref{Thm-Dhinfatzero-mono-asympt} for the value $W'_+(0)$,
and using $u=t^{1/\alpha}$, we obtain
\[
\frac{y(t)-y(0)}t
=\frac{x(u)-x_0}{u^\alpha}
=\frac{1}{u^\alpha\Gamma(\alpha)}\int_0^u(u-s)^{\alpha-1}\phi(s)\,\mathrm{d}s
\longrightarrow\frac{\phi(0)}{\alpha\Gamma(\alpha)}
=\frac{f(0,x_0)}{\Gamma(\alpha+1)},\;\mbox{as }t\to0,
\]
where $\frac1{u^\alpha}\int_0^u(u-s)^{\alpha-1}\phi(s)\mathrm{d}s\to\phi(0)\int_0^1(1-\sigma)^{\alpha-1}\mathrm{d}\sigma=\phi(0)/\alpha$ is used.
Thus $y'_+(0)=f(0,x_0)/\Gamma(\alpha+1)=\lim_{t\to0^+}y'(t)$. Hence, $y\in C^1[0,\delta^\alpha]$.
\hfill$\Box$

\begin{remark}
For any fixed $a>0$, let  $(t-a)_{+}=\max\{t-a,0\}$.
It is trivial to verify that
$x(t)=(\Gamma(1+\frac{\alpha}{2})/\Gamma(1+\frac{3\alpha}{2}))^{3/2}
(t-a)_{+}^{\frac{3\alpha}{2}}$
solves fractional differential equation
 ${}^C_0D_t^{\alpha}x(t)=x^{1/3}(t)$ with $x(0)=0$.
This shows that the solution of \dref{sys-plantFrac} need not be unique
if $f(t,x)$ does not satisfy \eqref{eq:fun-f-sat-Lip}.
\end{remark}

\begin{theorem}\label{tm-Vconvex-Lya-zero}
Let $x = 0$ be an equilibrium point for the system \dref{sys-plantFrac}
and $\Omega$ be a domain containing the origin.
Let $f(t,x)$ satisfy \eqref{eq:fun-f-sat-Lip}.
Let $V:\mathbb{R}^n\to\mathbb{R}^+$ be a function satisfying the following conditions:
$V(\cdot)$ is convex and differentiable on $\mathbb{R}^n$.
Assume that there exist three class-$\mathcal{K}$ functions  $\alpha_i$, $i = 1, 2, 3$
satisfying
\begin{equation}\label{nabla-V-ineq-al-tm-genal}
\begin{array}{l}
\alpha_1(\|x\|)\leq V(x(t))
\leq\alpha_2(\|x\|),\\
\langle \nabla V(x),f(t,x)\rangle\leq-\alpha_3(\|x\|).
\end{array}
\end{equation}
Then the equilibrium point of system \dref{sys-plantFrac} is asymptotically stable.
\end{theorem}
\noindent{\bf Proof.}
Since $V(\cdot)$ is convex and differentiable on $\mathbb{R}^n$,
by Theorem \ref{thm-gen-localsolution} and Lemma \ref{lem:Vx-Caputo-ine-min-HC}, we know that
 ${}^C_0D_t^{\alpha}V(x(t))\leq \langle \nabla V(x),f(t,x)\rangle$.
 Combining this with \dref{nabla-V-ineq-al-tm-genal} and $\|x(t)\|\geq\alpha_2^{-1}(V(x(t)))$,
 we have ${}^C_0D_t^{\alpha}V(x(t))\leq-\alpha_3(\alpha_2^{-1}(V(x(t))))$.
 Then the desired result follows from Theorem \ref{thm-DV-hV-gencase}. 
\hfill$\Box$

\begin{remark}
Compared  with the classic Lyapunov theory \cite{YoshizawaBook},
Theorem \ref{tm-Vconvex-Lya-zero} requires the Lyapunov function $V(x)$ to be convex,
which is not assumed in the classic Lyapunov theory
for ordinary differential equations or partial differential equations.
The reason is the lack of chain rule for fractional derivative.
\end{remark}

\section{Existence, uniqueness, regularity and decay properties of solution under class-$\mathcal K$ only }\label{Sec-h-classK-property}

In this section, for fixed $A>0$ and $0<\alpha<1$,
we continue to consider the weakly singular Volterra integral equation:
\begin{equation}\label{eq:Vh-rewriten-standard}
V(t)=A-\frac{1}{\Gamma(\alpha)}\int_0^t(t-s)^{\alpha-1}h(|V(s)|)\,\mathrm ds,
\qquad t\ge 0,
\end{equation}
where $h:[0,\infty)\to[0,\infty)$ is a class-$\mathcal K$ function.
In Section \ref{Sec-h-diff-property}, under the assumption that
the derivative function $h'$ satisfies local Lipschitz continuity on $(0,+\infty)$,
Theorem \ref{Thm-Dhinfatzero-mono-asympt} and Remark \ref{rem-exist-uniq-Thm41}
prove that \eqref{eq:Vh-rewriten-standard} admits a unique global solution $V\in C[0,+\infty)$ and that it satisfies \eqref{eq:V-properties}.
Now, no differentiability or Lipschitz assumption on $h'$ 
 only the class-$\mathcal K$ properties of $h$ is required, 
 by utilizing the results of Theorem \ref{Thm-Dhinfatzero-mono-asympt},
 we prove that  \eqref{eq:Vh-rewriten-standard} has a unique global solution $V(t)$
 which is positive, bounded, non-increasing, H\"older continuous.

We state the global positive solution under the class-$\mathcal K$ assumption only.
\begin{theorem}
\label{thm:global-classK}
Let $h$ be a class-$\mathcal K$ function, $A>0$, and $0<\alpha<1$. Then
equation \eqref{eq:Vh-rewriten-standard} possesses a solution $V\in C^{0,\alpha}_{\rm loc}[0,\infty)$ satisfying
\begin{equation}\label{eq:V-properties}
(i) 0<V(t)\le A,\ \forall t\ge 0,\quad (ii)
V(t)\ \text{is decreasing on }[0,+\infty), \quad
(iii)\lim_{t\to\infty}V(t)=0.
\end{equation}
\end{theorem}
\noindent{\bf Proof.}
The proof proceeds in six steps: we smooth $h$, solve the regularized
equations, obtain estimates uniform in the smoothing parameter, extract a
limit by compactness, and then verify positivity and decay.\\
\noindent{\bf Step 1:} We construct a smooth class-$\mathcal K$ approximation.
Extend $h$ to the whole real line by zero: $\widetilde h(x):=h(x)$ for $x\geq0$ and $\widetilde h(x):=0$ for $x<0$.
Then $\widetilde h$ is continuous, non-decreasing, and non-negative.
Let $\rho\in C_c^\infty(\mathbb R)$ be a standard mollifier:
$\rho\ge0$, $\operatorname{supp}\rho\subset[-1,1]$, $\rho>0$ on $(-1,1)$,
and $\int_{\mathbb R}\rho(x)\,\mathrm dx=1$. Set $\rho_n(x):=n\rho(nx)$ and
define the convolution $(\widetilde h*\rho_n)(x):=\int_{\mathbb R}\widetilde h(x-y)\rho_n(y)\,\mathrm dy$.
Define
$
h_n(x):=(\widetilde h*\rho_n)(x)-(\widetilde h*\rho_n)(0), \ x\geq 0.
$ 
We record the properties of $h_n$.

(a) $h_n$ is the restriction of a $C^\infty$ function and $h_n(0)=0$ by
construction.

(b) $h_n\ge0$. Indeed,
\begin{equation}\label{eq:h_n-int}
h_n(x)=\int_{\mathbb R}\rho_n(y)\big[\widetilde h(x-y)-\widetilde h(-y)\big]\,\mathrm dy,
\end{equation}
and for $x\ge0$ the integrand is non-negative because $\widetilde h$ is
non-decreasing.

(c) $h_n$ is strictly increasing. For $x_2>x_1\ge0$,
\begin{equation}\label{eq:h_n-strict}
h_n(x_2)-h_n(x_1)
=\int_{\mathbb R}\rho_n(y)\big[\widetilde h(x_2-y)-\widetilde h(x_1-y)\big]\,\mathrm dy\ge0.
\end{equation}
Note that the only $y$ that contributes in the integrand \eqref{eq:h_n-strict}
satisfies $y\in\operatorname{supp}\rho_n\subset[-1/n,1/n]$.
Let $S_n:=\bigl(-\tfrac1n,\tfrac1n\bigr)\cap(-\infty,x_1]=(-\tfrac1n,\min\{x_1,\tfrac1n\})$.
This interval has a positive measure, and $\rho_n(y)>0$ on it.
For every $y\in S_n$, one has $y\le x_1$, hence $x_1-y\ge0$ and $x_2-y>x_1-y\ge0$.
Since $\widetilde h=h$ is strictly increasing on $[0,\infty)$, we obtain 
$
\widetilde h(x_2-y)-\widetilde h(x_1-y)
=h(x_2-y)-h(x_1-y)>0
$
(when $x_1-y=0$, this reads $h(x_2-y)>h(0)=0$). Thus the integrand in
\eqref{eq:h_n-strict} is strictly positive on $S_n$ and non-negative
everywhere, so $h_n(x_2)-h_n(x_1)>0$.

Thus each $h_n$ is itself a smooth class-$\mathcal K$ function.

(d) $h_n$ converges to $h$ uniformly on $[0,A]$, and $h_n$ is uniformly
dominated there. Since $\widetilde h(x-y)\leq \widetilde h(x+\tfrac1n)=h(x+\tfrac1n)$ for $y\ge-\tfrac1n$,
\eqref{eq:h_n-int} gives
$
0\le h_n(x)\le h\!\left(x+\frac1n\right),\ x\ge0,
$ 
so that in particular, for $n\ge1$,
$
h_n(A)\le h(A+1)=:M.
$ 
Moreover, the standard uniform convergence of mollifications on compact sets
and the estimate $(\widetilde h*\rho_n)(0)=\int_{y\le0}\rho_n(y)h(-y)\,\mathrm dy
\le h(1/n)\to0$ give
\begin{equation}\label{eq:h_n-uniform}
\sup_{0\le x\le A}|h_n(x)-h(x)|\longrightarrow0\qquad(n\to\infty).
\end{equation}

\noindent{\bf Step 2:} We claim that the regularized equations have global positive solutions.
For each $n$, $h_n$ is smooth, hence $h_n'$ is locally Lipschitz on
$(0,\infty)$ (in fact on $[0,\infty)$, since $h_n'(0)$ exists and is finite).
Thus $h_n$ satisfies the hypotheses of Theorem \ref{Thm-Dhinfatzero-mono-asympt}.
Applying Theorem \ref{Thm-Dhinfatzero-mono-asympt}(i)--(iii)
together with Remark \ref{rem-exist-uniq-Thm41}, we have that
\begin{equation}\label{eq:Vn}
V_n(t)=A-(\Ial h_n(|V_n|))(t),\qquad t\ge0,
\end{equation}
admits a unique solution $V_n\in C[0,\infty)$ with
$
0<V_n(t)\le A,\ t\ge0,
$ 
and $V_n$ is decreasing on $[0,\infty)$ and converges to $0$. Moreover,
Theorem \ref{Thm-Dhinfatzero-mono-asympt}(ii) asserts that $W_n(t):=V_n(t^{1/\alpha})$ is differentiable
with $W_n'(t)\le0$ for $t\ge0$ (equivalently, $V_n$ is non-increasing).

\noindent{\bf Step 3:} We show the uniform H\"older estimate for $V_n(t)$.
Set $g_n(s):=h_n(|V_n(s)|)$. 
By $0<V_n(t)\le A$ for all $t\geq0$ and $h_n(A)\le h(A+1)=M$,
 we have
$0\le g_n(s)\le h_n(A)\le M=h(A+1),\ s\ge0.$
For $0\le t_1\le t_2\le T$ and every $n$, it is easy to see
$$
|(\Ial g_n)(t_2)-(\Ial g_n)(t_1)|
\le \frac{M}{\Gamma(\alpha)}
\left[\int_{t_1}^{t_2}(t_2-s)^{\alpha-1}\,\mathrm ds
+\int_0^{t_1}\big|(t_2-s)^{\alpha-1}-(t_1-s)^{\alpha-1}\big|\,\mathrm ds\right].
$$
Clearly, the first integral equals $(t_2-t_1)^\alpha/\alpha$ and the second integral
$\int_0^{t_1}\big|(t_2-s)^{\alpha-1}-(t_1-s)^{\alpha-1}\big|\,\mathrm ds
={(t_1^\alpha-t_2^\alpha+(t_2-t_1)^\alpha)}/{\alpha}
\leq (t_2-t_1)^\alpha/\alpha$. It follows from $V_n=A-\Ial g_n$ that for $0\le t_1\le t_2\le T$ and every $n$,
\begin{equation}\label{eq:Holder-HC-G}
|V_n(t_2)-V_n(t_1)|
\le \frac{2M}{\alpha\,\Gamma(\alpha)}\,(t_2-t_1)^\alpha.
\end{equation}

\noindent{\bf Step 4:} We show that \eqref{eq:Vh-rewriten-standard} has a solution
$V\in C^{0,\alpha}_{\rm loc}[0,\infty)$ and it satisfies $0\leq V(t)\leq A$ and (ii).
For every fixed $T>0$, the family $\{V_n\}$ is uniformly bounded in
$C[0,T]$ by $0<V_n(t)\le A,\ t\ge0$ 
and equicontinuous by
\eqref{eq:Holder-HC-G}. By the Arzel\`a--Ascoli theorem, a subsequence converges
uniformly on $[0,T]$. Applying this successively on $[0,k]$ for
$k=1,2,\dots$ and taking a diagonal subsequence (still denoted $\{V_n\}$), we
obtain
$
V_n\longrightarrow V$ uniformly on every compact interval of $[0,\infty)$.
For fixed $0\le t_1<t_2\le T$, by passing to the limit in \eqref{eq:Holder-HC-G}, we have
\begin{equation}\label{eq:V-Holder}
|V(t_2)-V(t_1)|
=\lim_{n\to\infty}|V_n(t_2)-V_n(t_1)|
\le \frac{2M}{\alpha\,\Gamma(\alpha)}\,(t_2-t_1)^\alpha.
\end{equation}
Thus $V\in C^{0,\alpha}[0,T]$ for every $T$, hence $V\in C^{0,\alpha}_{\rm loc}[0,\infty)$.
The pointwise bounds $0<V_n(t)\le A,\ t\ge0$ and the non-increasing property $V_n(t_1)\geq V_n(t_2)$ for $t_1\leq t_2$
are likewise preserved by the uniform limit, so $0\leq V(t)\leq A$ and $V$
is non-increasing.

Let $\omega_h(\delta)=\sup\{|h(x)-h(y)|:x,y\in[0,A],\ |x-y|\le\delta\}$ be the
modulus of continuity of $h$ on $[0,A]$. Since $h$ is continuous on $[0,+\infty)$, we have $\lim_{\delta\to0+}\omega_h(\delta)=0$.
Since $V_n,V\in[0,A]$ and $\|V_n-V\|_{C[0,T]}\to0$, it follows from \eqref{eq:h_n-uniform} that 
$$
\sup_{0\le s\le T}|h_n(|V_n(s)|)-h(V(s))|
\le \sup_{0\le x\le A}|h_n(x)-h(x)|
+\omega_h\!\left(\|V_n-V\|_{C[0,T]}\right)\longrightarrow0
$$ 
where we used $V_n>0$ and $h(|V_n|)=h(V_n)$. The estimate
$
\big|\Ial g(t)\big|\le \frac{t^\alpha}{\Gamma(\alpha+1)}\|g\|_{C[0,T]},
\ 0\le t\le T,
$
then shows that the fractional integral terms in \eqref{eq:Vn} converge
uniformly on $[0,T]$. Passing to the limit yields
\begin{equation}\label{eq:V-limit-eq}
V(t)=A-\frac{1}{\Gamma(\alpha)}\int_0^t(t-s)^{\alpha-1}h(V(s))\,\mathrm ds
=A-\frac{1}{\Gamma(\alpha)}\int_0^t(t-s)^{\alpha-1}h(|V(s)|)\,\mathrm ds,
\end{equation}
where the second equality again uses $V\ge0$.
This yields that $V(t)$ is the solution of \eqref{eq:Vh-rewriten-standard}.

{\bf Step 5:} We show $\lim_{t\to\infty}V(t)=0$.
Since $V$ is non-increasing and non-negative, the limit
$\ell:=\lim_{t\to\infty}V(t)$ exists and $\ell\ge0$. If $\ell>0$, then
$V(s)\ge\ell$ and $h(V(s))\ge h(\ell)>0$ for all $s\ge0$, so
\eqref{eq:V-limit-eq} gives
\begin{equation}\label{eq:tail}
V(t)\le A-\frac{h(\ell)}{\Gamma(\alpha)}
\int_0^t(t-s)^{\alpha-1}\,\mathrm ds
=A-\frac{h(\ell)}{\Gamma(\alpha+1)}\,t^\alpha,
\end{equation}
which is negative for all sufficiently large $t$, contradicting $V(t)\ge0$.
Hence $\ell=0$. 

{\bf Step 6:} We show the strict positivity of $V(t)$ for all $t\geq0$.
Suppose, for contradiction, that $V(t_0')=0$ for some $t_0'>0$. Let $t_0=\inf\{t_0'\geq0: V(t_0')=0\}$.
Since $V$ is non-increasing and non-negative,
we obtain $V(t)=0$ for all $t\ge t_0$ and thus $h(V(t))=0$ for $t\ge t_0$, which gives
\begin{equation}\label{eq:t0-balance-N}
A=\frac{1}{\Gamma(\alpha)}\int_0^{t_0}(t-s)^{\alpha-1}h(V(s))\,\mathrm ds.
\end{equation}
Letting $t\to+\infty$ in \dref{eq:t0-balance-N} gives $A=0$, contradicting
$A>0$. Thus, we have $V(t)>0$ for all $t\ge0$.
\hfill$\Box$

Here, there is a 
uniqueness lemma for positive solutions of \eqref{eq:Vh-rewriten-standard}.
\begin{lemma}
\label{lem:comparison}
Let $u,v\in C[0,T]$ be two positive solutions of \eqref{eq:Vh-rewriten-standard} on $[0,T]$,
so that $u(t)=A-(\Ial h(u))(t)$ and $v(t)=A-(\Ial h(v))(t)$. Then
$u\equiv v$ on $[0,T]$.
\end{lemma}
\noindent{\bf Proof.}
Define $z(t):=u(t)-v(t)$, $\varphi(t):=h(v(t))-h(u(t))$.
Both $z$ and $\varphi$ are continuous, $z(0)=0$, and subtracting the two
equations gives $z(t)=(\Ial\varphi)(t)$.
Since $h$ is strictly increasing, $\varphi$ has the opposite sign to $z$, which implies that
$z(t)\varphi(t)\leq 0$ for all $t\in[0,T]$ and $z(t)\varphi(t)<0$ whenever $z(t)\neq 0$.
%
Now we claim that the quadratic form
$$
Q:=\int_0^T z(t)\varphi(t)\,\mathrm dt
=\int_0^T (\Ial\varphi)(t)\varphi(t)\mathrm dt=0.
$$
For this purpose, we only need to show $Q\geq0$ because $Q=\int_0^T z(t)\varphi(t)\mathrm dt\leq 0$ due to $z(t)\varphi(t)\leq 0$ for all $t\in[0,T]$.
By $z(t)=(\Ial\varphi)(t)=\frac{1}{\Gamma(\alpha)}\int_0^t(t-s)^{\alpha-1}\varphi(s)\mathrm ds$,
we first rewrite $Q$ as a double integral:
$$
Q
=\int_0^T\!\left[\frac{1}{\Gamma(\alpha)}
 \int_0^t(t-s)^{\alpha-1}\varphi(s)\,\mathrm ds\right]\varphi(t)\,\mathrm dt
=\frac{1}{\Gamma(\alpha)}\iint_{0\le s\le t\le T}
 (t-s)^{\alpha-1}\varphi(t)\varphi(s)\,\mathrm ds\,\mathrm dt.
$$ 
Since $0<\alpha<1$, the reflection identity
$\Gamma(\alpha)\Gamma(1-\alpha)=\pi/\sin(\pi\alpha)$ and
$\int_0^\infty e^{-r\tau}r^{-\alpha}\,\mathrm dr=\tau^{\alpha-1}\Gamma(1-\alpha)$
give, for every $\tau>0$, $\frac{\tau^{\alpha-1}}{\Gamma(\alpha)}
=\frac{\sin(\pi\alpha)}{\pi}\int_0^\infty e^{-r\tau}r^{-\alpha}\,\mathrm dr$.
Based on this expression, we obtain
\begin{equation}\label{eq:Q-mixed}
\begin{aligned}
Q
&=\frac{\sin(\pi\alpha)}{\pi}
 \int_0^\infty r^{-\alpha}
 \left[\iint_{0\le s\le t\le T}
 e^{-r(t-s)}\varphi(t)\varphi(s)\,\mathrm ds\,\mathrm dt\right]\mathrm dr\\
 &=\frac{\sin(\pi\alpha)}{\pi}
 \int_0^\infty r^{-\alpha}
 \left(\int_0^T\varphi(t)Y_r(t)\,\mathrm dt\right)\mathrm dr.
\end{aligned}
\end{equation}
where $Y_r(t)$ is defined by $Y_r(t):=\int_0^t e^{-r(t-s)}\varphi(s)\mathrm ds$ for $0\leq t\le T$.
Clearly, $Y_r$ is the unique solution of $Y_r'(t)=\varphi(t)-rY_r(t)$ with $Y_r(0)=0$.
A simple computation shows
$$
\int_0^T \varphi(t)Y_r(t)\,\mathrm dt
=\int_0^T\big(Y_r'(t)+rY_r(t)\big)Y_r(t)\,\mathrm dt
=\frac{1}{2} Y_r(T)^2+r\int_0^T Y_r(t)^2\,\mathrm dt\ge 0.
$$
which, together with \dref{eq:Q-mixed}, gives
\begin{equation}\label{eq:Q-ge0}
\begin{aligned}
Q
&=\frac{\sin(\pi\alpha)}{\pi}\int_0^\infty r^{-\alpha}
   \left[\frac{1}{2}Y_r(T)^2+r\int_0^T Y_r(t)^2\,\mathrm dt\right]\mathrm dr
\geq 0.
\end{aligned}
\end{equation}
Hence, $Q=0$.
Finally, if $z(t_0)\neq 0$ for some $t_0\in[0,T]$, then by continuity of $z$ on $[0,T]$, there
exists a non-degenerate interval $J\subset[0,T]$ on which $z(t)\neq 0$ and
$z(t)\varphi(t)<0$ for all $t\in J$.
This gives $Q=\int_0^T z(t)\varphi(t)\,\mathrm dt<0$, contradicting $Q=0$.
Hence $z(t)=0$ for all $t\in[0,T]$, i.e.\ $u\equiv v$  on $[0,T]$.
\hfill$\Box$

\begin{theorem}
\label{thm:all-solutions}
Let $h$ be a class-$\mathcal K$ function, $A>0$, and $0<\alpha<1$.
Let $V$ be the global positive solution of \eqref{eq:Vh-rewriten-standard} given by Theorem \ref{thm:global-classK}.
If $W$ is any continuous solution of \eqref{eq:Vh-rewriten-standard} on its maximal interval of
existence $[0,T_*)$ with $W(0)=A$, then $T_*=\infty$ and
$
W(t)=V(t)$ for all $t\ge0$.
Consequently, \eqref{eq:Vh-rewriten-standard} has exactly one continuous
solution under the pure class-$\mathcal K$ hypothesis.
\end{theorem}
\noindent{\bf Proof.}
We first show that no solution can reach zero.
Since $W(0)=A>0$ and $W$ is continuous, $W>0$ on some initial interval.
Suppose that $W$ reaches a non-positive value and let
$
t_*:=\inf\{t\ge0:W(t)\le0\}.
$
Then $t_*>0$ and $W(t)>0$ for $t\in[0,t_*)$. On every compact interval
$[0,t']$ with $t'<t_*$, both $W$ and the global positive solution $V$ are
positive solutions of \eqref{eq:Vh-rewriten-standard}. Lemma~\ref{lem:comparison} gives
$W(t)=V(t)$ for  $t\in[0,t']$. Letting $t'\to t_*$ and using continuity, we obtain
$
W(t_*)=V(t_*)>0,
$
which contradicts $W(t_*)\le0$. Hence $W(t)>0$ for all $t$ in its
maximal interval.

Next, with $W(t)>0$ for  $t\in[0,T_*)$, Lemma~\ref{lem:comparison} applied on each
$[0,t']\subset[0,T_*)$ yields $W(t)=V(t)$ there.
If $T_*<\infty$, then $W(t)=V(t)$ for $t<T_*$, and $V$ is continuous at
$T_*$. Defining $W(t):=V(t)$ for $t\ge T_*$ extends $W$ to a global
continuous solution, since both sides of \eqref{eq:Vh-rewriten-standard} pass continuously to
the limit at $T_*$ and $V$ solves \eqref{eq:Vh-rewriten-standard} beyond $T_*$. This
contradicts maximality of the interval. Hence, we obtain that $T_*=\infty$ and
$W(t)=V(t)$ for all $t\ge0$ holds.
\hfill$\Box$

\begin{remark}
When $A=0$, the solution of \dref{eq:Vh-rewriten-standard} is not unique.
Let $h(x)=|x|^p$ with $p\in(0,1)$.
Let $K=(\frac{\Gamma(\beta+1-\alpha)}{\Gamma(\beta+1)})^{1/(1-p)}$ and
$V(t)=-K\max\{t-a,0\}^\beta$ with $\beta=\alpha/(1-p)$ where $a$ is any nonnegative constant.
It is straightforward to verify that $V(t)$ solves 
 $V(t)=-\frac{1}{\Gamma(\alpha)}\int_0^t(t-s)^{\alpha-1}|V(s)|^p\,\mathrm ds$.
\end{remark}

We state the regularity result under the class-$\mathcal K$ assumption.
\begin{theorem}
\label{thm:regularity}
Let $h$ be a class-$\mathcal K$ function, $A>0$, and $0<\alpha<1$.
Let $V$ be the unique global positive solution of \eqref{eq:Vh-rewriten-standard} ensured
by Theorem \ref{thm:global-classK} and Theorem \ref{thm:all-solutions}.
Let $g=h(V)$.
Then:
\begin{enumerate}
\item[\rm(i)] for every $T>0$, $V\in C^{0,\alpha}[0,T]$ and
$
|V(t_2)-V(t_1)|\le \frac{2h(A+1)}{\Gamma(\alpha+1)}\,|t_2-t_1|^\alpha$
for $0\leq t_1<t_2\leq T$,
with the sharp expansion at the origin
\begin{equation}\label{eq:origin}
V(t)=A-\frac{h(A)}{\Gamma(\alpha+1)}\,t^\alpha+o(t^\alpha),
\end{equation}
so the H\"older exponent $\alpha$ is optimal and $V'(0^+)=-\infty$;
\item[\rm(ii)] the Caputo derivative exists and is continuous on
$[0,\infty)$:
${}_0^C D_t^\alpha V(t)=-g(t)=-h(V(t))$;
\item[\rm(iii)] the rescaled function $W(s):=V(s^{1/\alpha})$
satisfies
$
W'(0)=-\frac{h(A)}{\Gamma(\alpha+1)}.
$
\end{enumerate}
\end{theorem}
\noindent{\bf Proof.}
{\bf Step 1:} The global $\alpha$-H\"older continuity follows from  Theorem \ref{thm:global-classK}.
We show that the H\"older exponent $\alpha$ is optimal.
For the expansion at the origin, substituting $s=tr$ gives
$
\frac{(\Ial g)(t)}{t^\alpha}
=\frac1{\Gamma(\alpha)}\int_0^1(1-r)^{\alpha-1}g(tr)\,\mathrm dr.
$
Since $g$ is continuous at $0$, it is seen that $g(tr)\to g(0)=h(A)$ uniformly for
$r\in[0,1]$. This, together with the integrability of $(1-r)^{\alpha-1}$ on $[0,1]$, yields
$$
\frac{(\Ial g)(t)}{t^\alpha}\longrightarrow
\frac{g(0)}{\Gamma(\alpha)}\int_0^1(1-r)^{\alpha-1}\,\mathrm dr,
=\frac{g(0)}{\Gamma(\alpha+1)} \mbox{ as }t\to0^+,
$$
which gives \eqref{eq:origin}. Since $g(0)=h(A)>0$, for any $\beta>\alpha$,
\[
\frac{V(t)-A}{t^\beta}
=-\frac{h(A)}{\Gamma(\alpha+1)}\,t^{\alpha-\beta}(1+o(1))
\longrightarrow-\infty,
\]
so no better H\"older exponent than $\alpha$ is possible at $0$, and the
classical derivative $V'(0^+)$ does not exist finitely.

\noindent{\bf Step 2:} We show the existence and continuity of Caputo derivative for $V(t)$.
By the semigroup property ${}_0I^{1-\alpha}_t\Ial={}_0I_t^1$ for
continuous functions, we obtain
$$
{}_0I^{1-\alpha}_t\big(V(t)-V(0)\big)
=-{}_0I^{1-\alpha}_t\Ial g(t)=-{}_0I^1_t g(t)
=-\int_0^t g(s)\,\mathrm ds.
$$
The right-hand side is $C^1$ with derivative $-g(t)$, so
${}_0^C D_t^\alpha V(t)=-g(t)=-h(V(t))$ for $t>0$, and the Caputo
derivative extends continuously to $t=0$ with limit
${}_0^C D_0^\alpha V(0)=-h(A)$.

\noindent{\bf Step 3:} We show the differentiability of $W$ at the origin.
Set $t=s^{1/\alpha}$ and substitute $u=tr$ in the fractional integral to give
\[
\frac{W(s)-W(0)}{s}
=-\frac{1}{s\,\Gamma(\alpha)}\int_0^t(t-u)^{\alpha-1}g(u)\,\mathrm du
=-\frac1{\Gamma(\alpha)}\int_0^1(1-r)^{\alpha-1}g(rs^{1/\alpha})\,\mathrm dr.
\]
Uniform continuity of $g$ at $0$ and integrability of 
$(1-r)^{\alpha-1}$ give
$
\lim_{s\to0^+}\frac{W(s)-W(0)}{s}
=-\frac{g(0)}{\Gamma(\alpha)}\int_0^1(1-r)^{\alpha-1}\,\mathrm dr
=-\frac{g(0)}{\Gamma(\alpha+1)}
=-\frac{h(A)}{\Gamma(\alpha+1)},
$
which yields $W'(0)=-\frac{h(A)}{\Gamma(\alpha+1)}$.
\hfill$\Box$

\begin{theorem}
\label{thm:regularity-diffable}
Let $h$ be a class-$\mathcal K$ function, $A>0$, and $0<\alpha<1$.
Let $V$ be the unique global positive solution of \eqref{eq:Vh-rewriten-standard} ensured
by Theorem \ref{thm:global-classK} and Theorem \ref{thm:all-solutions}.
Assume that $h$ is locally H\"older on $(0,\infty)$ with exponent
$\beta>1-\alpha$.
Then
$V\in C^1(0,\infty)$ and
\begin{equation}\label{eq:Vprime-reg}
V'(t)=-\frac1{\Gamma(\alpha)}
\left[h(V(t))\,t^{\alpha-1}
+(\alpha-1)\int_0^t(t-s)^{\alpha-2}
\big(h(V(s))-h(V(t))\big)\,\mathrm ds\right],\quad t>0.
\end{equation}
Moreover, $W(t)=V(t^{1/\alpha})\in C^1[0,\infty)$, with
$W'(0)=-h(A)/\Gamma(\alpha+1)$ and $W'\leq0$.
\end{theorem}
\noindent{\bf Proof.}
From Step 6 in the proof of Theorem \ref{Thm-Dhinfatzero-mono-asympt},
in order to show \dref{eq:Vprime-reg},
it suffices to prove $g=h(V)\in C^{0,\beta'}_{\rm loc}(0,+\infty)$ with some $\beta'>1-\alpha$.
Fix $0<\varepsilon<T<\infty$.
By Theorem \ref{thm:global-classK} and Theorem \ref{thm:all-solutions},
the solution $V(t)\in[\delta_T,A]\subset(0,\infty)$
with $\delta_T=V(T)>0$ for $t\in[\varepsilon,T]$, on which $h$ is
H\"older of exponent $\beta$ with constant $C_T$.
Since $V$ is non-increasing and $h$ is non-decreasing, $g=h(V)$ is
non-increasing. Applying the same monotonicity comparison as in
\eqref{eq:V-t-s-alpha} to the Volterra equation $V=A-{}_0I_t^\alpha g$ gives,
for $0\le s<t$,
\begin{equation}\label{eq:V-t-s-alp}
0\le V(s)-V(t)\le \frac{h(A)}{\Gamma(\alpha+1)}(t^\alpha-s^\alpha).
\end{equation}
On $[\varepsilon,T]$, the mean value theorem and
$r^{\alpha-1}\le\varepsilon^{\alpha-1}$ for $r\in[\varepsilon,T]$
give $t^\alpha-s^\alpha=\alpha\int_s^t r^{\alpha-1}\,\mathrm dr
\le\alpha\varepsilon^{\alpha-1}(t-s)$. This, together with \eqref{eq:V-t-s-alp}, implies that
 $V$ is Lipschitz:
$$
|V(t)-V(s)|\le L_T|t-s|,\ \
L_T:=\frac{h(A)\varepsilon^{\alpha-1}}{\Gamma(\alpha)}, \ \ \forall t,s\in[\varepsilon,T]
$$
Combining this with the $\beta$-H\"older map $h$ on $[\delta_T,A]$, we obtain
$$
|g(t)-g(s)|\le C_T|V(t)-V(s)|^\beta\le C_T L_T^\beta|t-s|^\beta,
$$
so $g\in C^{0,\beta}[\varepsilon,T]$. Since $\beta>1-\alpha$, the singular
integral \dref{eq:Vprime-reg} is convergent on $[\varepsilon,T]$. The
arbitrariness of $\varepsilon,T$ and Lemma \ref{lem-weakly-singular-cont}
yield $g\in C^{0,\beta}_{\rm loc}(0,+\infty)$ with exponent $\beta>1-\alpha$,
and consequently $V\in C^1(0,+\infty)$. From $W(s)=V(s^{1/\alpha})$, we have
$W\in C^1(0,\infty)$.

Next we show $\lim_{t\to0^+}W'(t)=W'(0)$. Set
$
I(t):=\int_0^t(t-s)^{\alpha-2}\big(g(s)-g(t)\big)\,\mathrm ds.
$
From \eqref{eq:Vprime-reg},  it is seen that
\[
W'(t^\alpha)=-\frac{g(t)}{\Gamma(\alpha+1)}
+\frac{1-\alpha}{\alpha\Gamma(\alpha)}\,t^{1-\alpha}I(t).
\]
Since by Theorem \ref{thm:regularity},  $W'_+(0)=-g(0)/\Gamma(\alpha+1)$ and $g(t)\to g(0)$,
it suffices to prove
$\lim_{t\to0^+}t^{1-\alpha}I(t)=0$.
Since $h$ is locally $\beta$-H\"older on $(0,+\infty)$ and $V(t)\in[V(1),A]$ for $t\in[0,1]$,  there is a
constant $C_1>0$ such that for $0\leq s<t\leq1$, $|g(s)-g(t)|=|h(V(s))-h(V(t))|
\le C_1\,|V(s)-V(t)|^\beta
\le C_2(t^\alpha-s^\alpha)^\beta$ with $C_2=C_1[h(A)/\Gamma(\alpha+1)]^\beta$,
where \dref{eq:V-t-s-alp} is used.
Substituting this into $I(t)$
and writing $s=tr$, we obtain
$$
|I(t)|
\le C_2\int_0^t(t-s)^{\alpha-2}(t^\alpha-s^\alpha)^\beta\,\mathrm ds
= C_2t^{\alpha-1}\,t^{\alpha\beta}
  \int_0^1(1-r)^{\alpha-2}(1-r^\alpha)^\beta\,\mathrm dr.
$$
The last integral is finite: near $r=0$ it is bounded, near $r=1$ we use
$1-r^\alpha\sim\alpha(1-r)$, so the integrand behaves like
$(1-r)^{\alpha-2+\beta}$, which is integrable precisely because of
$\alpha-2+\beta>-1$.
This implies that
$
t^{1-\alpha}|I(t)|
\le C_2K_\alpha\,t^{\alpha\beta}\to0
$ as $t\to0^+$,
where $K_\alpha=\int_0^1(1-r)^{\alpha-2}(1-r^\alpha)^\beta\mathrm dr<+\infty$.
Therefore $W'(t)\to -g(0)/\Gamma(\alpha+1)=W'_+(0)$, so $W'$ is
continuous at $0$ and $W\in C^1[0,+\infty)$.
Finally, since, by Theorem \ref{thm:global-classK}, $V$ is decreasing on $[0,\infty)$ and $V\in C^1(0,\infty)$,
we have $V'(t)\le0$ for $t>0$, and $W'(t)=V'(t^{1/\alpha})\frac{1}{\alpha}\,t^{\frac{1}{\alpha}-1}$ for $t>0$
 together with $W'_+(0)\le0$ shows $W'(t)\le0$ for all $t\ge0$
\hfill$\Box$

\begin{remark}
It should be pointed out that the fractional cusp of $V(t)$ at the origin is unavoidable.
The failure of $V'(0)$ to exist finitely is not caused by a lack of regularity of $h$,
it is the universal feature of a Caputo equation with non-zero initial forcing.
Even for $h(x)=x$, where the solution is $V(t)=A E_\alpha(-t^\alpha)$, one has
$V(t)=A-\frac{A}{\Gamma(\alpha+1)}t^\alpha+O(t^{2\alpha})$. The rescaling
$W(s)=V(s^{1/\alpha})$ removes precisely this fractional cusp, which is
why the natural regularity statement for $V$ is $C^{0,\alpha}$ at the origin
but $C^1$ in the interior, and for $W$ it is $C^1[0,\infty)$.
\end{remark}

\begin{remark}
The comparison principle in Lemma~\ref{lem:comparison} uses only the
monotonicity of $h$ and continuity of the solutions, no Lipschitz or
differentiability condition on $h$ is involved. The local Lipschitz
assumption on $h'$ imposed in
Theorem \ref{Thm-Dhinfatzero-mono-asympt}
is therefore not needed for existence, uniqueness, or decay of the positive
continuous solution. It is used for the regularity statements, in
particular the differentiability of the rescaled function
$W(t)=V(t^{1/\alpha})$.
\end{remark}

Next, we show that the decay-rate estimate for the solution of \eqref{eq:Vh-rewriten-standard}
with \emph{no assumption whatsoever on $h'$}. 
The decay argument only uses the monotonicity of $g$, which is inherited from the
monotonicity of $h$ and $V$, together with the limit $V(t)\to0$.

\begin{theorem}
\label{thm:decay-classK}
Let $h$ be a class-$\mathcal K$ function, $A>0$, and $0<\alpha<1$.
Let $V$ be the unique global positive solution of \eqref{eq:Vh-rewriten-standard} ensured
by Theorem \ref{thm:global-classK} and Theorem \ref{thm:all-solutions}.
Then, as $t\to\infty$,
\begin{enumerate}
\item[\rm(i)]
\begin{equation}\label{eq:decay-g}
h(V(t))=g(t)\sim \frac{A}{\Gamma(1-\alpha)}\,t^{-\alpha};
\end{equation}
\item[\rm(ii)] if $h(x)\sim\kappa x^p$ as $x\to0^+$ for some $p>0$ and
$\kappa>0$, then
\begin{equation}\label{eq:decay-power}
V(t)\sim
\left(\frac{A}{\kappa\,\Gamma(1-\alpha)}\right)^{1/p}t^{-\alpha/p};
\end{equation}
\item[\rm(iii)] if there are constants $0<c_1\le c_2$ and exponents
$0<p_1\le p_2$ and $\delta>0$ such that
$c_1 x^{p_2}\le h(x)\le c_2 x^{p_1}$ for $0\le x\le\delta$,
then for all sufficiently large $t$,
\begin{equation}\label{eq:decay-bounds}
\left(\frac{A}{c_2\Gamma(1-\alpha)}\right)^{1/p_1}(1+o(1))\,t^{-\alpha/p_1}
\le V(t)\le
\left(\frac{A}{c_1\Gamma(1-\alpha)}\right)^{1/p_2}(1+o(1))\,t^{-\alpha/p_2};
\end{equation}
in particular, $h(x)\asymp x^p$ on $[0,\delta]$ implies
$V(t)\asymp t^{-\alpha/p}$.
\end{enumerate}
\end{theorem}
\noindent{\bf Proof.}
\noindent{\bf Step 1:} We show (i) holds.
Since $h$ is non-decreasing and $V$ is non-increasing, the composition
$g=h(V)$ is non-increasing. Moreover, $0\le g(t)\le h(A)$, so that $g$ is
bounded and locally integrable. By continuity of $h$ and $h(0)=0$, $\lim_{t\to\infty}g(t)=h(0)=0$.
From \eqref{eq:Vh-rewriten-standard}, we obtain $(\Ial g)(t)=A-V(t)\to A$ as $t\to\infty$.
Let $\widehat g(s):=\int_0^\infty e^{-st}g(t)\mathrm{d}t$ be the Laplace transform of $g$.
For $s>0$, taking the Laplace transform of \eqref{eq:Vh-rewriten-standard} and using
$\mathcal L(\Ial g)(s)=s^{-\alpha}\widehat g(s)$ gives
\begin{equation}\label{eq:laplace}
\widehat V(s)=\frac{A}{s}-s^{-\alpha}\widehat g(s),
\qquad\text{hence}\qquad
\widehat g(s)=s^{\alpha-1}\big(A-s\widehat V(s)\big).
\end{equation}
Since  $V(t)\in(0,A]$ and $\lim_{t\to\infty}V(t)=0$, it follows from
final value theorem that  $\lim_{s\to0^+}s\widehat V(s)=\lim_{t\to\infty}V(t)=0$.
Inserting this into \eqref{eq:laplace}, we derive
\begin{equation}\label{eq:g-laplace-asympt}
\widehat g(s)\sim A\,s^{\alpha-1}=A\,s^{-(1-\alpha)}, \ s\to0^+.
\end{equation}
Note that only the boundedness and the limit $V(t)\to0$ were used in \dref{eq:g-laplace-asympt}, no
differentiability of $h$ was involved.
Now we use  Hardy-Littlewood Tauberian theorem
for monotone functions \cite[Page 446, Theorem 4]{Feller1971Book} which states that
if $g\ge0$ is non-increasing and locally integrable on $[0,\infty)$ and, for some $\rho\in(0,1)$ and $C>0$,
$\widehat g(s)\sim C\,s^{-\rho}$ as $s\to0$, then $g(t)\sim\frac{C}{\Gamma(\rho)}\,t^{\rho-1}$ as $t\to\infty$.
Clearly, $g=h(V)$ satisfies the hypotheses of Hardy-Littlewood Tauberian theorem.
Taking $\rho=1-\alpha\in(0,1)$ and $C=A$ in \eqref{eq:g-laplace-asympt}, we obtain
\[
g(t)\sim\frac{A}{\Gamma(1-\alpha)}\,t^{(1-\alpha)-1}
=\frac{A}{\Gamma(1-\alpha)}\,t^{-\alpha},
\]
which is precisely \eqref{eq:decay-g}, proving part~(i).

\noindent{\bf Step 2:} We show (ii) holds.
Assume $h(x)\sim\kappa x^p$ as $x\to0^+$. Since $V(t)\to0$ and $g=h(V)$, we have
$
\frac{g(t)}{\kappa V(t)^p}=\frac{h(V(t))}{\kappa V(t)^p}\to1$ as $t\to+\infty,
$
so that $\kappa V(t)^p=g(t)\,(1+o(1))$. Using part~(i), we obtain
\[
V(t)=\left(\frac{g(t)}{\kappa}\right)^{1/p}(1+o(1))
\sim\left(\frac{A}{\kappa\Gamma(1-\alpha)}\right)^{1/p}t^{-\alpha/p},\mbox{ as }t\to+\infty,
\]
which proves \eqref{eq:decay-power}.

\noindent{\bf Step 3:} We show (iii) holds.
Let $t\ge T_\delta$ be so large that $V(t)\le\delta$.
From $c_1 x^{p_2}\le h(x)\le c_2 x^{p_1}$ for $0\le x\le\delta$
and $g(t)=h(V(t))$, we obtain
$
c_1 V(t)^{p_2}\le g(t)\le c_2 V(t)^{p_1}.
$
The right-hand inequality and part~(i) give
\[
V(t)\ge\left(\frac{g(t)}{c_2}\right)^{1/p_1}
\sim\left(\frac{A}{c_2\Gamma(1-\alpha)}\right)^{1/p_1}t^{-\alpha/p_1},
\]
while the left-hand inequality gives
\[
V(t)\le\left(\frac{g(t)}{c_1}\right)^{1/p_2}
\sim\left(\frac{A}{c_1\Gamma(1-\alpha)}\right)^{1/p_2}t^{-\alpha/p_2}.
\]
Together these yield \eqref{eq:decay-bounds}. When $p_1=p_2=p$, the two
bounds coincide in order and $V(t)\asymp t^{-\alpha/p}$; with
$c_1=c_2=\kappa$ they reproduce the exact asymptotic \eqref{eq:decay-power}.
\hfill$\Box$

\begin{remark}
For an arbitrary class-$\mathcal K$ function $h$, the inverse formula
$V(t)\sim h^{-1}\!\big(A\,\Gamma(1-\alpha)^{-1}t^{-\alpha}\big)$ is not a
genuine asymptotic equivalence unless $h$ satisfies an additional
regularity assumption on $h$.
The sharp universal statements under the pure class-$\mathcal K$ assumption are
\eqref{eq:decay-g} for $h(V)$ and the two-sided bounds \eqref{eq:decay-bounds}
for $V$. If $h$ is regularly varying at $0$ with index $p>0$, i.e.
$h(x)=x^p L(x)$ for a slowly varying function $L$ near $0$, then Karamata's
inversion theorem applies: $h^{-1}$ is regularly varying at $0$ with index
$1/p$, and
\[
V(t)\sim h^{-1}\!\left(\frac{A}{\Gamma(1-\alpha)}t^{-\alpha}\right),
\]
which contains part~(ii) as the constant case $L\equiv\kappa$.
\end{remark}

\begin{remark}
Part~(i) uses the  Hardy-Littlewood Tauberian theorem, but the upper bound in
\eqref{eq:decay-bounds} admits a direct proof that uses only monotonicity.
For $t\ge T_\delta$, $V$ non-increasing gives $V(s)\ge V(t)$ for
$0\le s\le t$, hence by $c_1 x^{p_2}\le h(x)$ 
 for $0\le x\le\delta$,
$
h(V(s))\ge h(V(t))\ge c_1 (V(t))^{p_2}.
$
Equation \eqref{eq:Vh-rewriten-standard} and $V(t)\ge0$ therefore give
$
A\ge A-V(t)=\frac1{\Gamma(\alpha)}\int_0^t(t-s)^{\alpha-1}h(V(s))\,\mathrm ds
\ge \frac{c_1 V(t)^{p_2}}{\Gamma(\alpha+1)}t^\alpha,
$
so that
\begin{equation}\label{eq:elementary-upper}
V(t)\le
\left(\frac{A\,\Gamma(\alpha+1)}{c_1}\right)^{1/p_2}t^{-\alpha/p_2}
\qquad(t\ge T_\delta).
\end{equation}
The matching lower bound relies on the Tauberian estimate \eqref{eq:decay-g}.
\end{remark}

\begin{theorem}\label{thm:CapD-Vine-p-stability}
Let $\alpha\in(0,1)$, $p\in(0,1)$, and $c>0$, and let $V\geq0$ and
$V\in C[0,\infty)$ satisfy
\begin{equation}\label{CapDVine-p}
{}^C_0D_t^{\alpha}V(t)\le -c\,V^{p}(t).
\end{equation}
Then $V$ is Mittag-Leffler stable. In general, \eqref{CapDVine-p}
neither entails nor precludes finite-time stability of $V$. Both
possibilities are consistent with the inequality.
\end{theorem}
{\bf Proof.} We first show that $V(t)$ is Mittag-Leffler stable.
Since $V(t)\equiv0$ if $V(0)=0$, we consider $V(0)>0$.
By ${}^C_0D_t^{\alpha}V(t)\leq-cV^{p}(t)\leq0$ and Step 1 in the proof of Theorem \ref{thm-Gen-fl},
we obtain $V(t)\leq V(0)$ for all $t\geq0$. This, together with $p\in(0,1)$, gives $V^{p-1}(t)\geq V^{p-1}(0)$.
By ${}^C_0D_t^{\alpha}V(t)\leq-cV^{p}(t)\leq0$ again, we have
$$
{}^C_0D_t^{\alpha}V(t)\leq-cV^{p}(t)
=-cV^{p-1}(t)V(t)\leq-cV^{p-1}(0)V(t).
$$
This, jointly with Lemma \ref{lem-fracsys-compare} and
the fact that ${}^C_0D_t^{\alpha}W(t)=-cV^{p-1}(0)W(t)$ with $W(0)=2V(0)$
has a solution $W(t)=2V(0)E_{\alpha}(-cV^{p-1}(0)t^{\alpha})$, gives
$
V(t)\leq 2V(0)E_{\alpha}(-cV^{p-1}(0)t^{\alpha}).
$

From Theorem \ref{Thm-Dhinfatzero-mono-asympt} or from Theorem \ref{thm:global-classK} and Theorem \ref{thm:all-solutions},
it is seen that ${}^C_0D_t^{\alpha}V(t)=-c\,V^{p}(t)$ with $V(0)>0$ has a unique solution
$V(t)\in C[0,+\infty)$ satisfying $V(t)\in(0,V(0)]$
and thus that \eqref{CapDVine-p} does not entail finite-time stability of $V$.

Next, let $A>0$ and take $V(t)=(A-C_1t^{\alpha})_{+}=\max\{A-C_1t^{\alpha},0\}$
for a sufficiently large constant $C_1>0$
satisfying $-C_1\Gamma(1+\alpha)+cA^p \leq0$.
Let $t_0=(A/C_1)^{1/\alpha}$. Clearly, $V(t)>0$ for $t\in[0,t_0)$ and $V(t)\equiv0$ for $t\geq t_0$.
Moreover, it is trivial to see that $V$ is locally absolutely continuous on the half line $[0,+\infty)$.
Hence, for $t\in[0,t_0]$, we have
\begin{eqnarray*}
{}^C_0D_t^{\alpha}V(t)&=&\frac1{\Gamma(1-\alpha)}
\int_0^t(t-s)^{-\alpha}V'(s)ds
=-\frac1{\Gamma(1-\alpha)}
\int_0^t(t-s)^{-\alpha}C_1\alpha s^{\alpha-1}ds\\
&=&
-\frac{C_1\alpha}{\Gamma(1-\alpha)}
\int_0^1(1-\tau)^{-\alpha} \tau^{\alpha-1}d\tau
=
-C_1\Gamma(1+\alpha).
\end{eqnarray*}
Since the maximum of $V^p(t)$ on $[0,t_0]$ equals $A^p$,
 it follows from ${}^C_0D_t^{\alpha}V(t)=-C_1\Gamma(1+\alpha)$ and $-C_1\Gamma(1+\alpha)+cA^p \leq0$
that ${}^C_0D_t^{\alpha}V(t)\leq-cV^p(t)$ for $t\in[0,t_0]$.
 On the other hand, if $t\geq t_0$ then $V(t)=V^p(t)=0$, 
 and ${}^C_0D_t^{\alpha}V(t)=-\frac{C_1\alpha}{\Gamma(1-\alpha)}
\int_0^{t_0}(t-s)^{-\alpha} s^{\alpha-1}ds\leq0=-cV^p(t)$ for $t\geq t_0$,
which means that $V(t)=(A-C_1t^{\alpha})_{+}$ satisfies inequality \eqref{CapDVine-p}.
Hence, \eqref{CapDVine-p} does not preclude finite-time stability of $V$.
\hfill$\Box$
\begin{remark}
For the equality ${}^C_0D_t^\alpha V(t)=-cV^p(t)$ with $V(0)>0$ and $0<p<1$,
by  Theorem \ref{thm:decay-classK}, the solution
is positive for all $t>0$ and decays  like
$({V(0)}/{c\Gamma(1-\alpha)})^{1/p}t^{-\alpha/p}$ 
(hence is not
finite-time stable). The inequality form, however, permits faster decay,
including solutions that vanish identically after a finite time, for
which the Caputo derivative beyond the extinction time is strictly
negative due to the memory term.
\end{remark}

\noindent{\bf Acknowledgments}
The second author is grateful to Professor Kou Chunhai for his continued
interest in and support of this work, which has sustained the second author's study,
begun in 2010, of the qualitative properties of solutions to a weakly singular Volterra
integral equation related to fractional Lyapunov asymptotic stability theory.
The second author also wishes to thank the first author for his deep analytical
insight into the qualitative behaviour of solutions of weakly singular Volterra
equations, on which the argument of Theorem 3.4 rests, and for many discussions
that have shaped the second author's understanding of these questions.
The authors would like to thank Prof. Yangquan Chen for his valuable discussions
on fractional calculus, and wish to express their deep sorrow over his sudden
passing while this paper was nearing completion.
The authors would also like to thank Prof. Yiheng Wei and Prof. Fei-Fei Du
for their helpful discussions on fractional stability and fractional inequalities,
and Prof. Jiongmin Yong and Prof. Jianguo Liu for their helpful advice,
which has led to this improved version.




\end{document}